\documentclass[journal,twosided,web]{IEEEtran}
\usepackage{cite}
\usepackage{amsmath,amssymb,amsfonts}
\usepackage{graphicx}
\usepackage{algorithm,algorithmicx}
\usepackage{algpseudocode}
\usepackage{textcomp}
\usepackage{adjustbox}
\usepackage{multicol}
\usepackage{siunitx}

\usepackage[pdfstartview=XYZ,
bookmarks=true,
colorlinks=true,
linkcolor=blue,
urlcolor=blue,
citecolor=blue,
pdftex,
bookmarks=true,
linktocpage=true, 
hyperindex=true
]{hyperref}
\usepackage{orcidlink}

\def\BibTeX{{\rm B\kern-.05em{\sc i\kern-.025em b}\kern-.08em
    T\kern-.1667em\lower.7ex\hbox{E}\kern-.125emX}}
\begin{document}
\title{Tracking EEG Thalamic and Cortical Focal Brain Activity using Standardized Kalman Filtering with Kinematics Modeling}
\author{Veikka Piispa \orcidlink{0009-0003-9732-772X}, Dilshanie Prasikala \orcidlink{0009-0005-3556-5583}, Joonas Lahtinen \orcidlink{0000-0002-9377-8713}, Alexandra Koulouri \orcidlink{0000-0001-7605-3844} and Sampsa Pursiainen \orcidlink{0000-0002-9131-9070}
\thanks{
This study was supported by the Research Council of Finland through the Center of Excellence in Inverse Modelling and Imaging 2018--2025 (353089), the  Flagship of Advanced Mathematics for Sensing, Imaging and Modelling (FAME) (359185), and VN/3137/2024 Doctoral Education Pilot for Mathematics of Sensing, Imaging and Modelling. The work of J.\ Lahtinen was supported by the Jenny and Antti Wihuri Foundation.  We thank Prof.\ Carsten H.\ Wolters (University of Münster, Münster, Germany) for the fruitful discussions and his support in mathematical modelling. We are also grateful to the German Academic Exchange Service (DAAD) and the Research Council of Finland (RCF; decisions 354976 and 367453) for supporting our travels to Münster. AK and SP were also supported by project 359198. AK was also supported by Institute for Mathematical Innovation, University of Bath, UK.} \thanks{Authors are with the Faculty of Information Technology 
and Communication Sciences, Tampere University, 33720 Finland. Corresponding author: Dilshanie Prasikala, E-mail:  dilshanie.wanniachchikankanamge@tuni.fi}
}

\maketitle

\begin{abstract}

Kalman filtering has proven to be effective for estimating brain activity using EEG recordings. In particular, the introduced post hoc standardization step of the algorithm, inspired by the sLORETA time-invariant method, reduces the depth bias and thus allows the estimation to appear at the correct depth from the electrode surface. In the current work, we propose first and second-order kinematic evolution models, where the state-space vector includes not only the dipolar source activity but also its velocity and acceleration. Compared to our previous study, this motion model yields smoother and more physically plausible estimates of brain activity even when the measurement noise is high, for both superficial and deep sources. In addition, we introduce a tunable power parameter that enhances the computational efficiency of the algorithm. Our simulation study, which involves thalamic and cortical activity in the somatosensory region, demonstrates that accurate estimation and tracking of both superficial and deep brain activity are feasible.
\end{abstract}


\section{Introduction}

Electroencephalography (EEG) records voltage differences on the human scalp and can be used to monitor electrical activity inside the brain noninvasively. This forms the basis for the EEG source imaging problem, which aims to reconstruct the brain activity sources from the observed EEG signals. Unlike many other neuroimaging modalities, EEG offers high temporal resolution, which can be effectively utilized in dynamic brain imaging. Additionally, its portability and low cost make it a highly attractive tool for clinical diagnostics.

Reconstructing source activity from EEG data presents several challenges, primarily due to the limited number of scalp measurements and the presence of high noise levels \cite{hamalainen1993}. Recent advancements in source localization algorithms have significantly enhanced their utility in the medical domain, aiding in the diagnosis and localization of critical functional brain areas and epilepsy sources, particularly for surgical planning \cite{MERTENS2000E,LealAlberto2008_4,deGooijer2013_14,Coito2019_18,Neugebauer2022}. However, EEG source imaging is inherently a time-varying problem, and despite these advances, many current methods do not fully leverage the time-series nature of EEG data to improve source reconstructions.

Among the most promising approaches to solve the time-dynamic source imaging problem are Bayesian filtering techniques, particularly Kalman filtering \cite{Kalman1960,Kalman1961,Sims1969}, which have been successfully applied in various fields, e.g., in \cite{Alsadik2019,Codrescu2018,Vauhkonen1998,Harveya,Koulouri2022}. To address the dynamic aspect of brain activity, Kalman filtering has been found to be beneficial for estimating the time-varying distribution of brain activity inside the brain \cite{GalkaAndreas2004KalmanEEG}. However, as electrical potential measurements are taken non-invasively on the scalp, unweighted EEG estimators (including the Kalman filter) tend to localize the source activity in the vicinity of the sensors, even if the actual activity is deeper; this is referred to as depth bias. 

A recent algorithm incorporating the temporal Kalman filtering state-space model, accompanied with a standardization step equivalent to the time-invariant sLORETA \cite{PascualMarqui2002},  proposed in \cite{Lahtinen2024SKF}. This algorithm allowed the EEG temporal information to be included and reduced the depth-bias effect in the state-space estimates by applying a post hoc standardization at each time step \cite{Lahtinen2024Onbias}. 

In the current work, we advance the previously proposed modeling of the Standardized Kalman filtering (SKF) \cite{Lahtinen2024SKF} by using a first and second-order kinematic system, where the state-space vector of the quantity of interest, i.e., brain activity, also includes the velocity in the first-order case and the velocity and acceleration in the second-order case. Kinematic modeling is well-known in navigation, e.g., \ GPS tracking, where it is known as a crucial concept to enhance the accuracy and robustness of the positioning outcome, see, e.g., \cite{schwarz1989comparison,wang2023stochastic,hong2015medium}.  With the present implementation (dipolar activity, its velocity, and acceleration), we aim at smoother and more physically plausible estimates of the brain activity. In particular, our simulated study aims to demonstrate that this formulation enhances the accuracy of activity strength tracking over time and facilitates a clear distinction between the temporal emergence of co-existing cortical and deep sources, even in cases with high measurement noise. Furthermore, in the current work, we propose to employ a flexible exponent parameter $p$ for the weighting matrix in the standardization step. This allows for higher computational efficiency and stability, compared to $p=0.5$ 
which is the standard value used originally \cite{Lahtinen2024SKF}. 
Particularly, 
$p=1$ seems to admit a higher focality, i.e., a greater concentration in a smaller region, even in a low-resolution source space that would not be possible otherwise. This can significantly accelerate (i) the computations, as we can use a smaller state-space vector, and (ii) the standardization step, since we need to deal with matrix multiplication and inversion operators. 

Based on the results, the proposed algorithm constitutes a promising approach for monitoring brain activity in challenging scenarios where the time evolution of the activity is crucial, such as in epilepsy, where an epileptic focus can trigger secondary activity in anatomically or effectively connected brain regions \cite{LemieuxLouis2011CoCa}. Namely, EEG source localization is an inverse problem in which, for example, the simultaneous reconstruction of deep and superficial sources is difficult or nearly impossible \cite{RezaeiA2021}. Only recently has it been demonstrated that a deep source can be reliably reconstructed \cite{seeber2019subcortical}. Therefore, noise-robust methods such as the standardized Kalman filter (SKF) investigated here are of primary importance.


\section{Dynamic EEG Brain Imaging}

\subsection{Observation Model}\label{AA}

The relation between EEG recordings  ${\bf y}_t\in \mathbb{R}^{m}$ and the dipole field represented as the ${\bf d}_t\in \mathbb{R}^{n}$ at time frame $t=1,\cdots T$ can be expressed as
\begin{equation}\label{eq:measmodel}
    {\bf y}_t=L{\bf d}_t+{\bf r}_t,
\end{equation}
where $L\in \mathbb{R}^{m\times n}$ is the lead field matrix and ${\bf r}_t\in \mathbb{R}^{m}$ represents the measurement noise where \({\bf r}_t\sim\mathbf{N}(0,R_t)\).

\subsection{Kalman Filter Models} 
In this modeling scheme, the state variable includes higher-order temporal derivatives of the source activity. Hence, we model the evolution of the hidden state as an $s$-order kinematic system with time-invariant dynamics:
\begin{equation}
    {\bf x}_t^{(s)} = A_t^{(s)} {\bf x}_{t-1}^{(s)} + {\bf q}_t^{(s)},
\end{equation}
where 
\[
{\bf x}_t^{(1)} = 
\begin{bmatrix}
{\bf d}_t &
{\bf v}_t 
\end{bmatrix}^\mathrm{T} \in \mathbb{R}^{2n}\]
and
\[{\bf x}_t^{(2)} = 
\begin{bmatrix}
{\bf d}_t &
{\bf v}_t &
{\bf a}_t
\end{bmatrix}^\mathrm{T} \in \mathbb{R}^{3n}
\]
is the state vector consisting of source activity, velocity, and acceleration components respectively, and ${\bf q}_t\sim\mathcal{N}(0,Q_t)$.
In particular, the transition matrices for these cases are
\begin{equation}
    A_t^{(1)} = \begin{bmatrix}
        I_{n} & \Delta t \, I_{n} \\
        0_{n} & I_{n} \\
    \end{bmatrix}\in \mathbb{R}^{2n \times 2n},
\end{equation}
    and
\begin{equation}
    A_t^{(2)} = \begin{bmatrix}
        I_{n} & \Delta t \, I_{n} & \frac{1}{2} \Delta t^2 \, I_{n} \\
        0_{n} & I_{n} & \Delta t \, I_{n} \\
        0_{n} & 0_{n} & I_{n}
    \end{bmatrix} \in \mathbb{R}^{3n \times 3n},
\end{equation}

where \( \Delta t > 0 \) is the time step, and \( I_n \), \( 0_n \) denote identity and zero matrices of size \( n \times n \), respectively.

Also, the process noise covariance is
\begin{equation}
    Q_t^{(s)} = \begin{bmatrix}
        0_{sn} & 0_{s n\times n}  \\
        0_{n\times s n} & \frac{s}{\Delta t^s} \, \varphi \, I_n
    \end{bmatrix} \in \mathbb{R}^{(s+1)n \times (s+1)n},
\end{equation}

where \( \varphi \) is the tunable process noise variance affecting the acceleration component.

Finally, the observation model becomes
\begin{equation}
    {\bf y}_t = 
    \begin{bmatrix}
        L & 0_{m \times sn}
    \end{bmatrix}
    {\bf x}_t + {\bf r}_t := H^{(s)} {\bf x}_t + {\bf r}_t,
\end{equation}
where \( L \in \mathbb{R}^{m \times n} \) is the lead field matrix, \( H^{(s)} \in \mathbb{R}^{m \times (s+1)n} \) is the full observation matrix, and \( {\bf r}_t \sim \mathcal{N}(0_m, R_t) \) is the Gaussian distributed measurement noise.

\subsection{Dynamical Standardization Algorithm}

The proposed Dynamical standardized Kalman filter (DSKF) is given in Alg. \ref{Algo:KalmanFilterStandardized}. Every state-space estimate at time $t$ is given by ${\bf x}_{t\mid t}$ and  the final standardized vector is $${\bf z}_{t\mid t} = W_t {\bf x}_{t\mid t},$$ where 
 $W_t = \mathrm{Diag}\left(P_{t\mid t-1}^{-1/2} K_t S_t K_t^\mathrm{T} P_{t\mid t-1}^{-1/2} \right)^{-p} P_{t\mid t-1}^{-1/2}$ is the weighting matrix, $K_t$ is the
 Kalman gain, $S_t$ is the innovation covariance, and $P_{t\mid t-1}$ is the predictive covariance. The assumption is that the computation is started at the time point, where no significant activity appears, i.e., ${\bf x}_0=0_{3n}$. The initial source covariance is set to be $P_0=\theta I_{3n}$, where $\theta$ is set following the approach in \cite{RezaeiAtena2020PtCG}.
 
\begin{algorithm}[H]
\caption{Dynamical Standardized Kalman Filter ($s$-DSKF)}
\begin{algorithmic}[1]
\State \textbf{Input:} Initial state estimate ${\bf x}_{0\mid 0}$, covariance $P_{0\mid 0}$, system matrices $A_t^{(s)}$, $H^{(s)}$, $Q_t^{(s)}$, $R_t$, time steps $t = 1, \ldots, T$, weighting power $p>0$
\For{$t = 1$ to $T$}
    \State \textbf{Prediction Step:}
    \State ${\bf x}_{t\mid t-1} = A_t^{(s)} {\bf x}_{t-1\mid t-1}$
    \State $P_{t\mid t-1} = A_t^{(s)} P_{t-1\mid t-1} \left(A_t^{(s)}\right)^\mathrm{T} + Q_t^{(s)}$
    
    \State \textbf{Kalman Gain Calculation:}
    \State $S_t = H^{(s)} P_{t\mid t-1} \left(H^{(s)}\right)^\mathrm{T} + R_t$
    \State $K_t = P_{t\mid t-1} \left(H^{(s)}\right)^\mathrm{T} S_t^{-1}$
    
    \State \textbf{Update Step:}
    \State ${\bf x}_{t\mid t} = {\bf x}_{t\mid t-1} + K_t({\bf y}_t - H^{(s)} {\bf x}_{t\mid t-1})$
    \State $P_{t\mid t} = P_{t\mid t-1} - K_t S_t K_t^\mathrm{T}$
    
    \State \textbf{Transformation:}
    \State $W_t = \mathrm{Diag}\left(P_{t\mid t-1}^{-1/2} K_t S_t K_t^\mathrm{T} P_{t\mid t-1}^{-1/2} \right)^{-p} P_{t\mid t-1}^{-1/2}$
    \State ${\bf z}_{t\mid t} = W_t {\bf x}_{t\mid t}$

   \EndFor
\end{algorithmic}\label{Algo:KalmanFilterStandardized}
\end{algorithm}

Furthermore, here Alg. \ref{Alg:RTS} presents the Rauch-Tung-Striebel (RTS) smoother \cite{RTS1965} that is used with SKF for SSKF and DSKF algorithms.
                      
\begin{algorithm}
  \caption{Rauch-Tung-Striebel (RTS) smoother for SSKF and  DSKF}
\begin{algorithmic}
\State \textbf{Initialize Smoother:}
\State $P_{T\mid T}$, ${\bf z}_{T\mid T} = W_T {\bf x}_{T\mid T}$ (from the filter)
\State ${\bf\bar z}_{T\mid T} = {\bf z}_{T\mid T}$ and $\bar{P}_{T\mid T}=P_{T\mid T}$
\For{$t = T-1$ down to $1$}
    \State ${\bf z}_{t+1\mid t}^- = A_{t}{\bf z}_{t\mid t}$
    \State $P_{t+1\mid t}^- = A_{t}P_{t\mid t} A_{t}^\mathrm{T}+Q_t$
    \State $C_t = P_{t\mid t} A_{t}^\mathrm{T} (P_{t+1\mid t}^-)^{-1}$
    \State ${\bf \bar z}_{t\mid t} = {\bf z}_{t\mid t} + C_t({\bf \bar z}_{t+1\mid t+1} - {\bf z}_{t+1\mid t}^-)$
    \State $\bar{P}_{t\mid t} = P_{t\mid t} + C_t(\bar{P}_{t+1\mid t} - P_{t+1\mid t}^-)C_t^\mathrm{T}$
\EndFor

\end{algorithmic}\label{Alg:RTS}
\end{algorithm}

We note that the basic standardized Kalman filter (SKF) algorithm is obtained by setting $A_t=I_n$, replacing $H$ by $L$. The standardization exponent is set to be $p=1$, which can be interpreted as standardization of reconstructed power, to guarantee sufficient distinguishability of the deep activity in all situations (Appendix \ref{app:exponent}).

\begin{figure}[htb!]
\centering
\begin{minipage}{0.3cm}
\centering
    \mbox{} 
\end{minipage}
\begin{minipage}{2.2cm}
\centering
    {\bf Side view}
\end{minipage}\begin{minipage}{2.2cm}
\centering
    {\bf Back view}
\end{minipage}\begin{minipage}{2.2cm}
\centering
    {\bf Top view}
\end{minipage}\begin{minipage}{2.2cm}
\centering
    {\bf Upper back left diagonal view}
\end{minipage}
\\

\begin{minipage}{0.3cm}
\rotatebox{90}{\small{Simulated sources}}
\end{minipage}\begin{minipage}{2.2cm}
    \centering
    \includegraphics[trim={0 0 1cm 0},clip,height=2.0 cm]
    {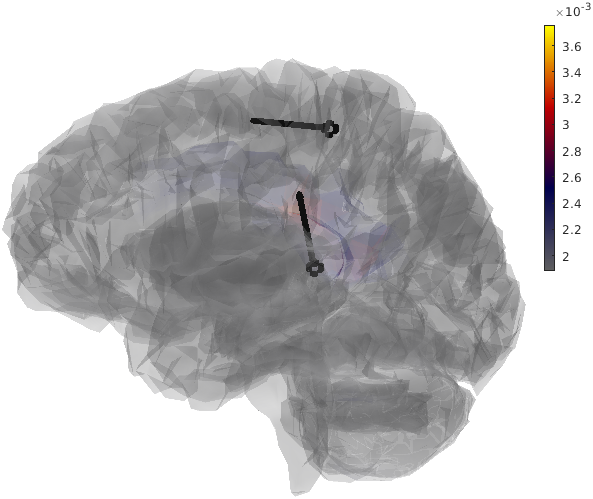}
\end{minipage}\begin{minipage}{2.2cm}
    \centering
    \includegraphics[trim={0 0 1cm 0},clip,height=2.0 cm]{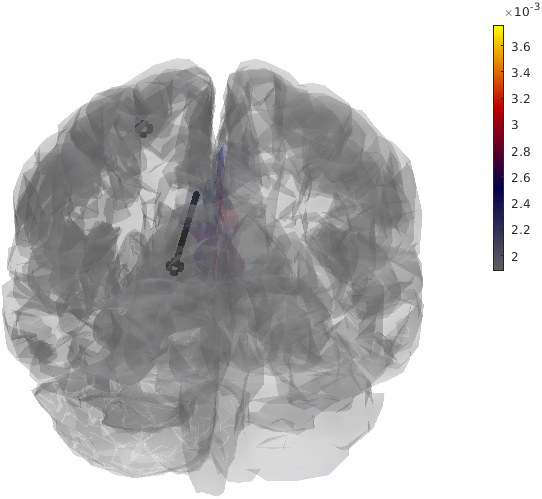}
\end{minipage}\begin{minipage}{2.0cm}
    \centering
    \includegraphics[trim={0 0 2cm 0},clip,height=1.8 cm]{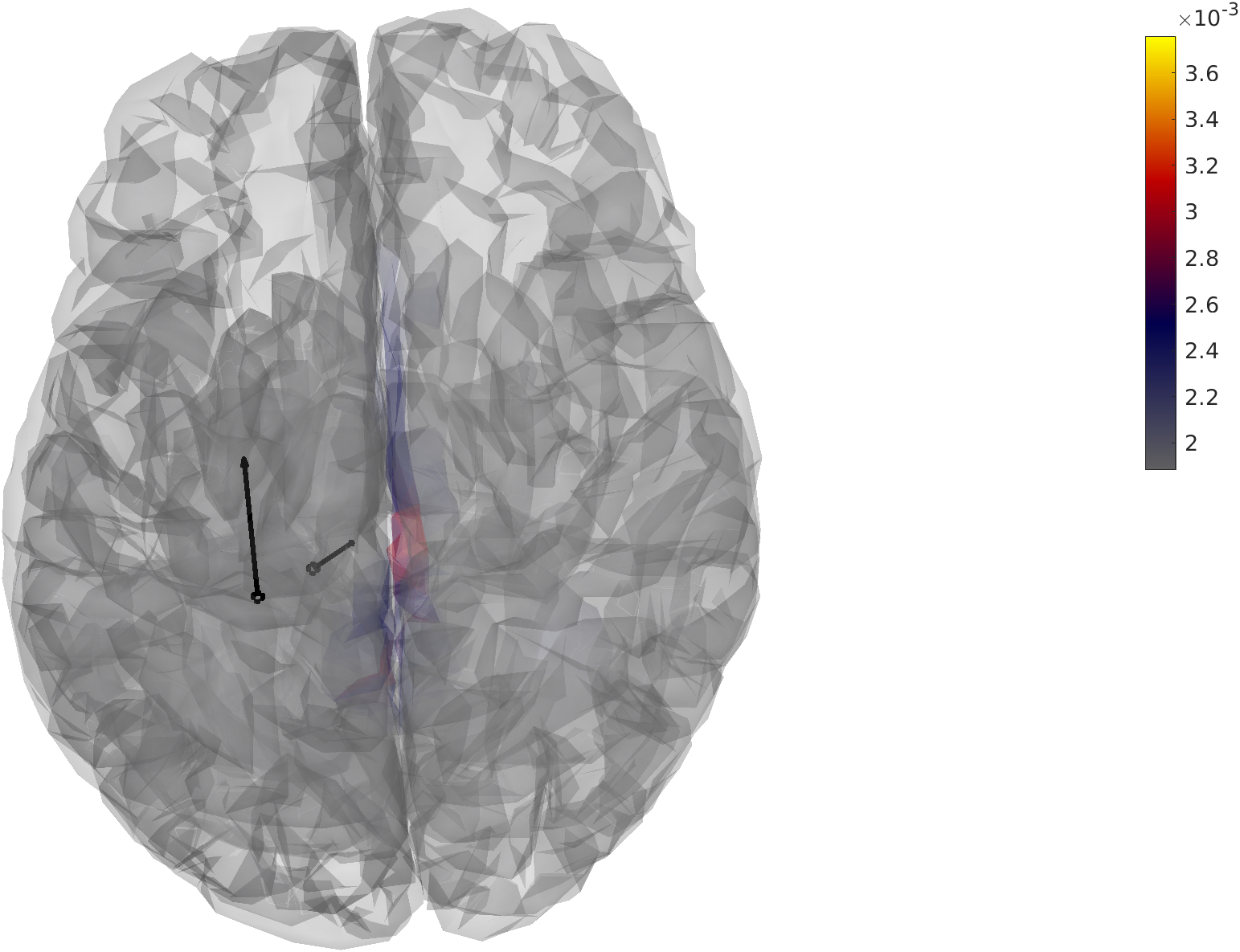}
\end{minipage}\begin{minipage}{2.2cm}
    \centering
    \includegraphics[trim={0 0 1cm 0},clip,height=2.0 cm]{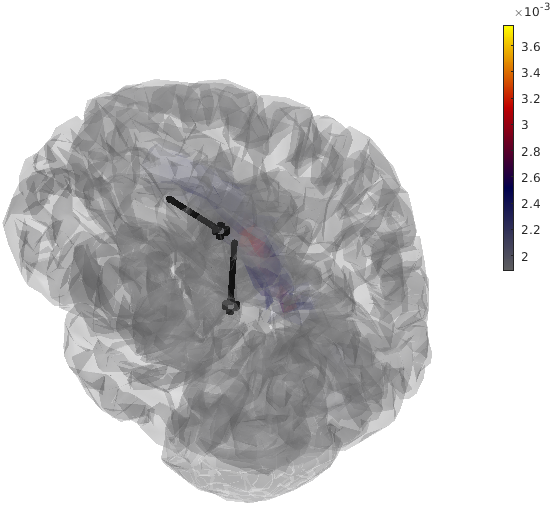}
\end{minipage}
\\
\begin{minipage}{0.3cm}
\rotatebox{90}{\small{ROI}}
\end{minipage}\begin{minipage}{2.2cm}
    \centering
    \includegraphics[trim={0 0 1cm 0},clip,height=2.0 cm]
    {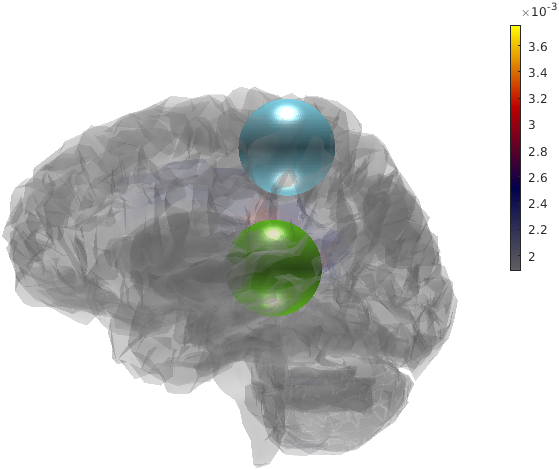}
\end{minipage}\begin{minipage}{2.2cm}
    \centering
    \includegraphics[trim={0 0 1cm 0},clip,height=2.0 cm]{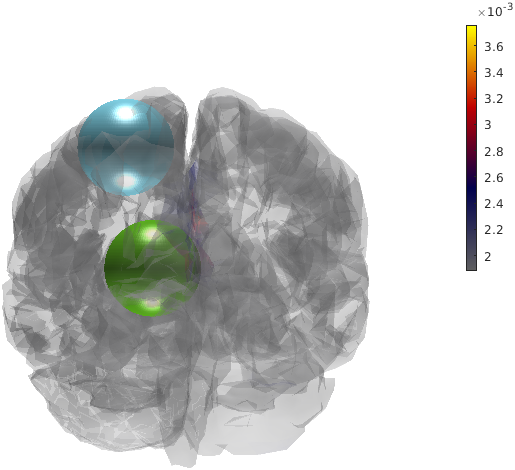}
\end{minipage}\begin{minipage}{2.2cm}
    \centering
    \includegraphics[trim={0 0 1cm 0},clip,height=2 cm]{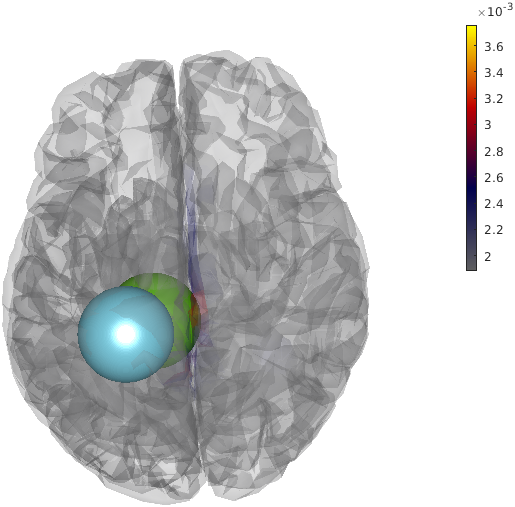}
\end{minipage}\begin{minipage}{2.2cm}
    \centering
    \includegraphics[trim={0 0 1cm 0},clip,height=2.0 cm]{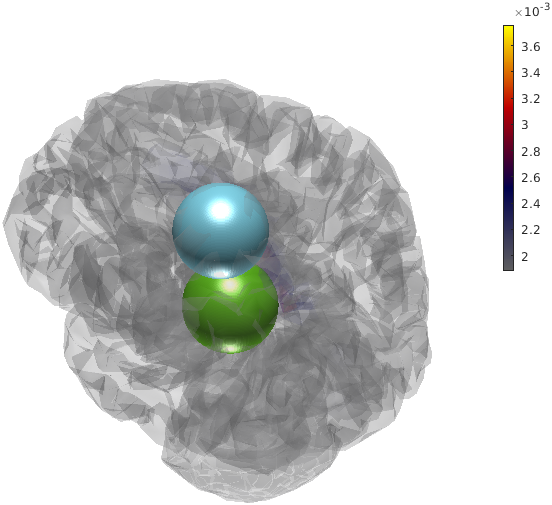}
\end{minipage}
\\
\begin{minipage}{0.3cm}
\rotatebox{90}{\small{Simulated sources}}
\end{minipage}\begin{minipage}{2.2cm}
    \centering
    \includegraphics[trim={0 0 0 0},clip,height=1.9 cm]
    {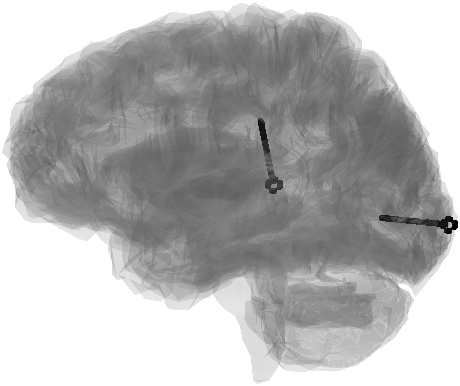}
\end{minipage}\begin{minipage}{2.2cm}
    \centering
    \includegraphics[trim={0 0 0 0},clip,height=2.0 cm]{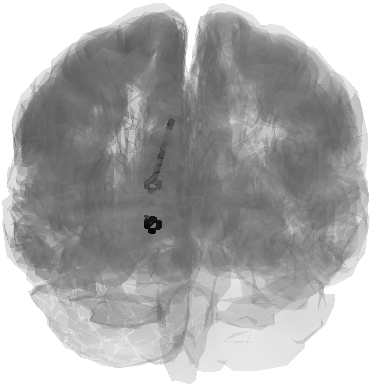}
\end{minipage}\begin{minipage}{2.2cm}
    \centering
    \includegraphics[trim={0 0 0 0},clip,height=2.0 cm]{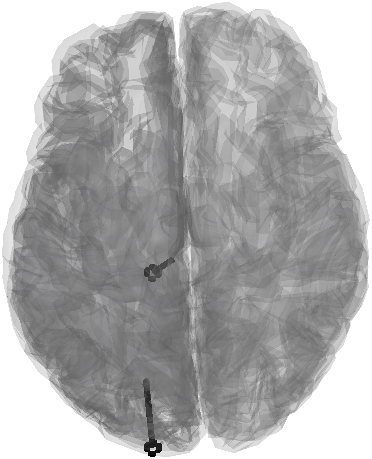}
\end{minipage}\begin{minipage}{2.2cm}
    \centering
    \includegraphics[trim={0 0 0 0},clip,height=2.0 cm]{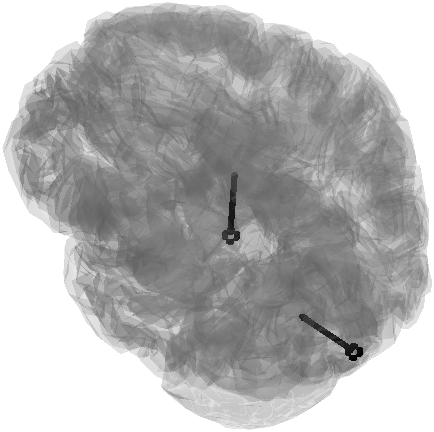}
\end{minipage}
\\
\begin{minipage}{0.3cm}
\rotatebox{90}{\small{ROI}}
\end{minipage}\begin{minipage}{2.2cm}
    \centering
    \includegraphics[trim={0 0 0 0},clip,height=2.0 cm]
    {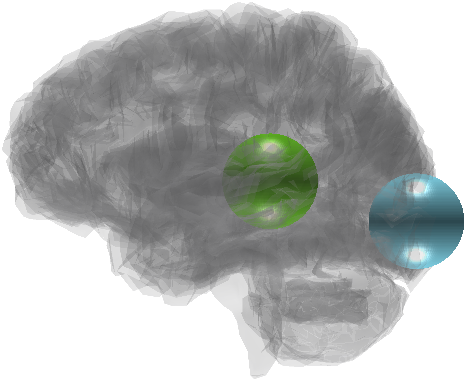}
\end{minipage}\begin{minipage}{2.2cm}
    \centering
    \includegraphics[trim={0 0 0 0},clip,height=2.0 cm]{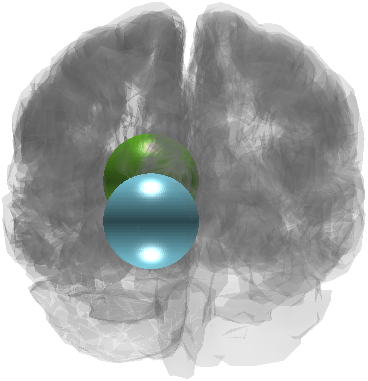}
\end{minipage}\begin{minipage}{2.2cm}
    \centering
    \includegraphics[trim={0 0 0 0},clip,height=2.0 cm]{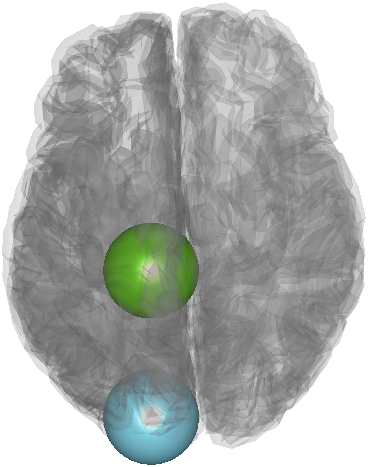}
\end{minipage}\begin{minipage}{2.2cm}
    \centering
    \includegraphics[trim={0 0 0 0},clip,height=2.0 cm]{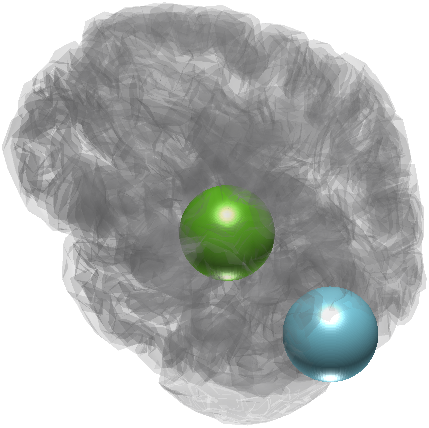}
\end{minipage}
    \caption{The first and third rows present the locations and directions of cortical (somatosensory cortex (first row), visual cortex (third row)) and deep (thalamic) dipolar sources used in data simulation. The sources are presented as black quivers within the brain model obtained from MRI data. The second and fourth rows show the regions of interest from which the activity strength is averaged for source strength tracking.}
    \label{fig:sourcesetup}
\end{figure}

\begin{figure*}
    \centering
    \begin{minipage}{9cm}
    \centering
    \includegraphics[height=6cm]{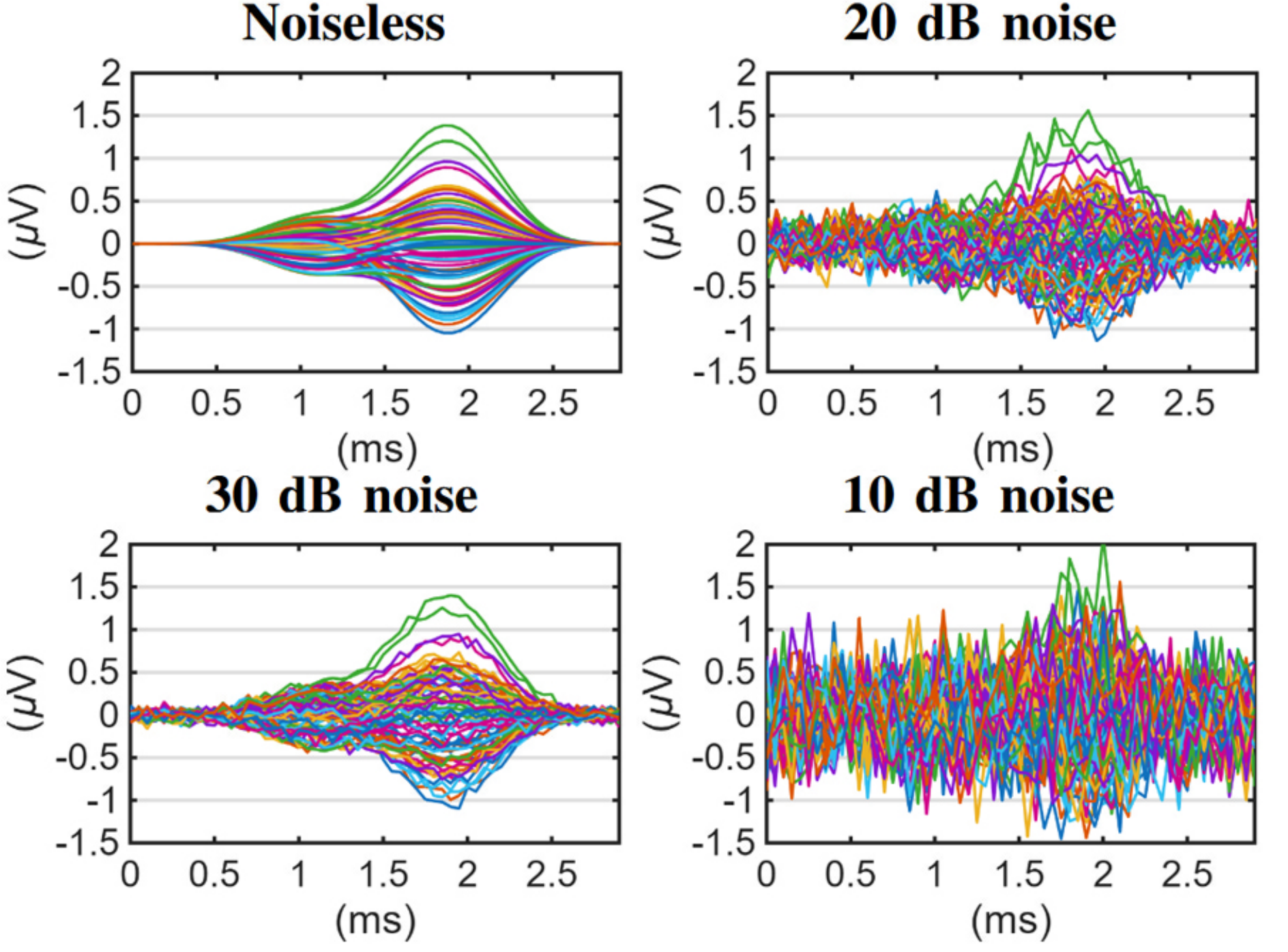}
     \end{minipage}
     \begin{minipage}{5cm}
     \centering
         \includegraphics[height=5cm]{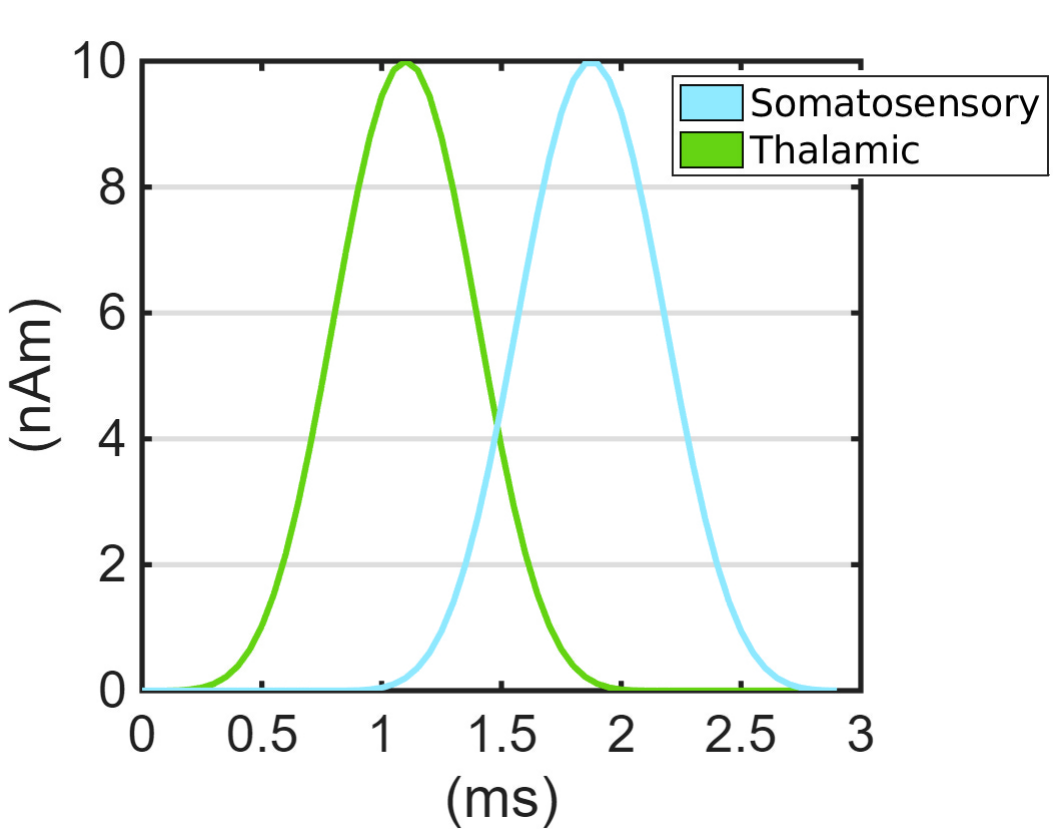}
     \end{minipage}
    \caption{{\bf Left:} Plots of the simulated EEG data without noise and with the experimented noise levels: 30, 20, and 10 dB of noise. The x-axis shows the time in milliseconds, and the y-axis shows the signal strength scaled to 1. {\bf Right:} The graph presents the time evolution of the activity strengths for the deep source placed in the thalamus (green) and the cortical source placed in the somatosensory cortex (turquoise). The x-axis shows the time in seconds, and the y-axis shows the activity strength in nanoampere-meters.}
    \label{fig:butterfly}
    \label{fig:trueevolution}
\end{figure*}


\begin{figure*}[h!]
\centering
\begin{footnotesize}
\begin{minipage}{0.5cm}
\mbox{}
\end{minipage}
\begin{minipage}{4.5cm}
\centering
30 dB
\end{minipage} 
\begin{minipage}{4.5cm}
\centering
20 dB
\end{minipage} 
\begin{minipage}{4.5cm}
\centering
10 dB
\end{minipage} \\ \vskip0.2cm
\begin{minipage}{0.5cm}
\rotatebox{90}{3-DSKF}
\end{minipage}
\begin{minipage}{4.5cm}
\centering
\includegraphics[trim={0 3cm 0 0},clip, width=4.2cm]{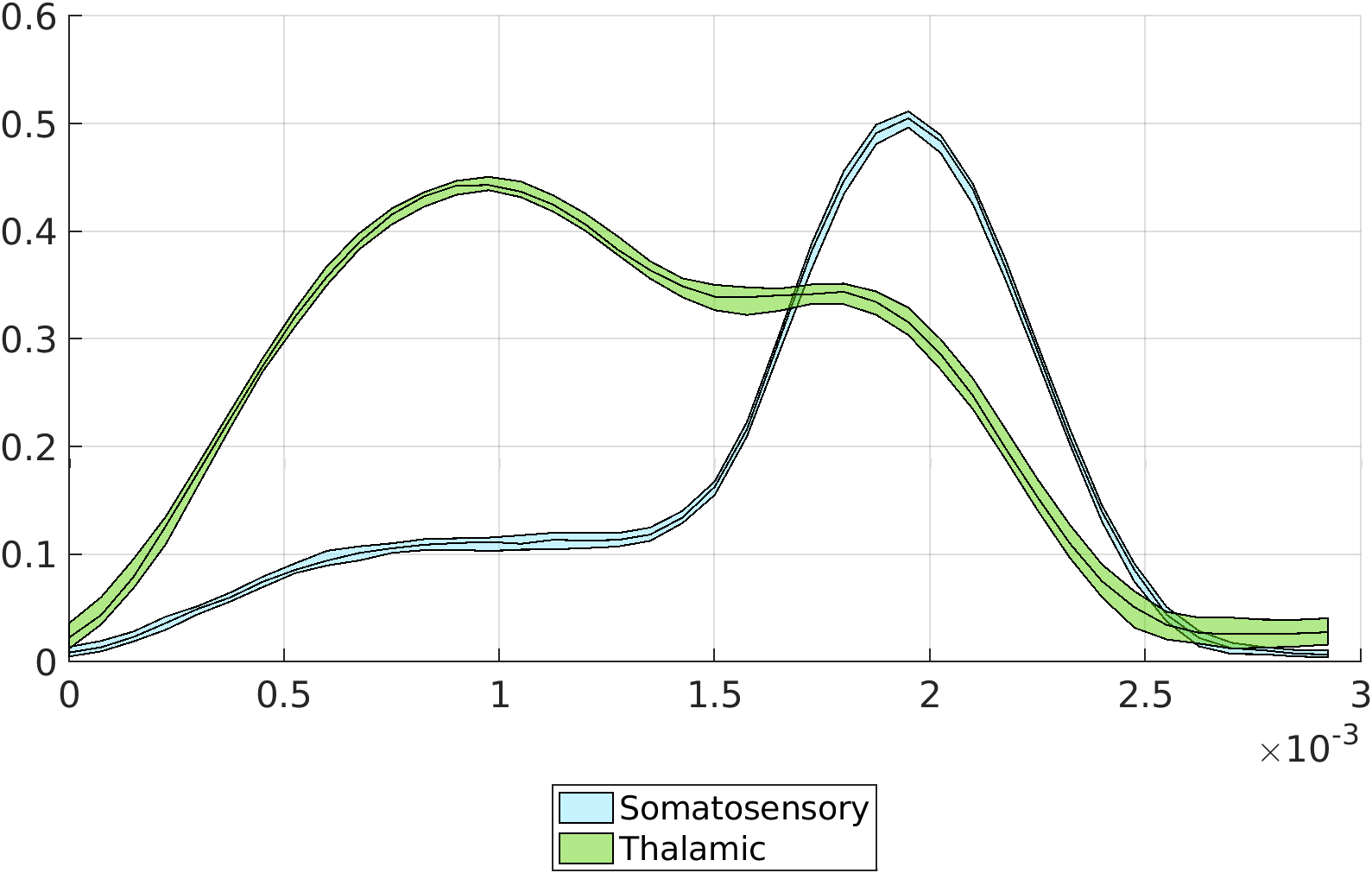}
\end{minipage}
\begin{minipage}{4.5cm}
\centering
\includegraphics[trim={0 3cm 0 0},clip, width=4.2cm]{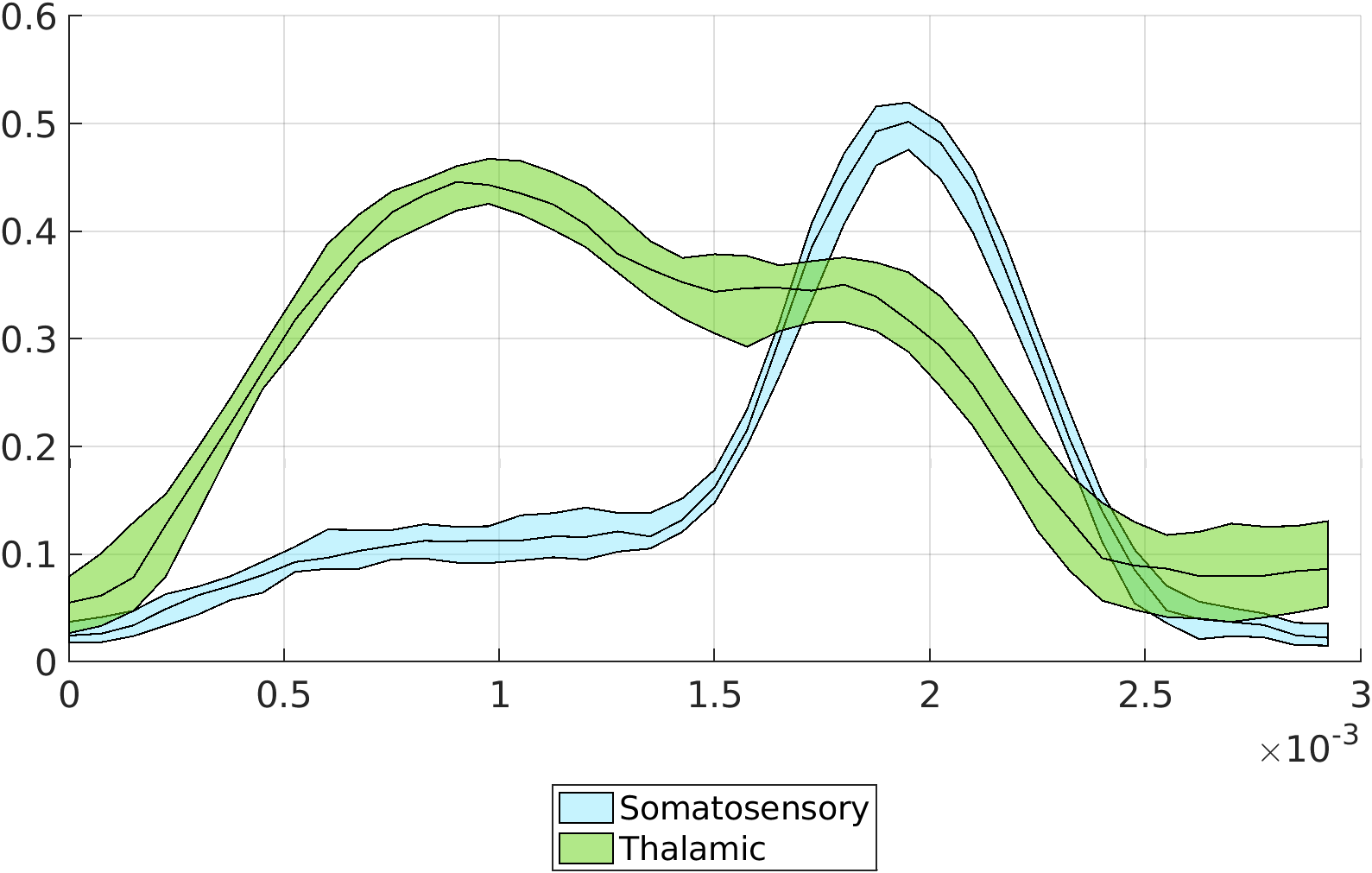}
\end{minipage}     
\begin{minipage}{4.5cm}
\centering
\includegraphics[trim={0 3cm 0 0},clip, width=4.2cm]{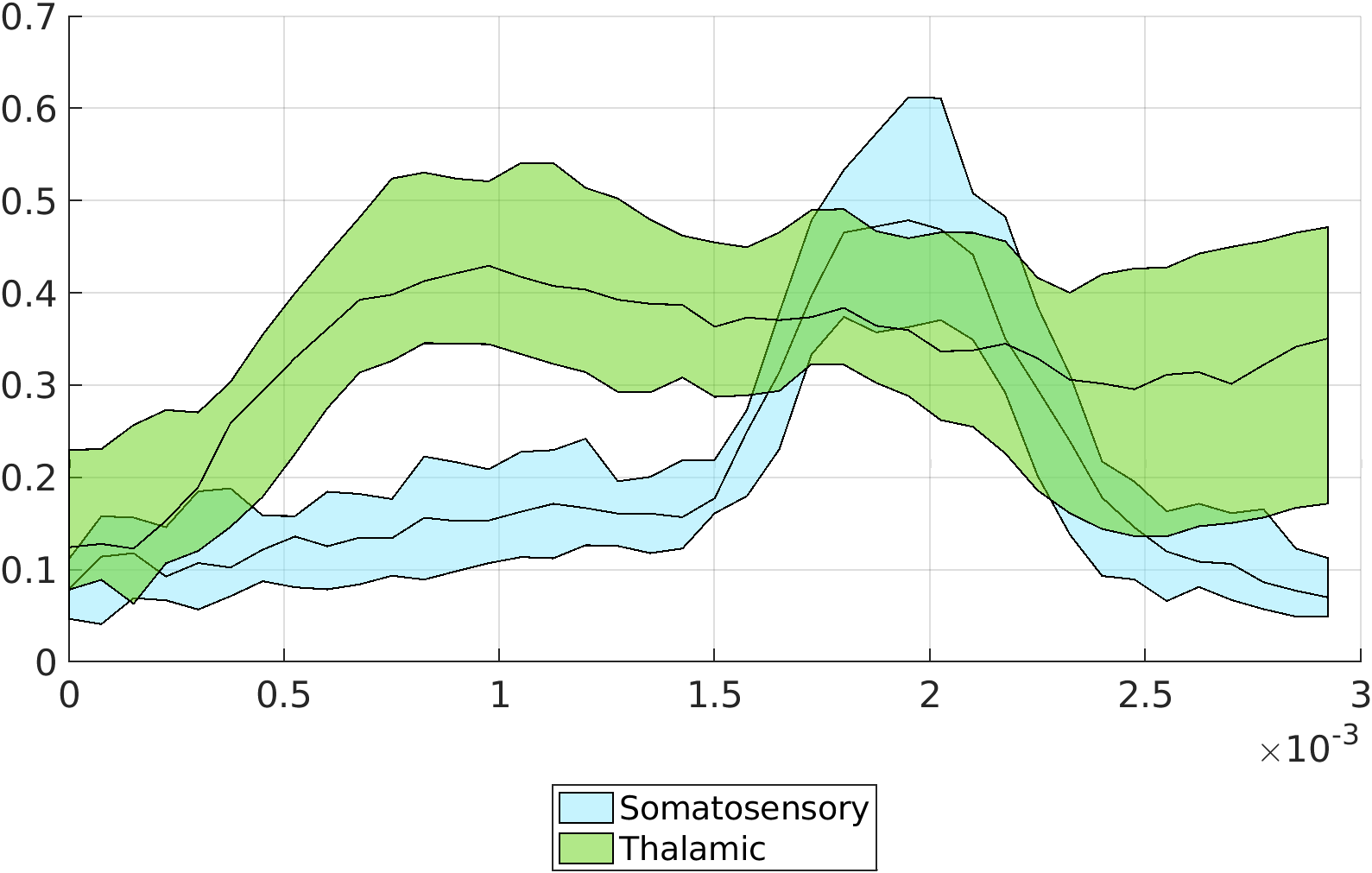}
\end{minipage}   \\ \vskip0.2cm
\begin{minipage}{0.5cm}
\rotatebox{90}{2-DSKF}
\end{minipage}
\begin{minipage}{4.5cm}
\centering
\includegraphics[trim={0 3cm 0 0},clip, width=4.2cm]{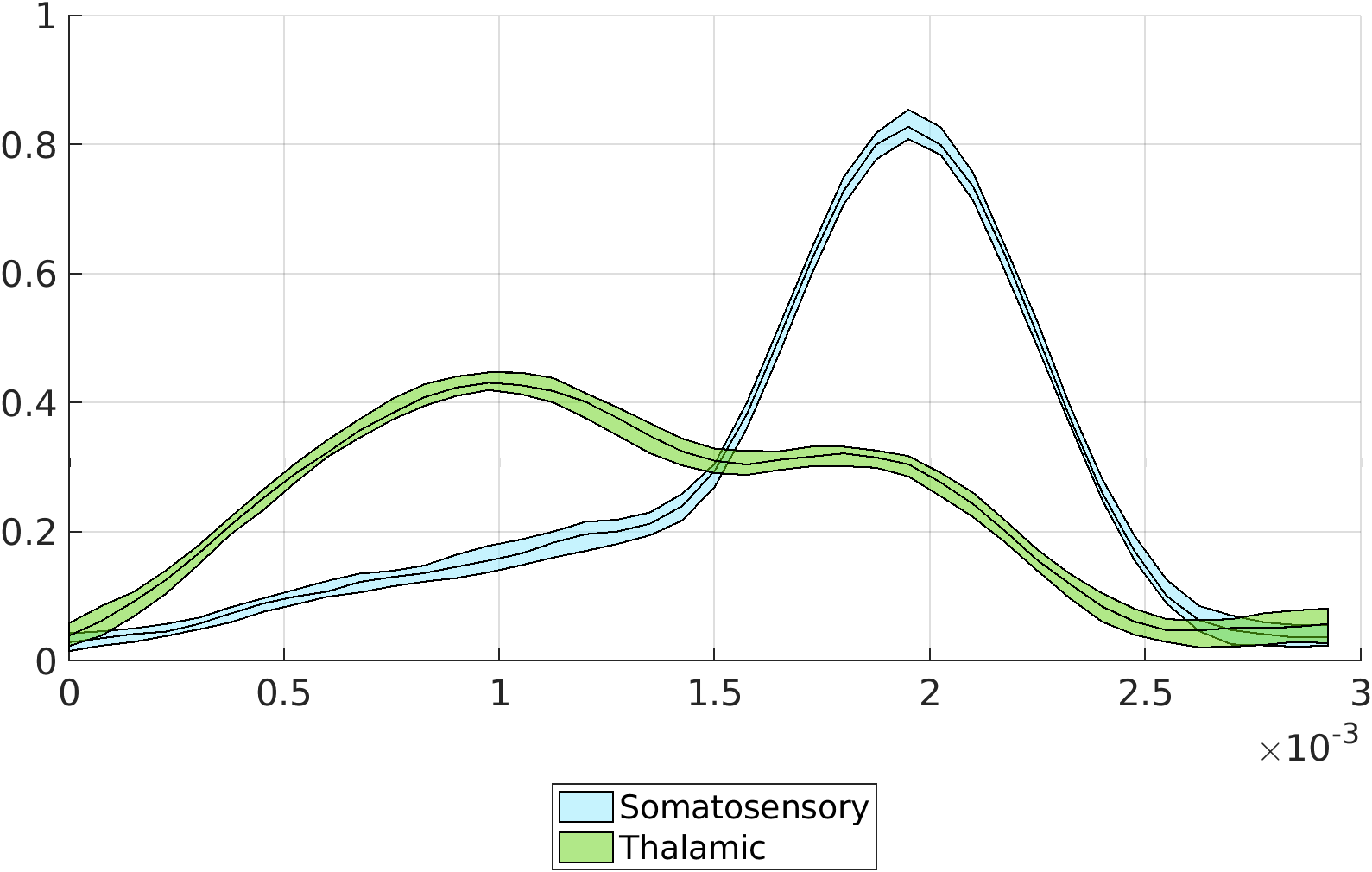}
\end{minipage}
\begin{minipage}{4.5cm}
\centering
\includegraphics[trim={0 3cm 0 0},clip, width=4.2cm]{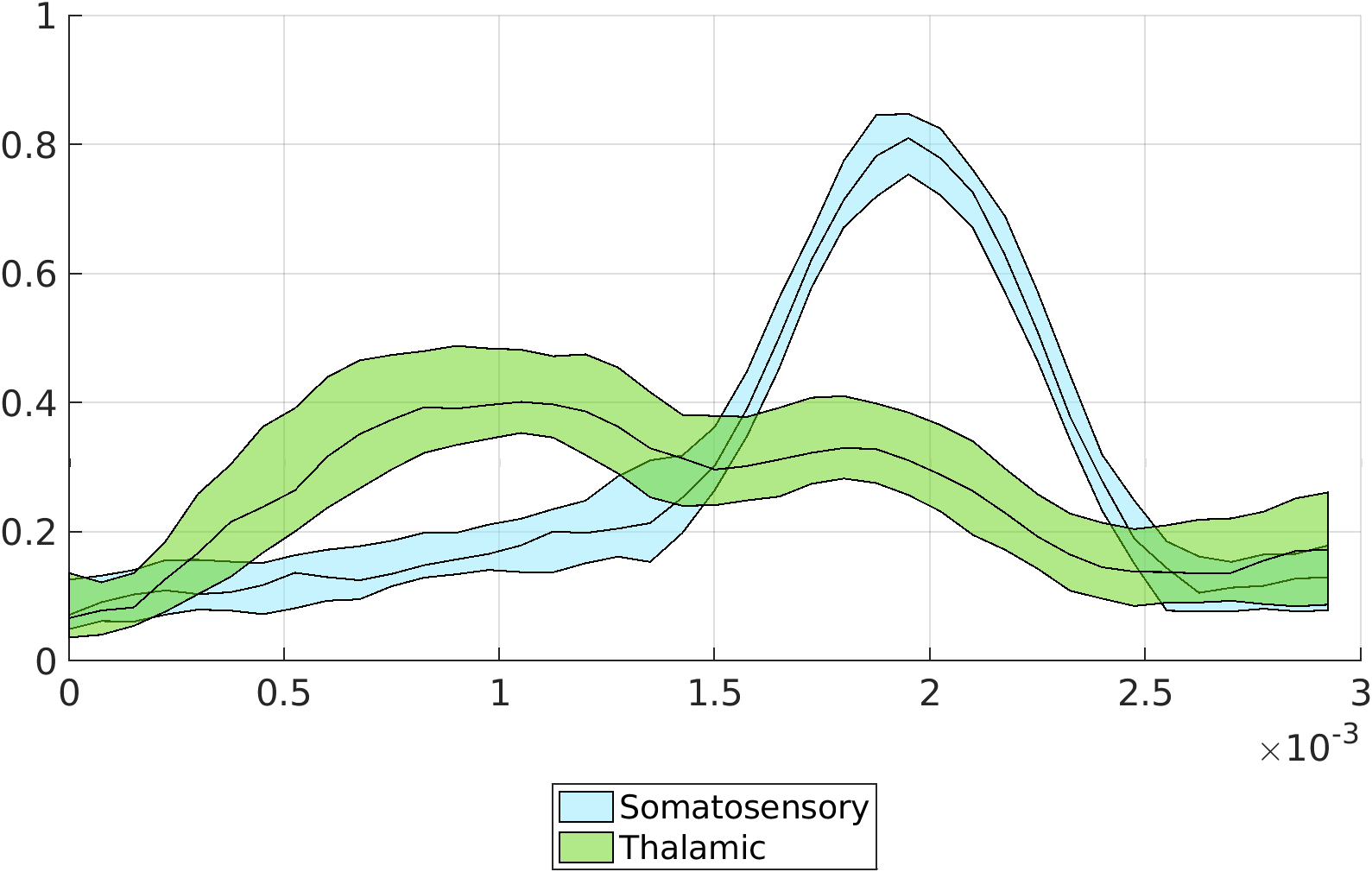}
\end{minipage}     
\begin{minipage}{4.5cm}
\centering
\includegraphics[trim={0 3cm 0 0},clip, width=4.2cm]{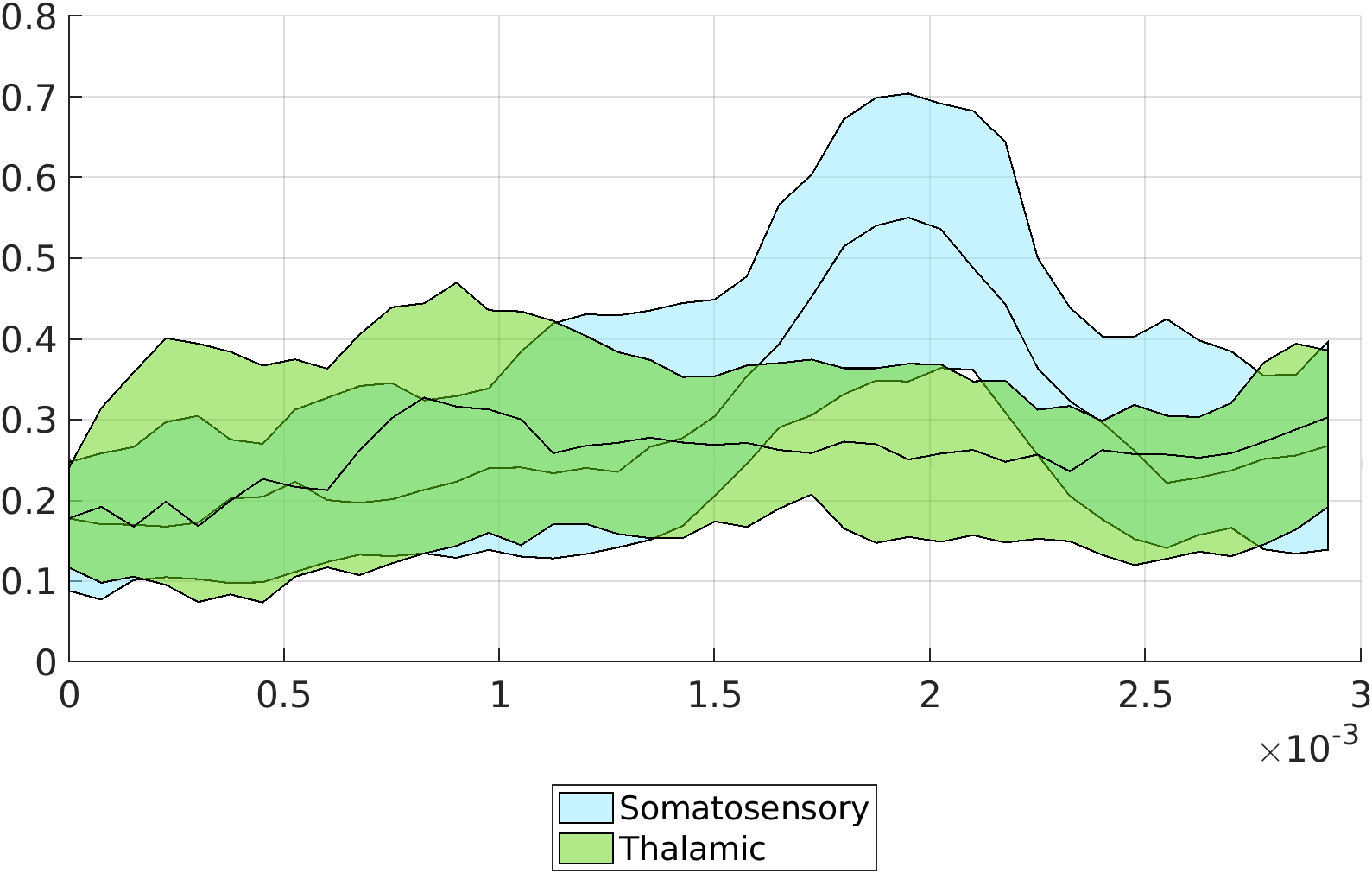}
\end{minipage}   \\ \vskip0.2cm
\begin{minipage}{0.5cm}
\rotatebox{90}{SSKF}
\end{minipage}
\begin{minipage}{4.5cm}
\centering
\includegraphics[trim={0 3cm 0 0},clip, width=4.2cm]{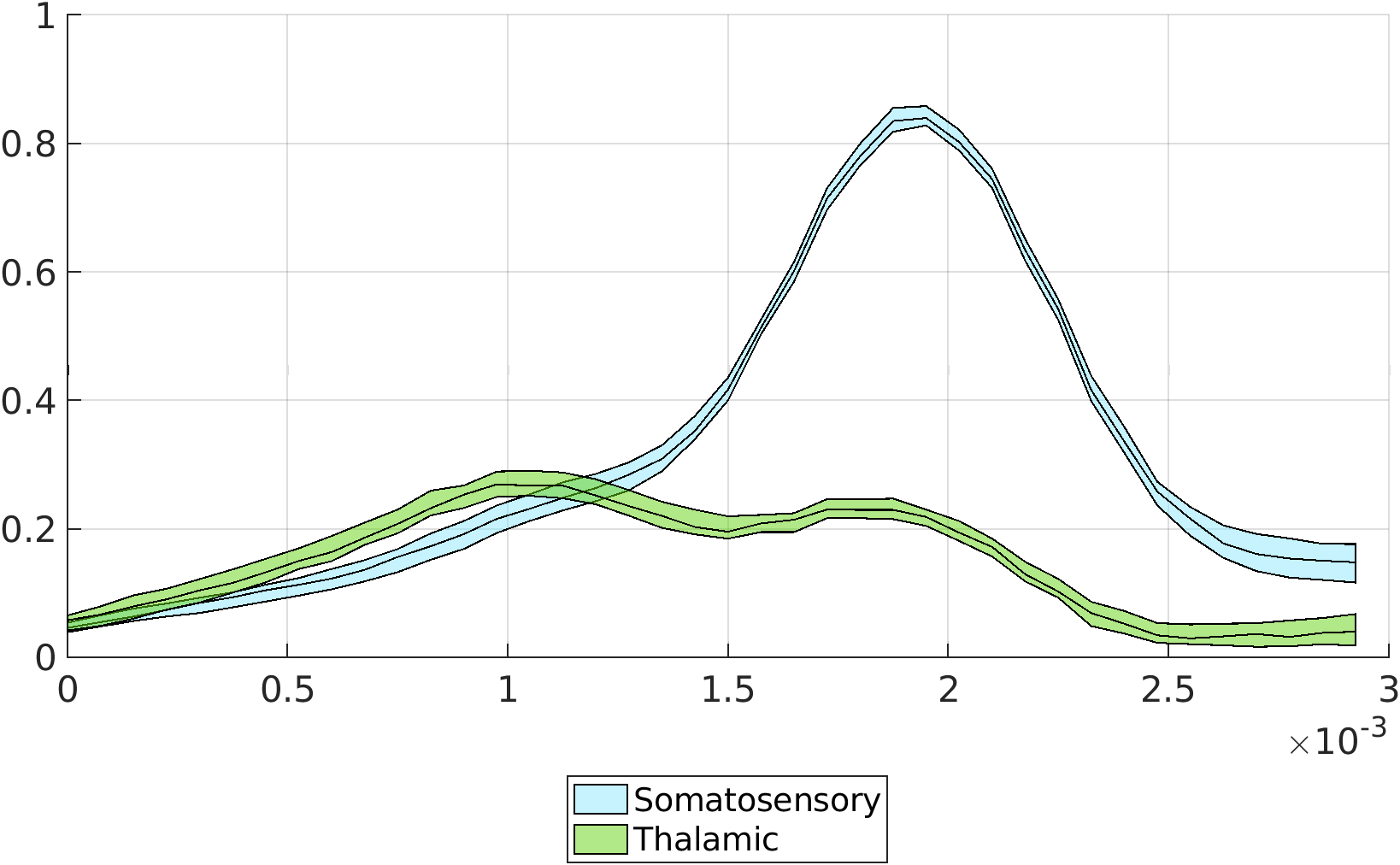}
\end{minipage}
\begin{minipage}{4.5cm}
\centering
\includegraphics[trim={0 3cm 0 0},clip, width=4.2cm]{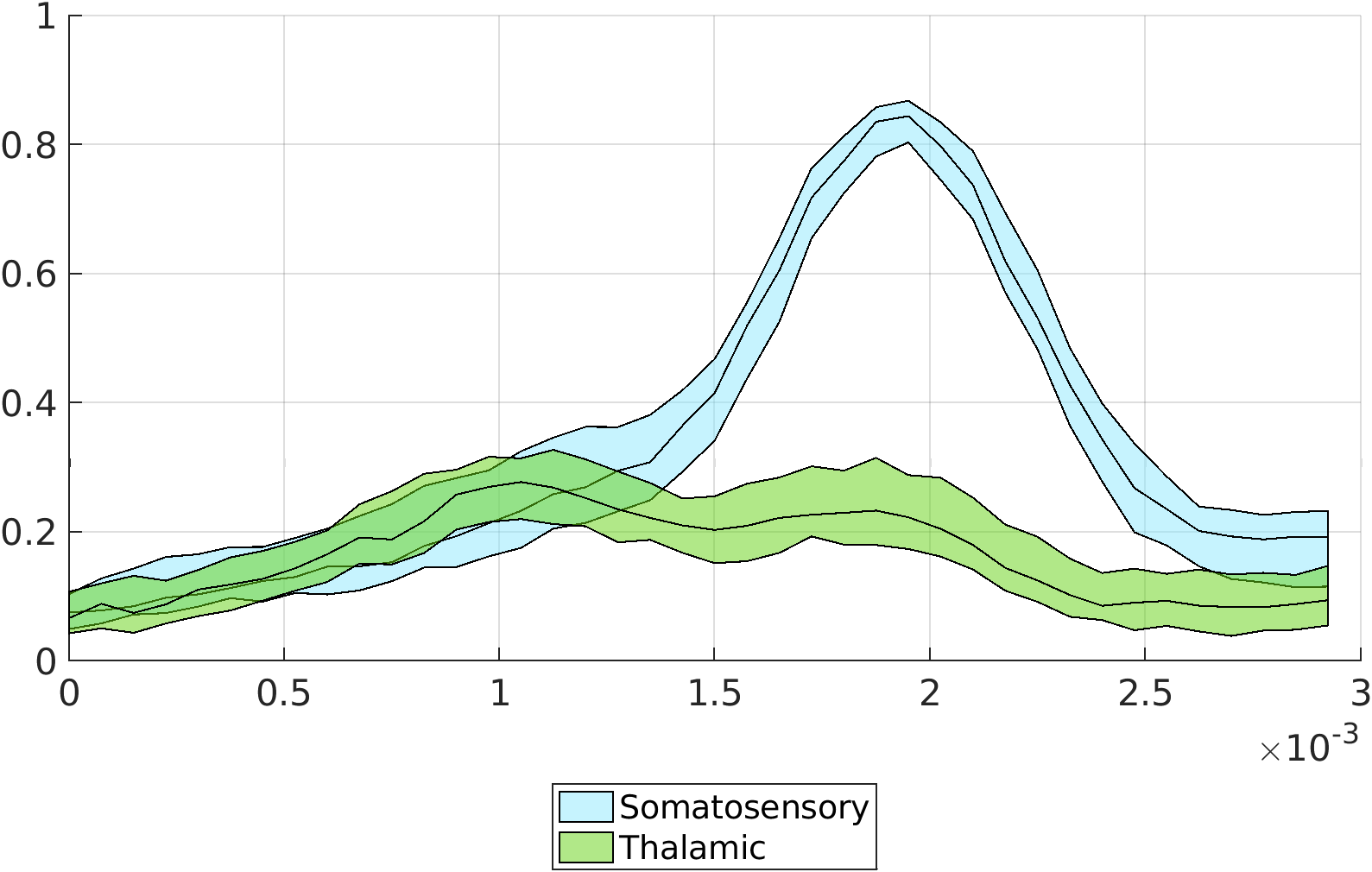}
\end{minipage}     
\begin{minipage}{4.5cm}
\centering
\includegraphics[trim={0 3cm 0 0},clip, width=4.2cm]{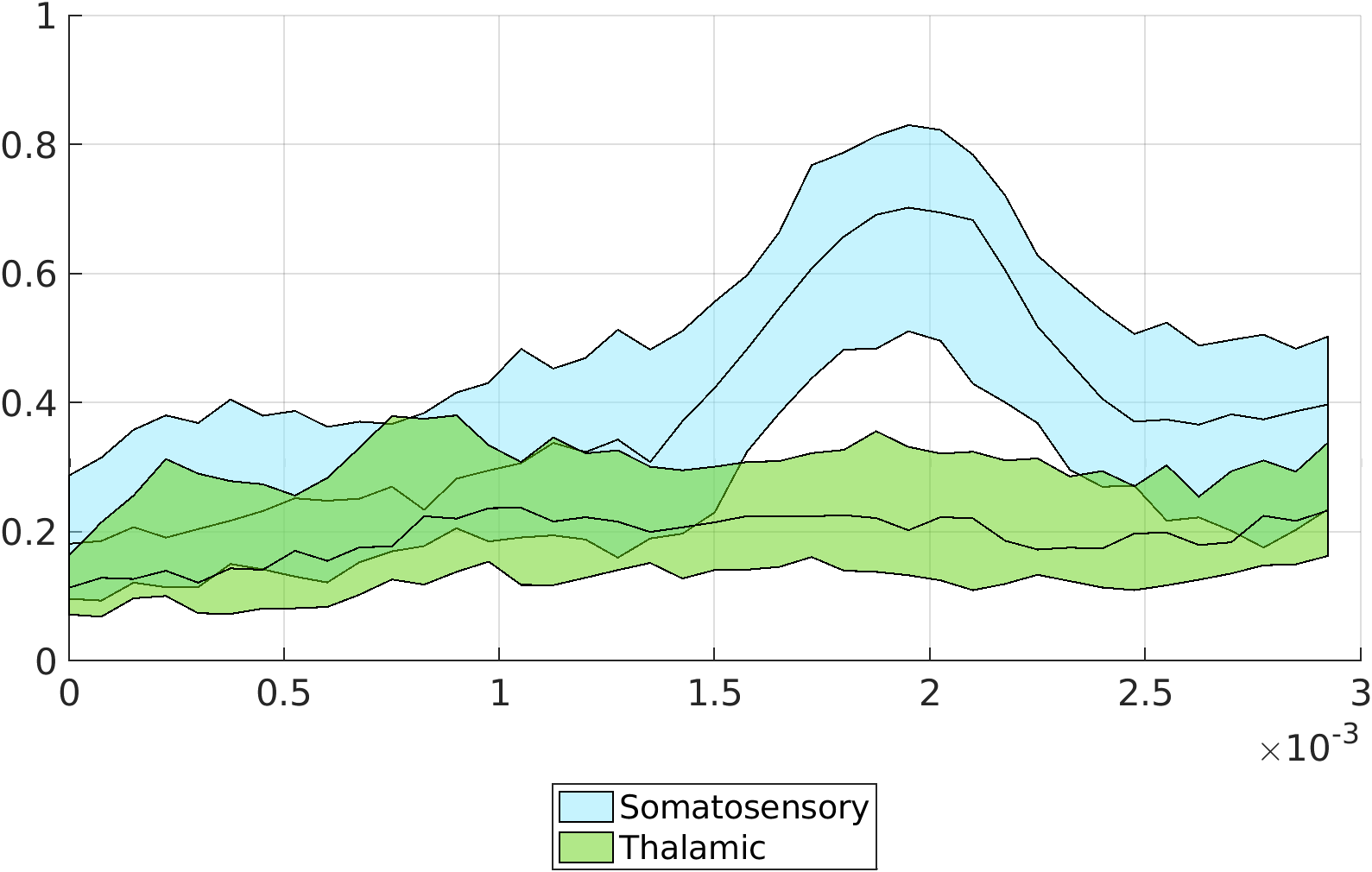}
\end{minipage}   \\ \vskip0.2cm
\begin{minipage}{0.5cm}
\rotatebox{90}{SKF}
\end{minipage}
\begin{minipage}{4.5cm}
\centering
\includegraphics[trim={0 3cm 0 0},clip, width=4.2cm]{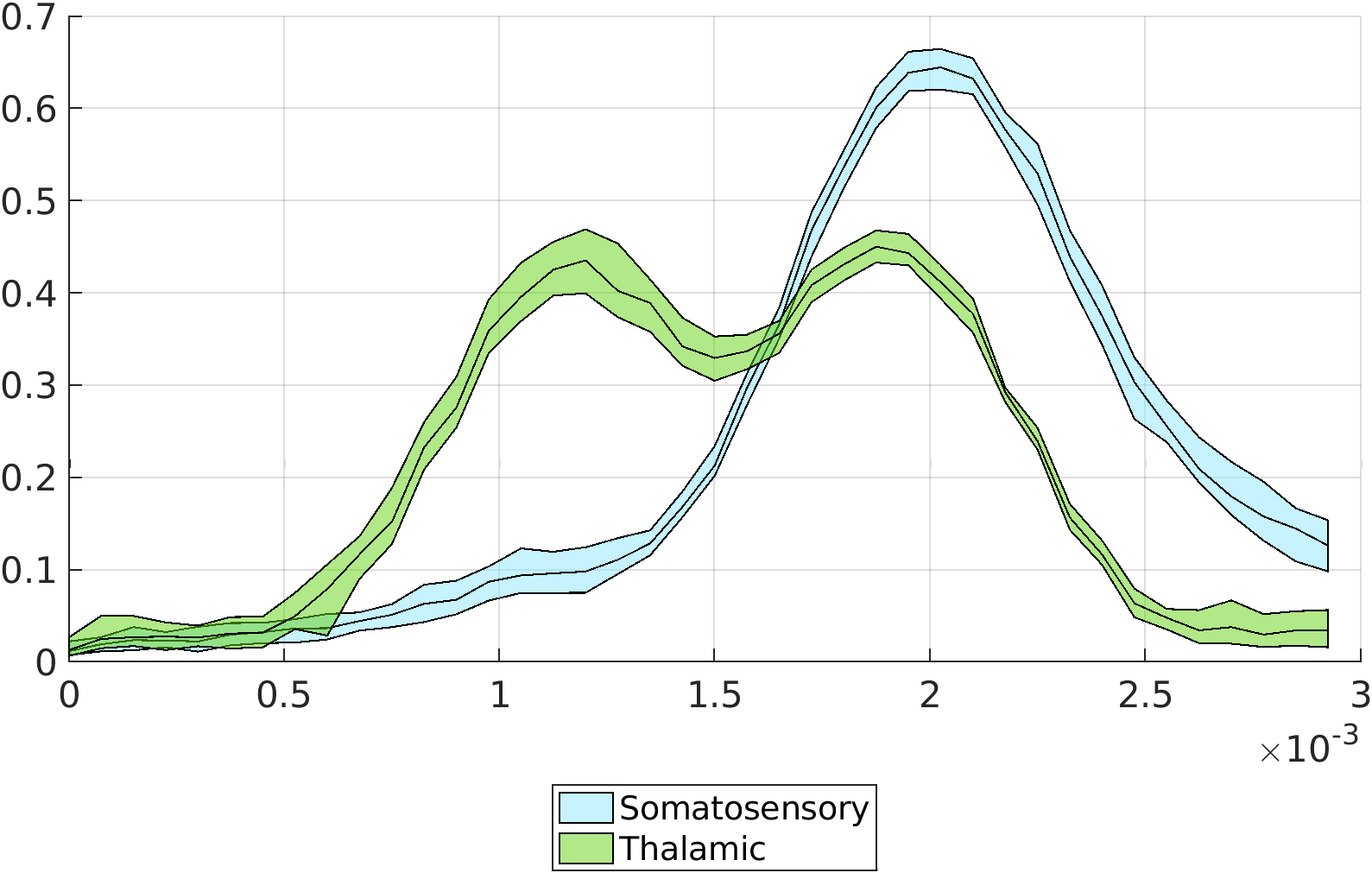}
\end{minipage}
\begin{minipage}{4.5cm}
\centering
\includegraphics[trim={0 3cm 0 0},clip, width=4.2cm]{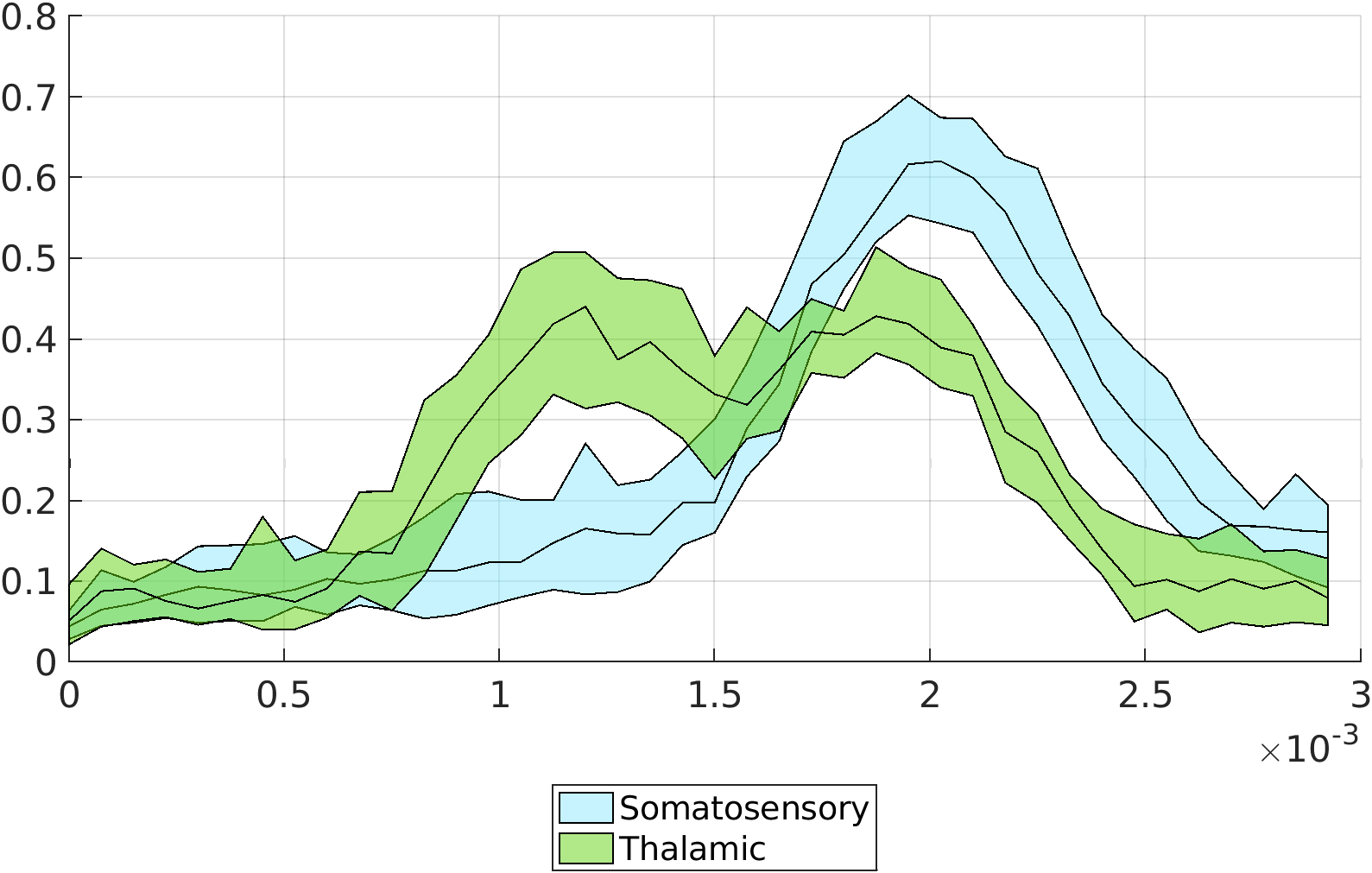}
\end{minipage}     
\begin{minipage}{4.5cm}
\centering
\includegraphics[trim={0 3cm 0 0},clip, width=4.2cm]{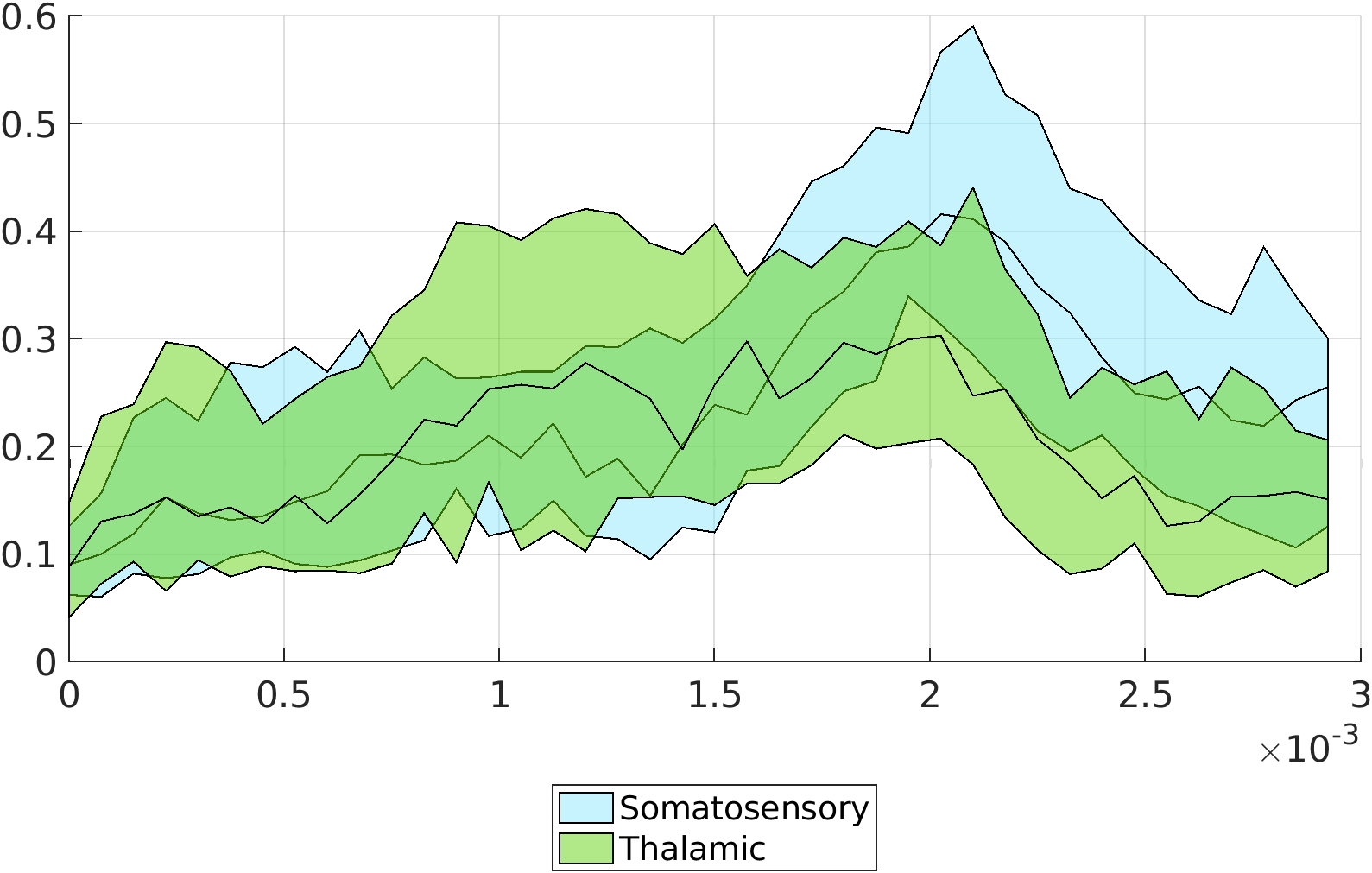}
\end{minipage}\\ 
\begin{minipage}{12cm}
\centering
\includegraphics[trim={0 0 0 15cm},clip, width=12cm]{VeikkaTimeseries/time_series_nl_10_evp_10_kf_1_sl_1_rts_0.png}
\end{minipage} 
\end{footnotesize}
    \caption{ Brain activity strength evolution in thalamus (green) and somatosensory cortex (blue). The solid, darker line represents the median over 20 estimations with different measurement noise realizations. Time runs on the x-axis, presented in milliseconds, and the y-axis shows the strength. }
    \label{fig:DynamicKFTracks}
\end{figure*}



\begin{figure}[h!]
\centering
\begin{footnotesize}
\begin{minipage}{0.5cm}
\mbox{}
\end{minipage}
\begin{minipage}{2cm}
\centering
30 dB
\end{minipage} 
\begin{minipage}{2cm}
\centering
20 dB
\end{minipage} 
\begin{minipage}{2cm}
\centering
10 dB
\end{minipage} \\ \vskip0.2cm
\begin{minipage}{0.5cm}
\rotatebox{90}{3-DSKF}
\end{minipage}
\begin{minipage}{2cm}
\centering
\includegraphics[trim={0 0cm 0 0},clip, width=1.7cm]{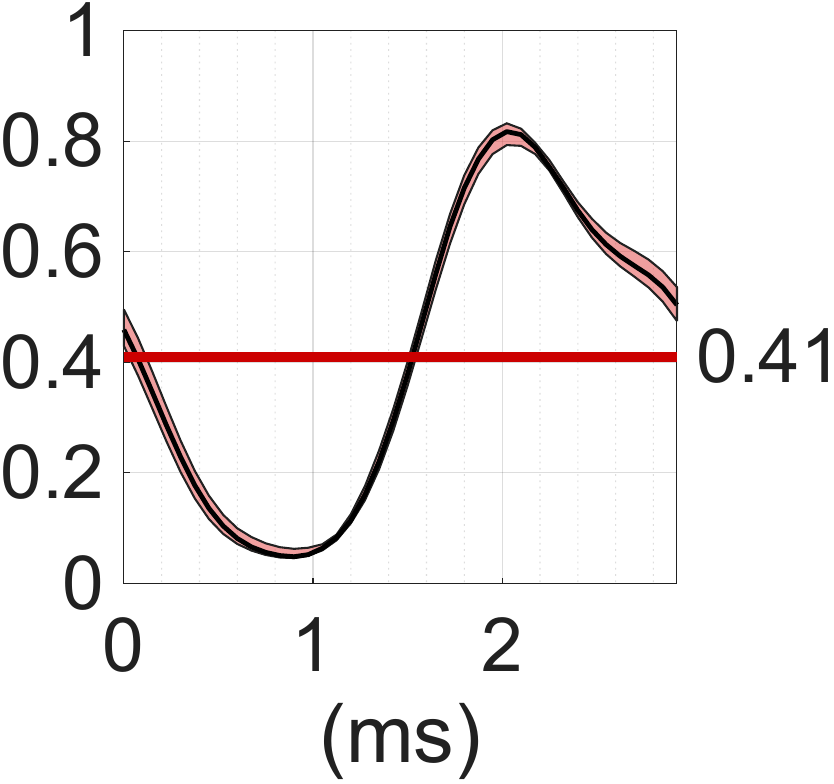}
\end{minipage}
\begin{minipage}{2cm}
\centering
\includegraphics[trim={0 0cm 0 0},clip, width=1.7cm]{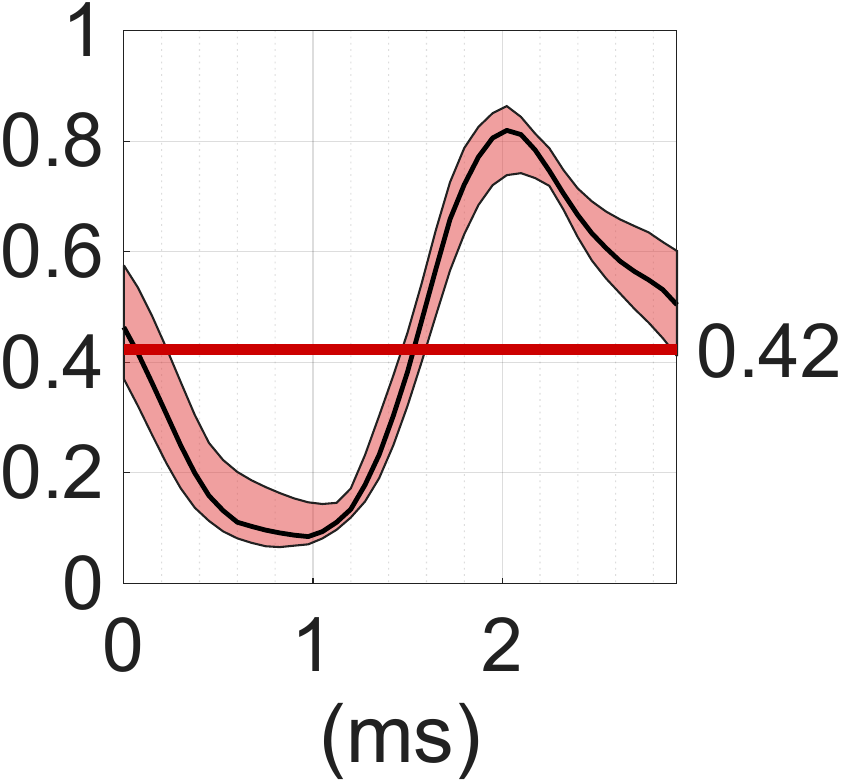}
\end{minipage}     
\begin{minipage}{2cm}
\centering
\includegraphics[trim={0 0cm 0 0},clip, width=1.7cm]{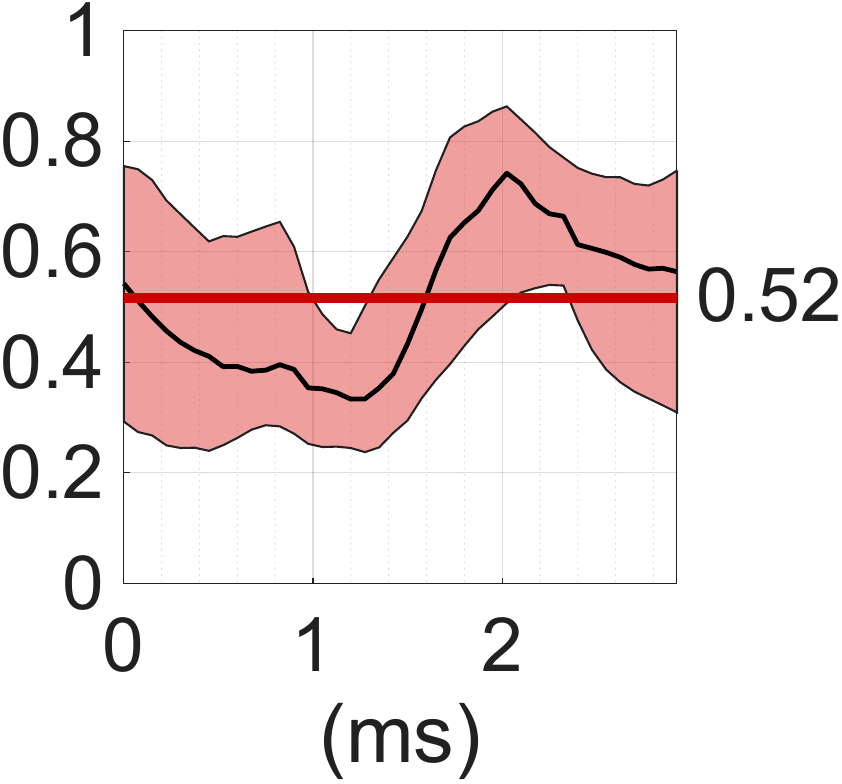}
\end{minipage}   \\ \vskip0.2cm
\begin{minipage}{0.5cm}
\rotatebox{90}{2-DSKF}
\end{minipage}
\begin{minipage}{2cm}
\centering
\includegraphics[trim={0 0cm 0 0},clip, width=1.7cm]{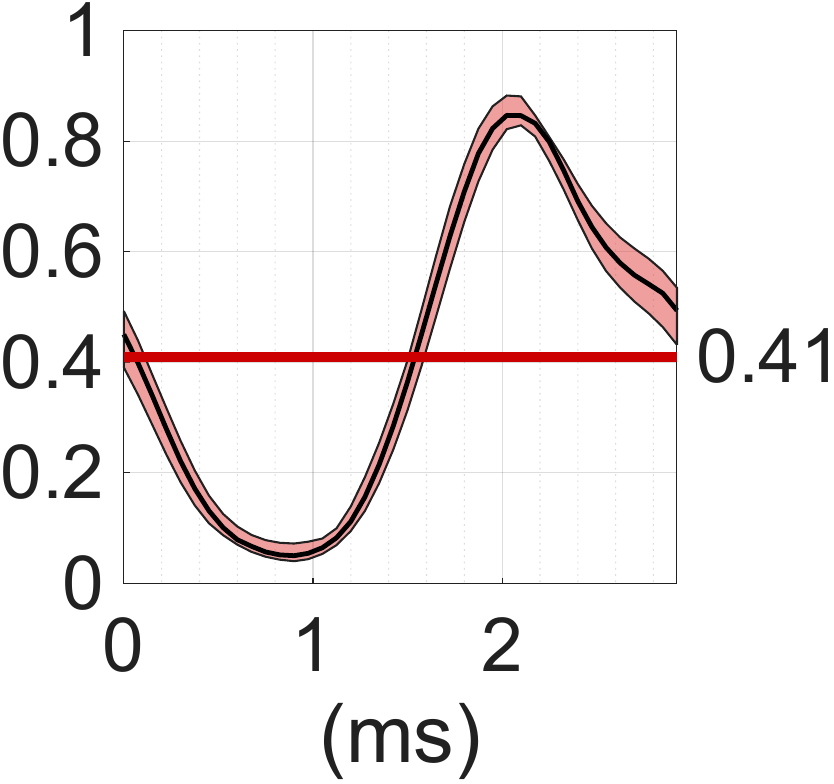}
\end{minipage}
\begin{minipage}{2cm}
\centering
\includegraphics[trim={0 0cm 0 0},clip, width=1.7cm]{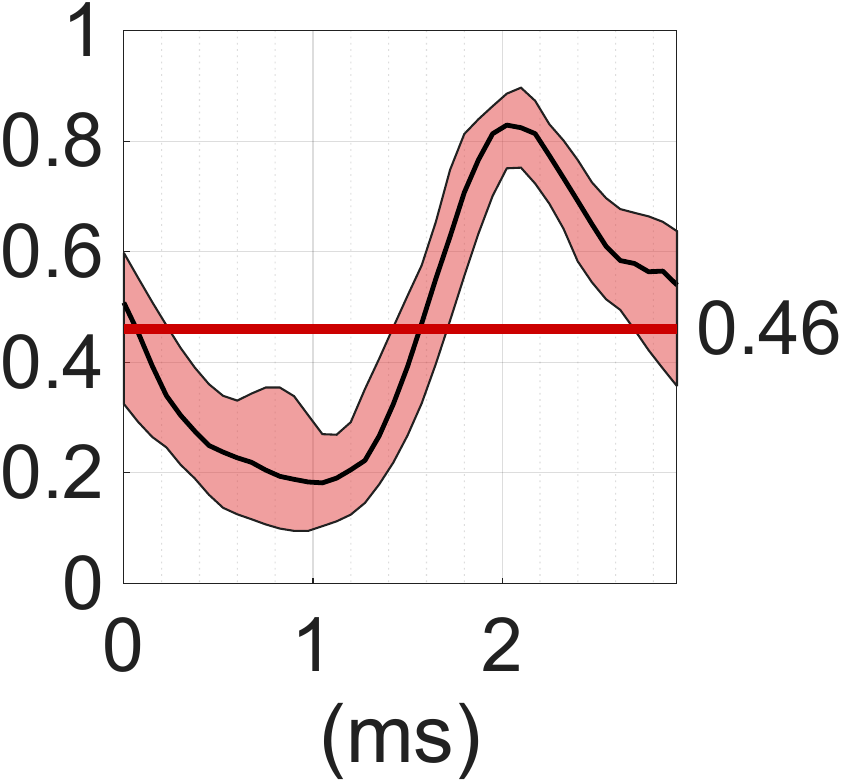}
\end{minipage}     
\begin{minipage}{2cm}
\centering
\includegraphics[trim={0 0cm 0 0},clip, width=1.7cm]{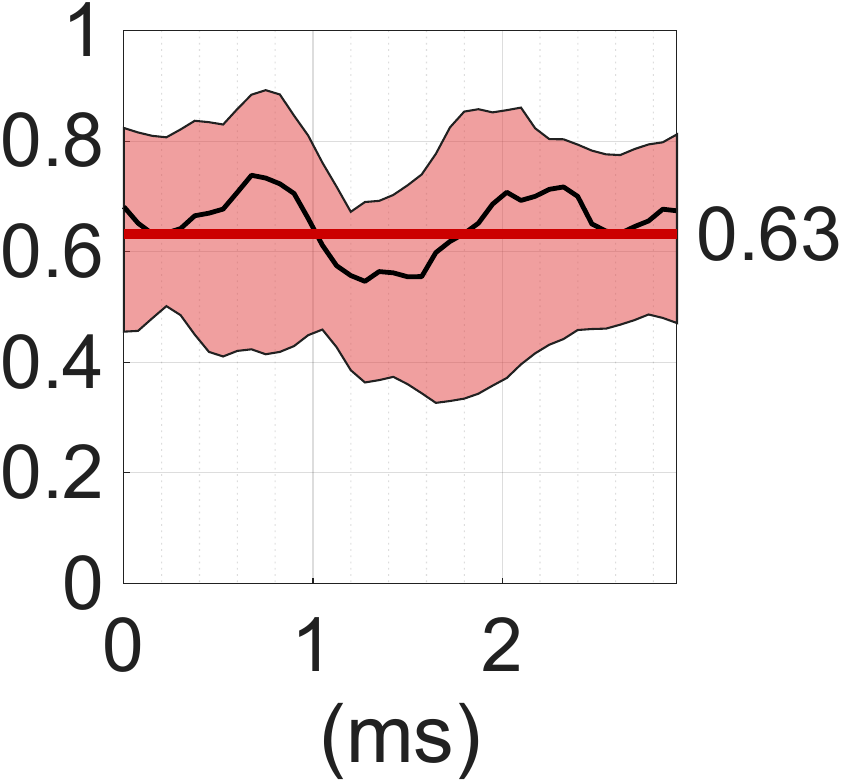}
\end{minipage}   \\ \vskip0.2cm
\begin{minipage}{0.5cm}
\rotatebox{90}{SSKF}
\end{minipage}
\begin{minipage}{2cm}
\centering
\includegraphics[trim={0 0cm 0 0},clip, width=1.7cm]{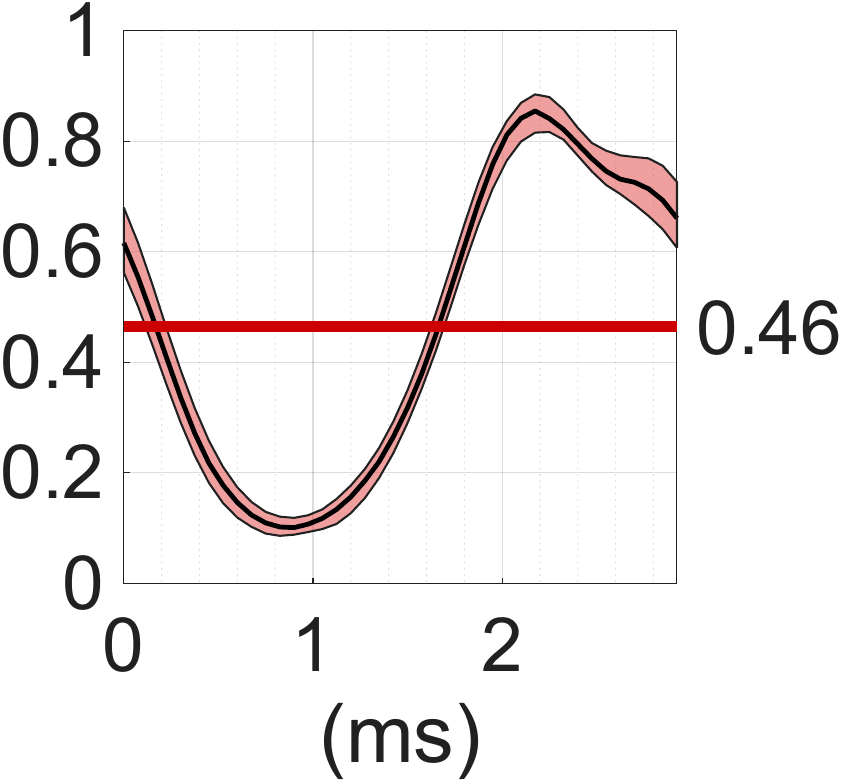}
\end{minipage}
\begin{minipage}{2cm}
\centering
\includegraphics[trim={0 0cm 0 0},clip, width=1.7cm]{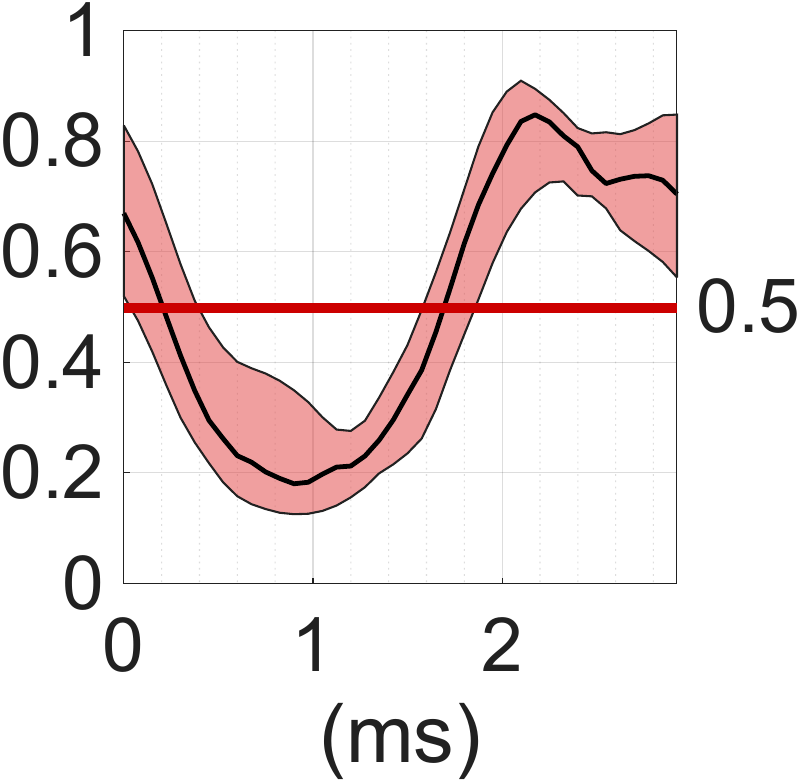}
\end{minipage}     
\begin{minipage}{2cm}
\centering
\includegraphics[trim={0 0cm 0 0},clip, width=1.7cm]{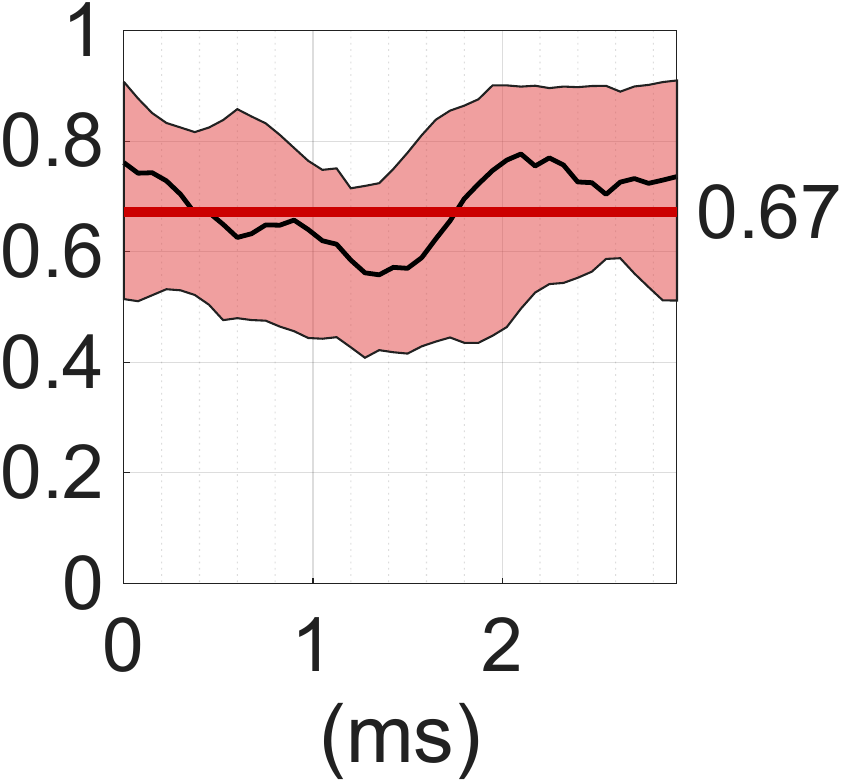}
\end{minipage}   \\ \vskip0.2cm
\begin{minipage}{0.5cm}
\rotatebox{90}{SKF}
\end{minipage}
\begin{minipage}{2cm}
\centering
\includegraphics[trim={0 0cm 0 0},clip, width=1.7cm]{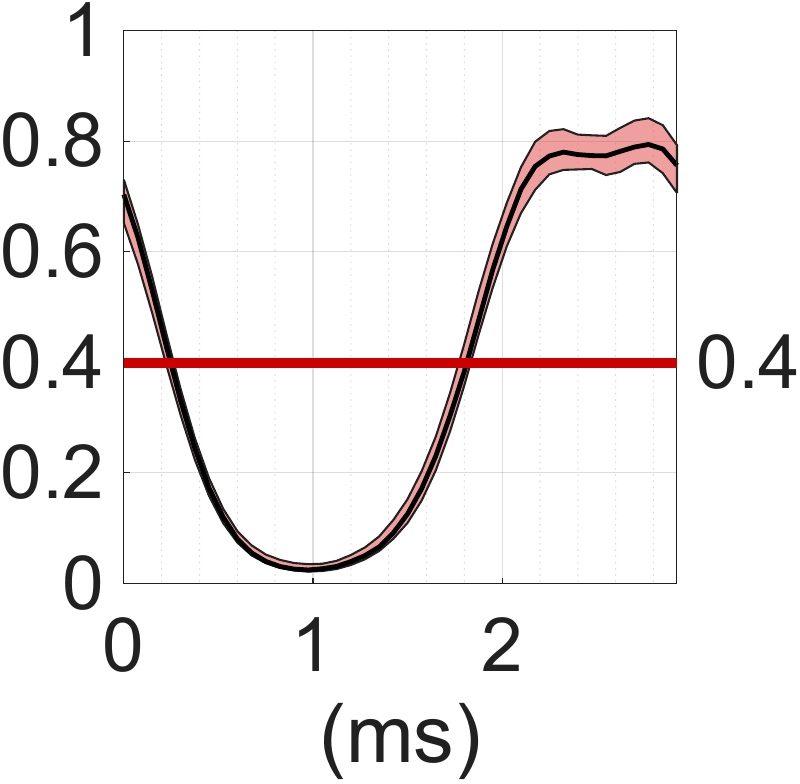}
\end{minipage}
\begin{minipage}{2cm}
\centering
\includegraphics[trim={0 0cm 0 0},clip, width=1.7cm]{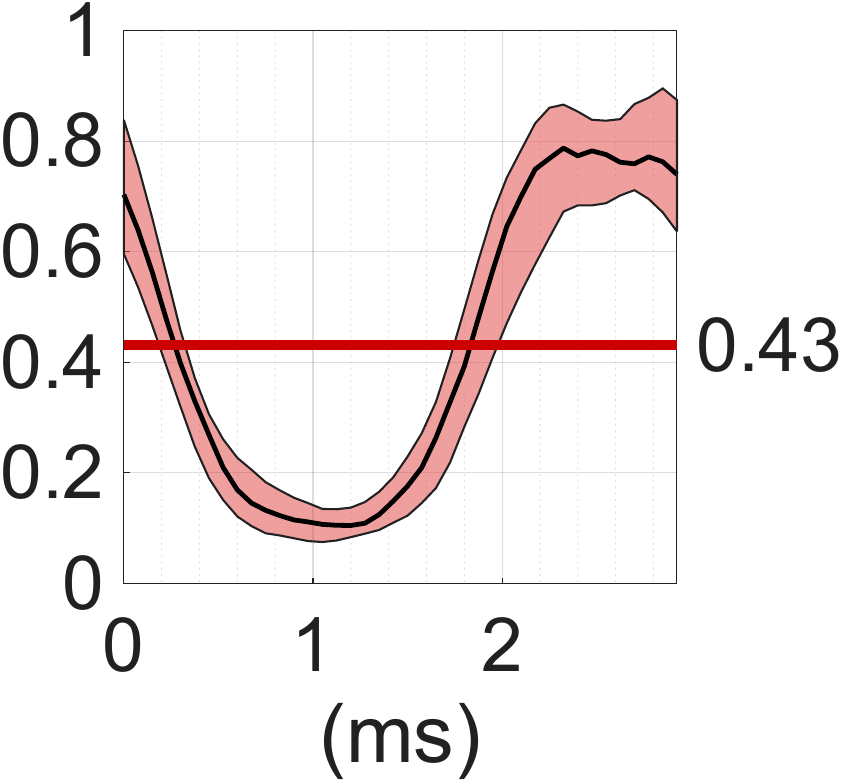}
\end{minipage}     
\begin{minipage}{2cm}
\centering
\includegraphics[trim={0 0cm 0 0},clip, width=1.7cm]{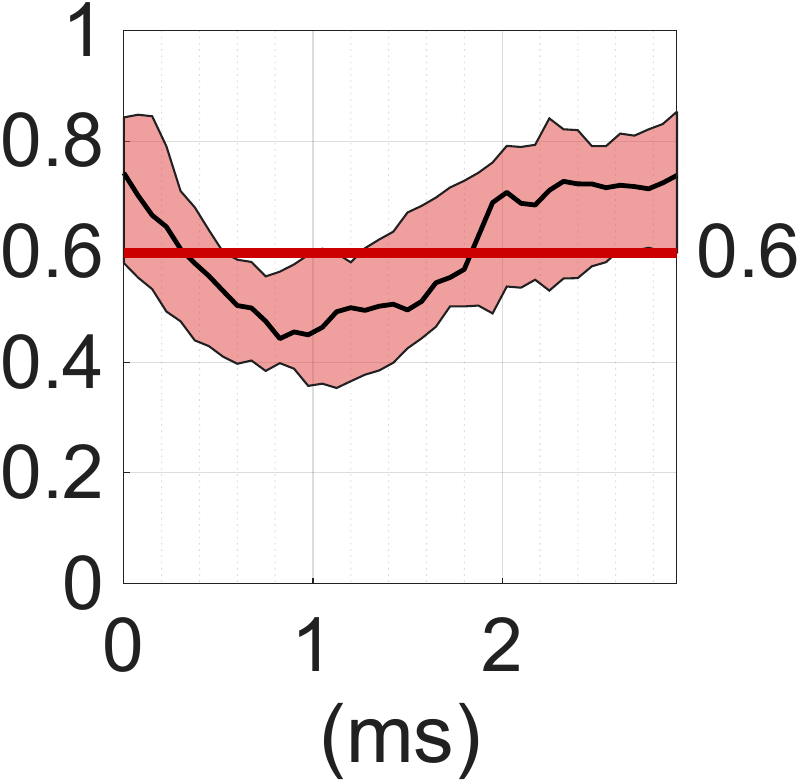}
\end{minipage}\\ 
\end{footnotesize}
    \caption{Normalized cross-correlations of the estimated tracks of the activity in somatosensory and thalamus presented in Figure \ref{fig:DynamicKFTracks}. The colored area covers the 10 and 90 \% interval around the median. A bold, solid curve presents the median track. The solid red horizontal line is the mean cross-correlation over the whole time period. This mean value is presented numerically on the right side of the graphs.}
    \label{fig:DynamicKFCrossCorr}
\end{figure}

\begin{figure}[h!]
\centering
\begin{footnotesize}
\begin{minipage}{0.5cm}
\mbox{}
\end{minipage}
\begin{minipage}{2cm}
\centering
30 dB
\end{minipage} 
\begin{minipage}{2cm}
\centering
20 dB
\end{minipage} 
\begin{minipage}{2cm}
\centering
10 dB
\end{minipage} \\ \vskip0.2cm
\begin{minipage}{0.5cm}
\rotatebox{90}{3-DSKF}
\end{minipage}
\begin{minipage}{2cm}
\centering
\includegraphics[trim={0 0cm 0 0},clip, width=1.7cm]{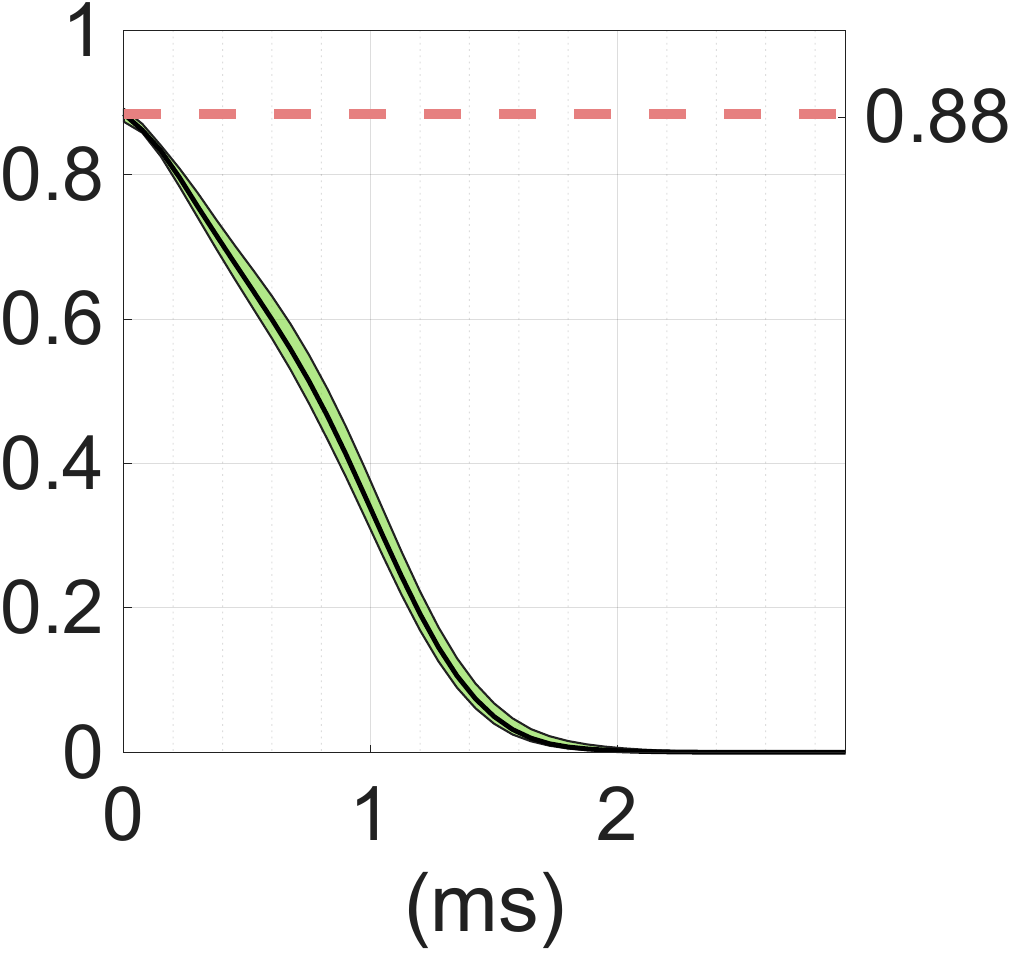}
\end{minipage}
\begin{minipage}{2cm}
\centering
\includegraphics[trim={0 0cm 0 0},clip, width=1.7cm]{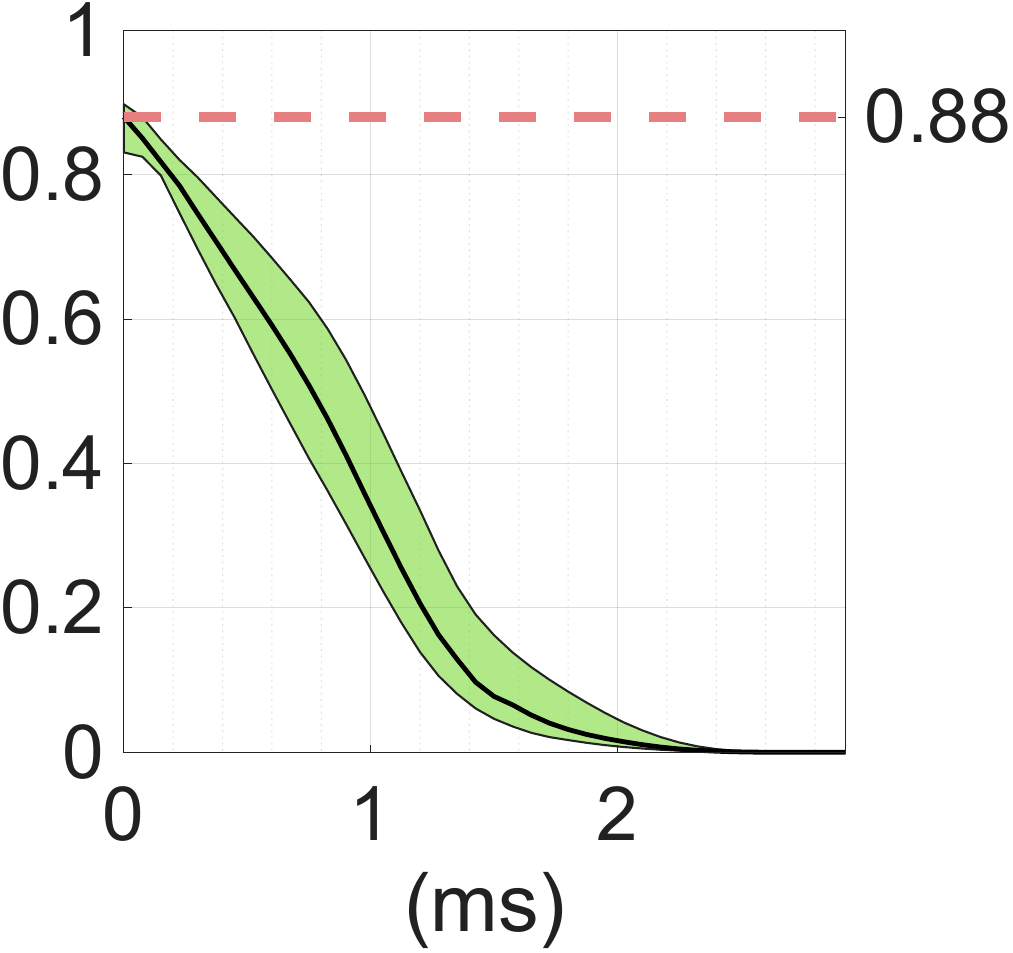}
\end{minipage}     
\begin{minipage}{2cm}
\centering
\includegraphics[trim={0 0cm 0 0},clip, width=1.7cm]{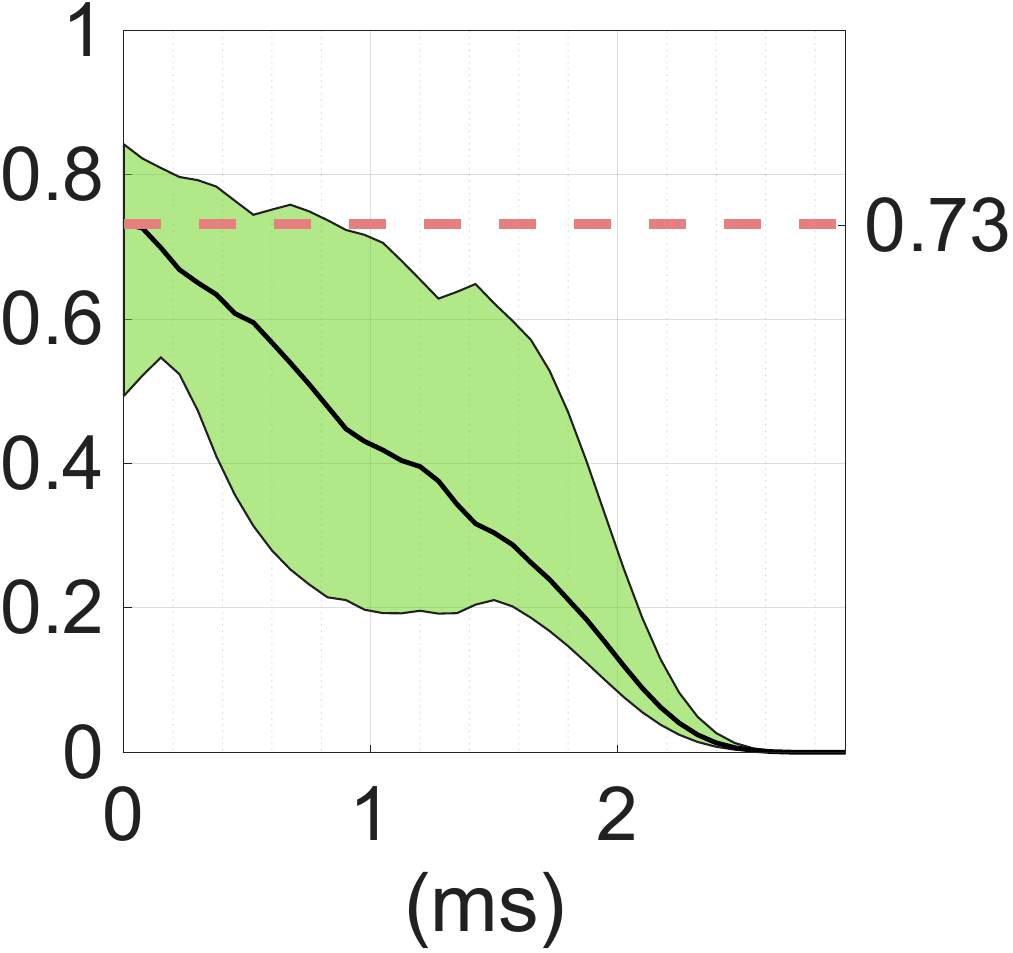}
\end{minipage}   \\ \vskip0.2cm
\begin{minipage}{0.5cm}
\rotatebox{90}{2-DSKF}
\end{minipage}
\begin{minipage}{2cm}
\centering
\includegraphics[trim={0 0cm 0 0},clip, width=1.7cm]{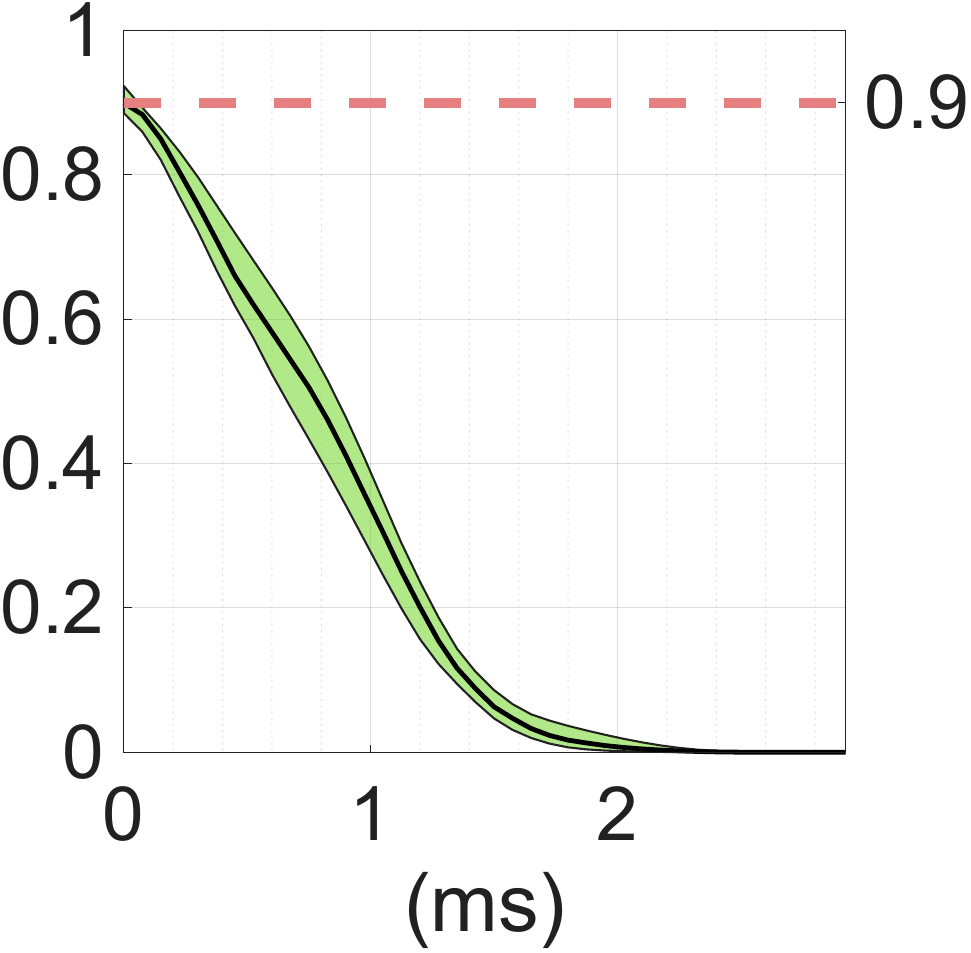}
\end{minipage}
\begin{minipage}{2cm}
\centering
\includegraphics[trim={0 0cm 0 0},clip, width=1.7cm]{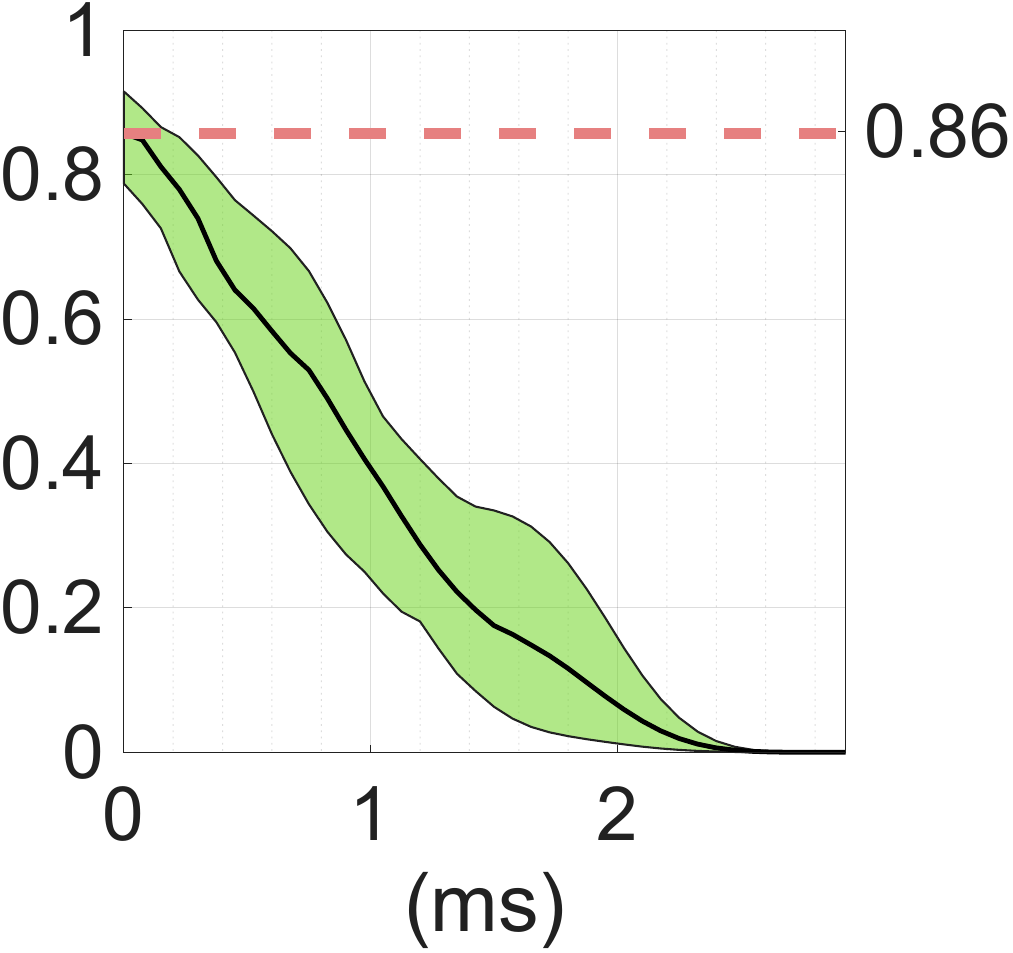}
\end{minipage}     
\begin{minipage}{2cm}
\centering
\includegraphics[trim={0 0cm 0 0},clip, width=1.7cm]{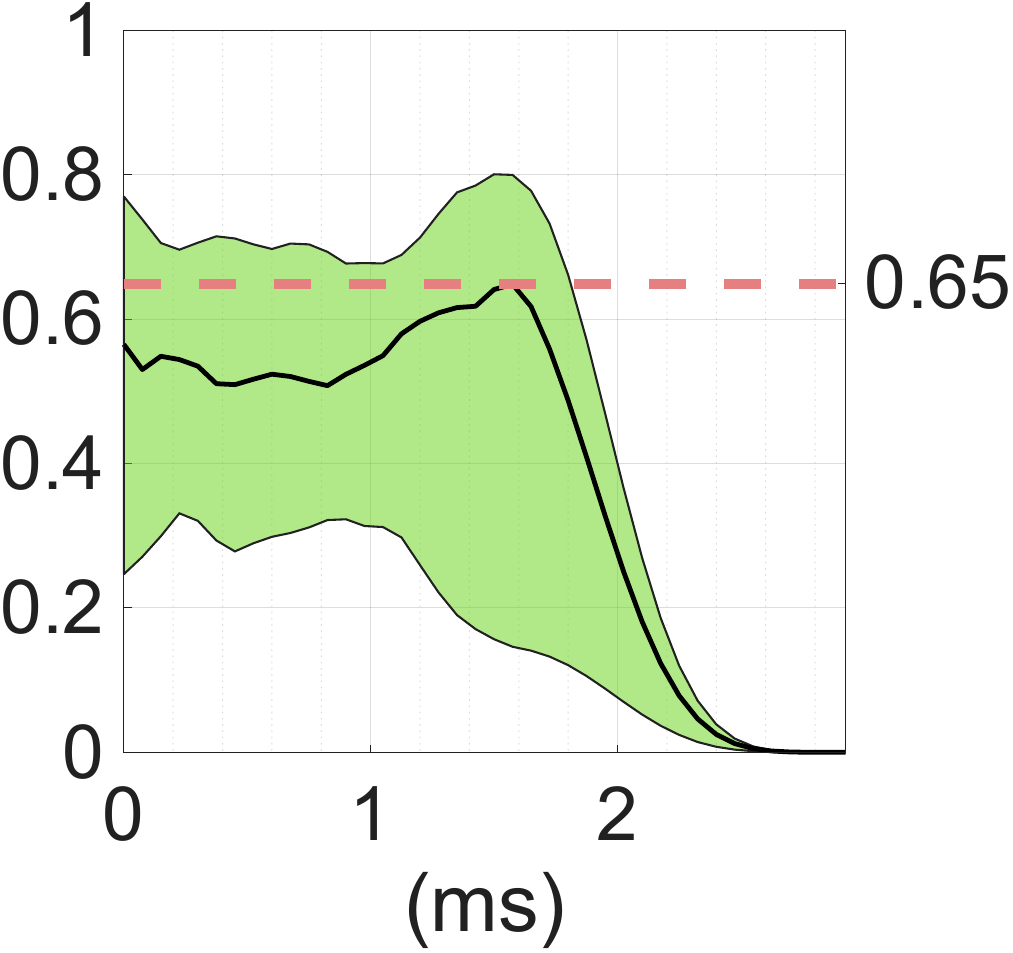}
\end{minipage}   \\ \vskip0.2cm
\begin{minipage}{0.5cm}
\rotatebox{90}{SSKF}
\end{minipage}
\begin{minipage}{2cm}
\centering
\includegraphics[trim={0 0cm 0 0},clip, width=1.7cm]{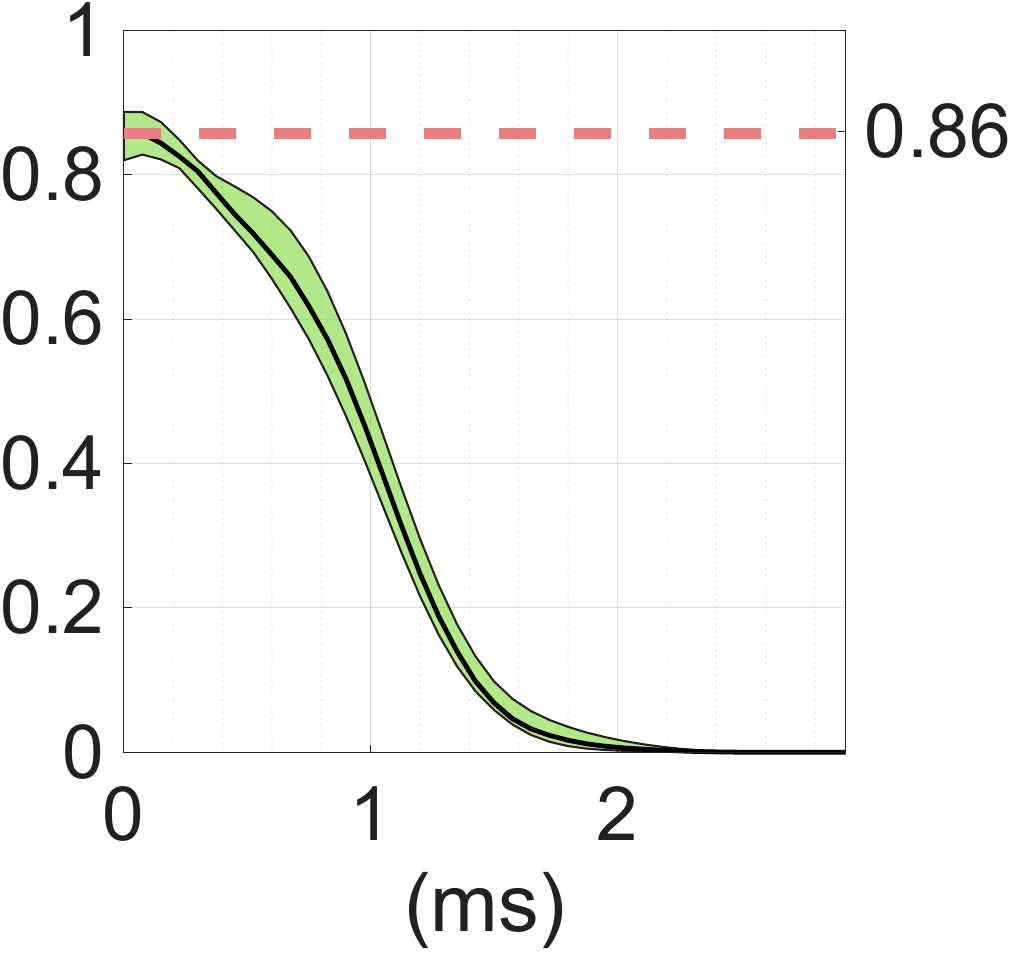}
\end{minipage}
\begin{minipage}{2cm}
\centering
\includegraphics[trim={0 0cm 0 0},clip, width=1.7cm]{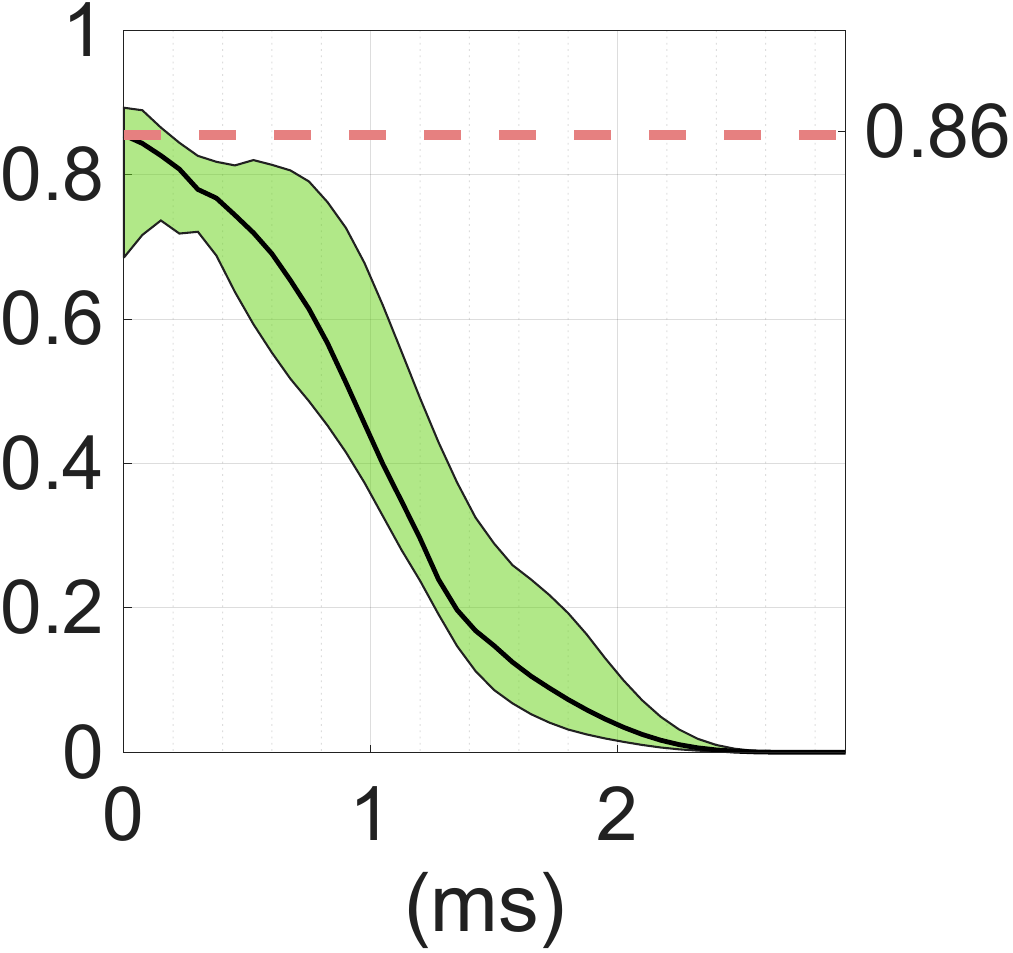}
\end{minipage}     
\begin{minipage}{2cm}
\centering
\includegraphics[trim={0 0cm 0 0},clip, width=1.7cm]{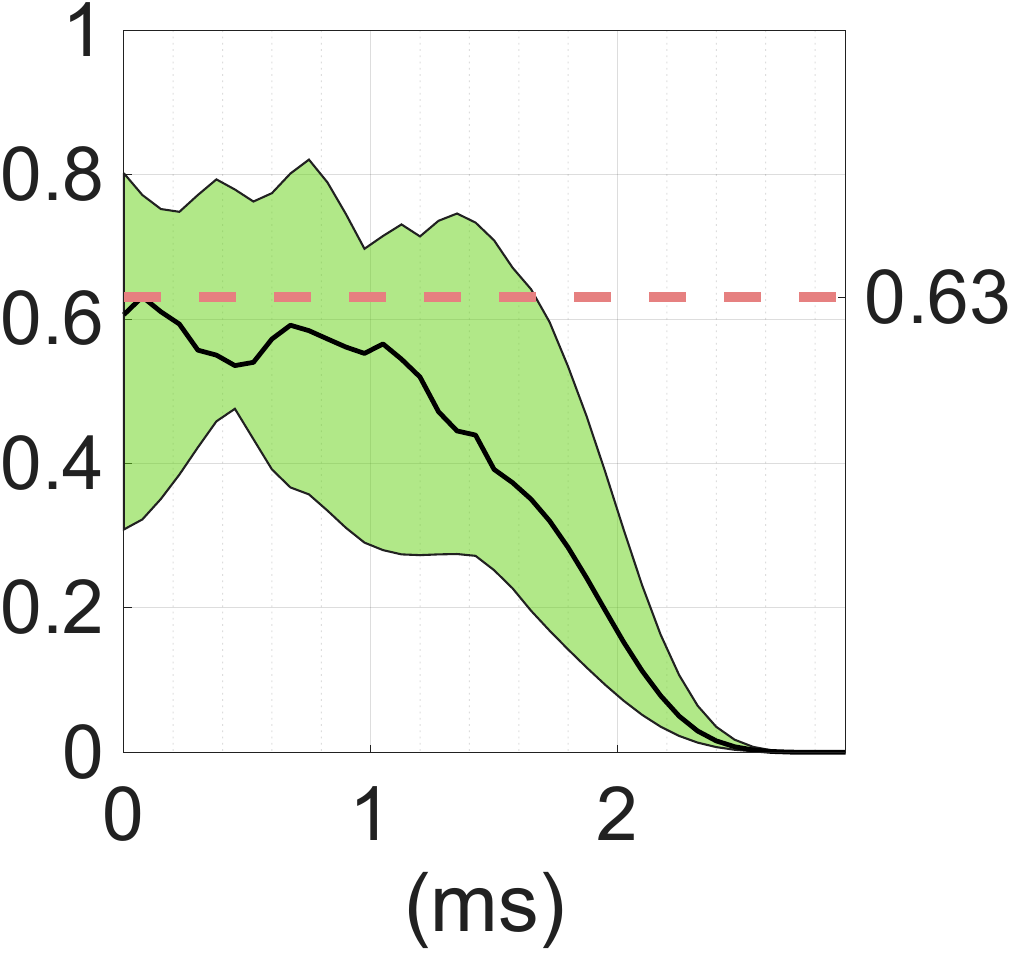}
\end{minipage}   \\ \vskip0.2cm
\begin{minipage}{0.5cm}
\rotatebox{90}{SKF}
\end{minipage}
\begin{minipage}{2cm}
\centering
\includegraphics[trim={0 0cm 0 0},clip, width=1.7cm]{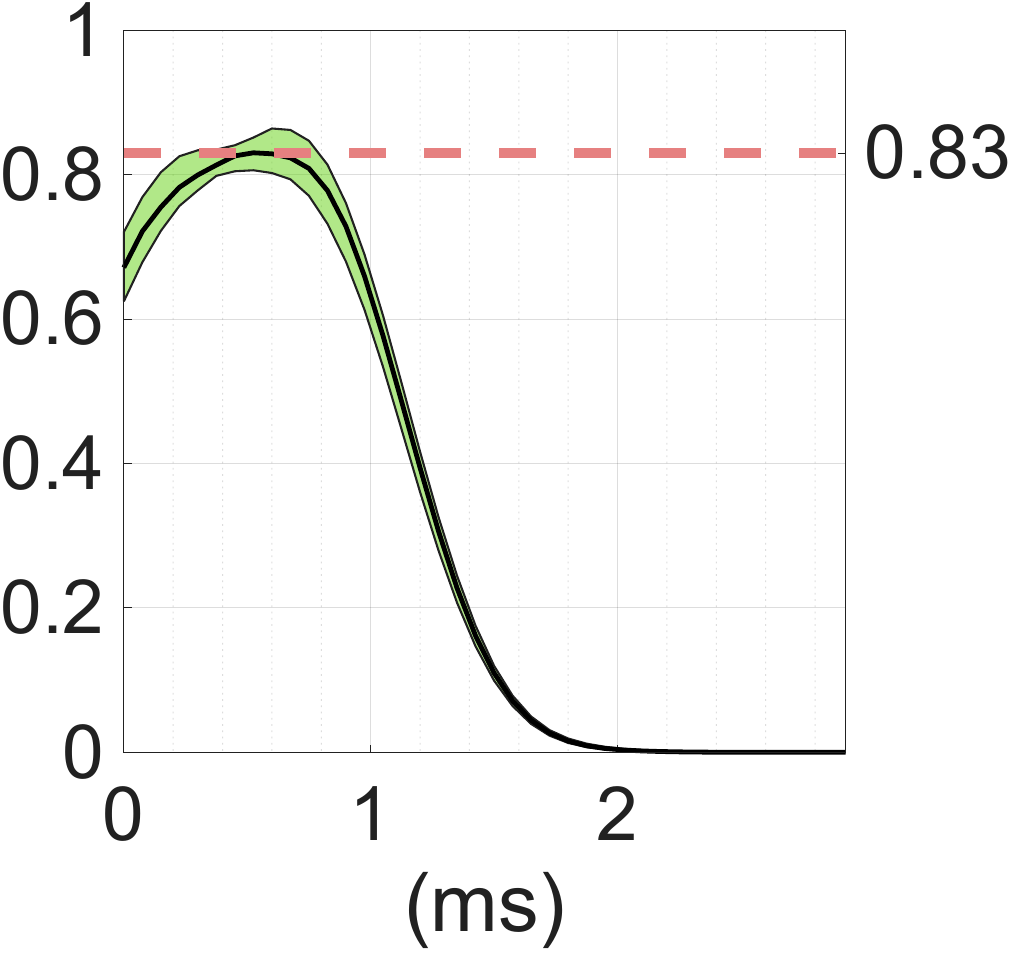}
\end{minipage}
\begin{minipage}{2cm}
\centering
\includegraphics[trim={0 0cm 0 0},clip, width=1.7cm]{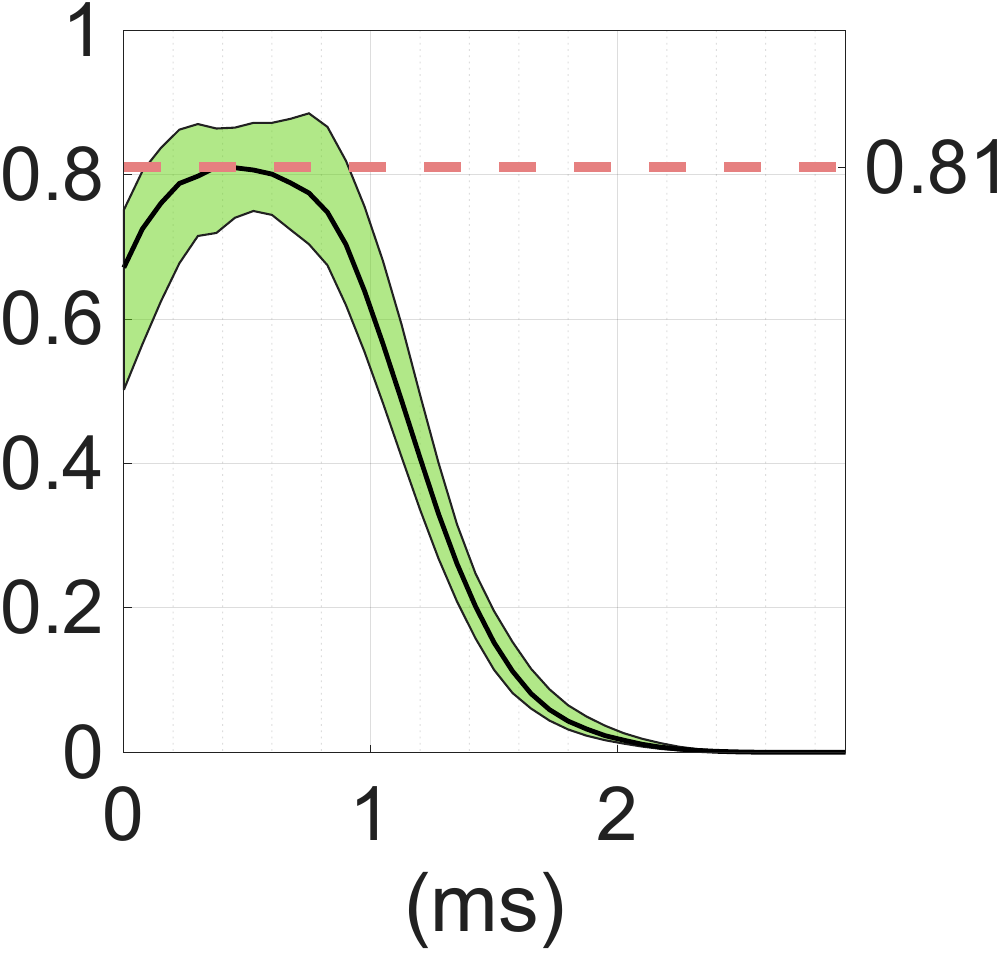}
\end{minipage}     
\begin{minipage}{2cm}
\centering
\includegraphics[trim={0 0cm 0 0},clip, width=1.7cm]{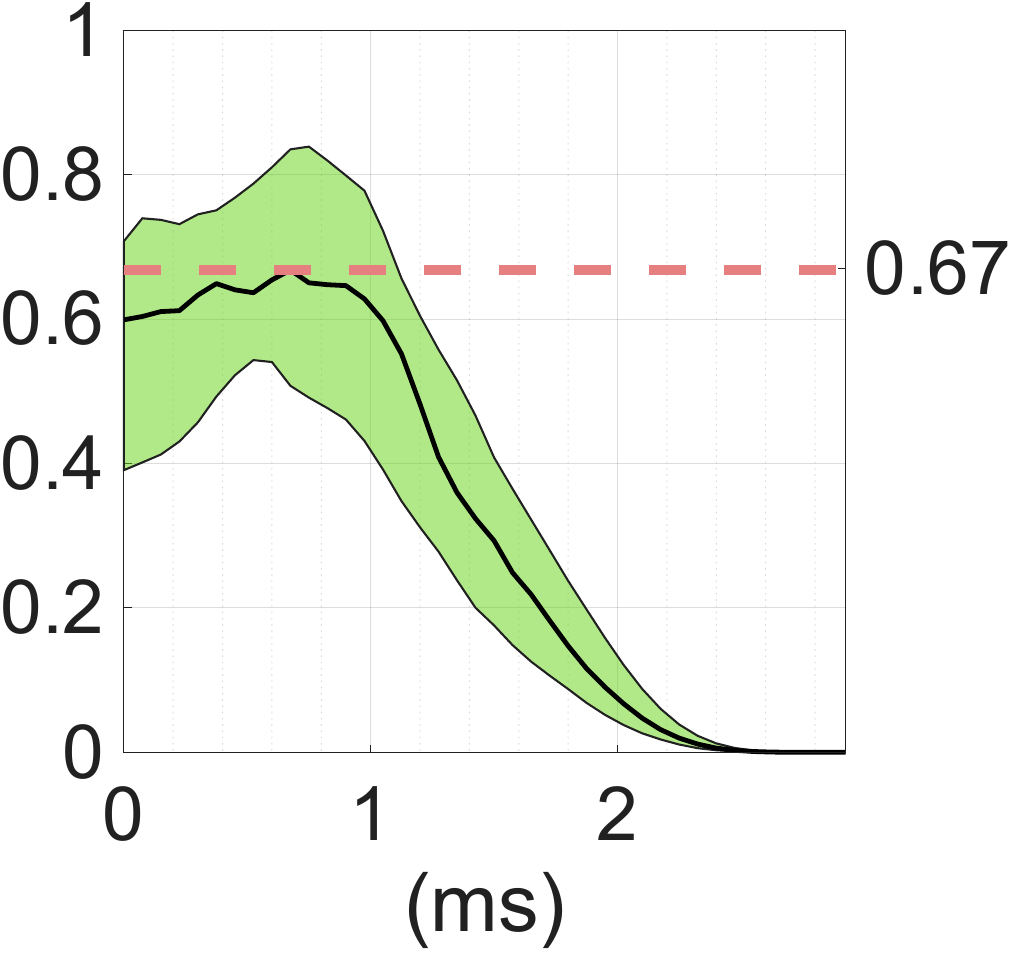}
\end{minipage}\\ 
\end{footnotesize}
    \caption{Normalized cross-correlations of the estimated tracks of the thalamic activity (presented by green in Figure \ref{fig:DynamicKFTracks}) against the true track presented in Figure \ref{fig:trueevolution}. The colored area covers the 10 and 90 \% interval around the median. A bold, solid curve presents the median track. The dashed light red horizontal line is the highest point of the median cross-correlation. This maximum value is presented numerically on the right side of the graphs.}
    \label{fig:DynamicKFCrossCorrTrue1}
\end{figure}

\begin{figure}[h!]
\centering
\begin{footnotesize}
\begin{minipage}{0.5cm}
\mbox{}
\end{minipage}
\begin{minipage}{2cm}
\centering
30 dB
\end{minipage} 
\begin{minipage}{2cm}
\centering
20 dB
\end{minipage} 
\begin{minipage}{2cm}
\centering
10 dB
\end{minipage} \\ \vskip0.2cm
\begin{minipage}{0.5cm}
\rotatebox{90}{3-DSKF}
\end{minipage}
\begin{minipage}{2cm}
\centering
\includegraphics[trim={0 0cm 0 0},clip, width=1.7cm]{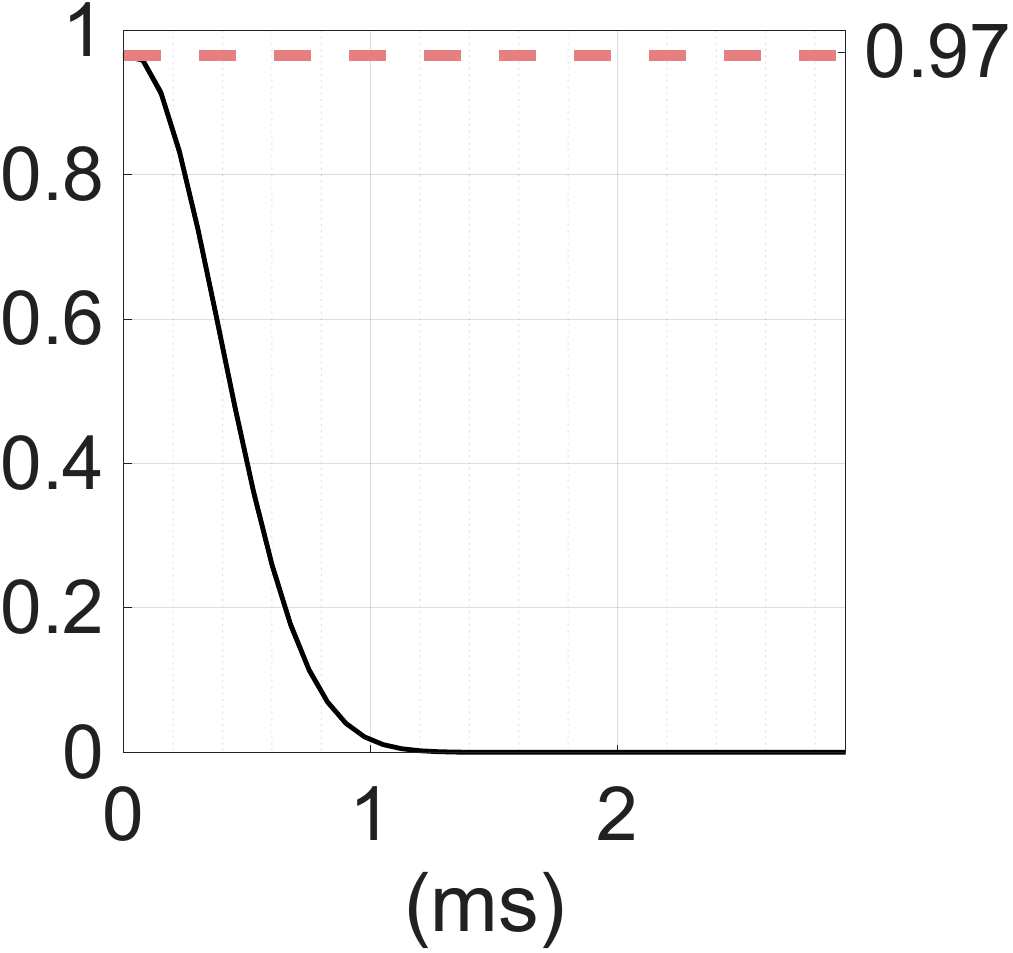}
\end{minipage}
\begin{minipage}{2cm}
\centering
\includegraphics[trim={0 0cm 0 0},clip, width=1.7cm]{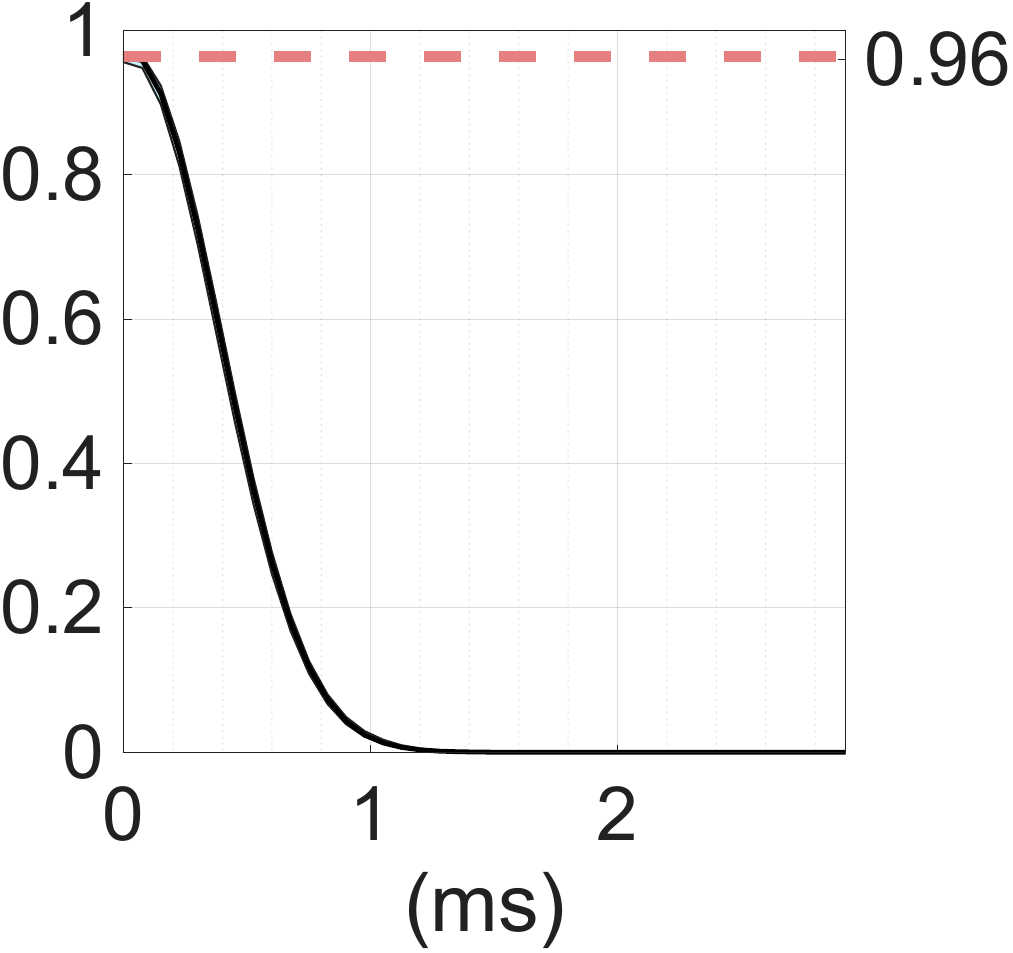}
\end{minipage}     
\begin{minipage}{2cm}
\centering
\includegraphics[trim={0 0cm 0 0},clip, width=1.7cm]{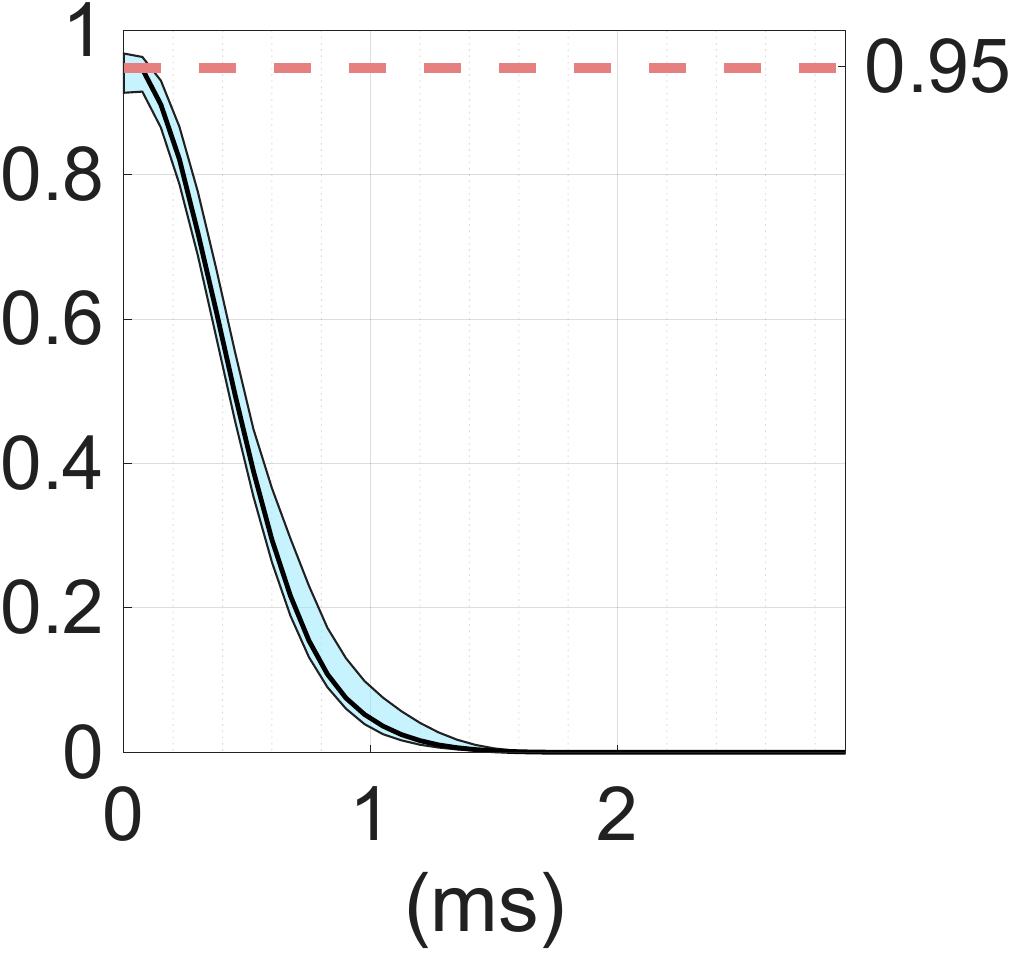}
\end{minipage}   \\ \vskip0.2cm
\begin{minipage}{0.5cm}
\rotatebox{90}{2-DSKF}
\end{minipage}
\begin{minipage}{2cm}
\centering
\includegraphics[trim={0 0cm 0 0},clip, width=1.7cm]{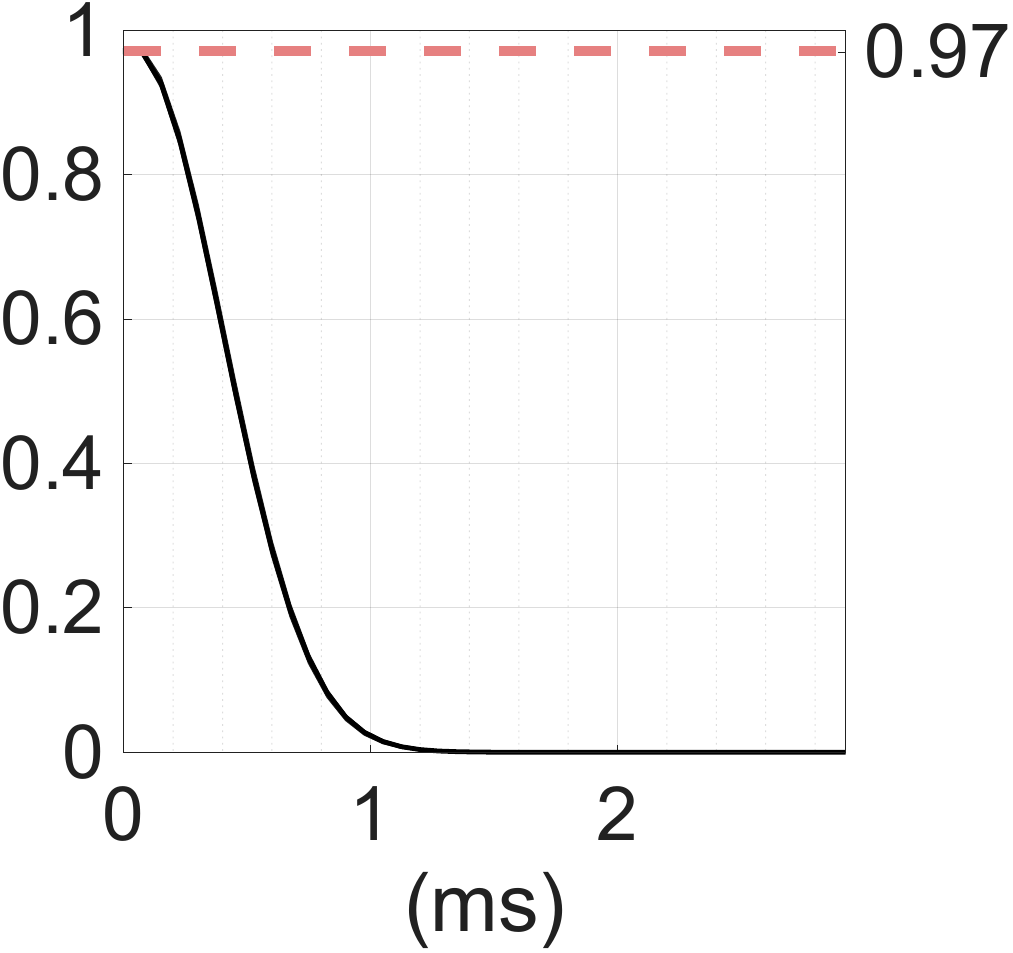}
\end{minipage}
\begin{minipage}{2cm}
\centering
\includegraphics[trim={0 0cm 0 0},clip, width=1.7cm]{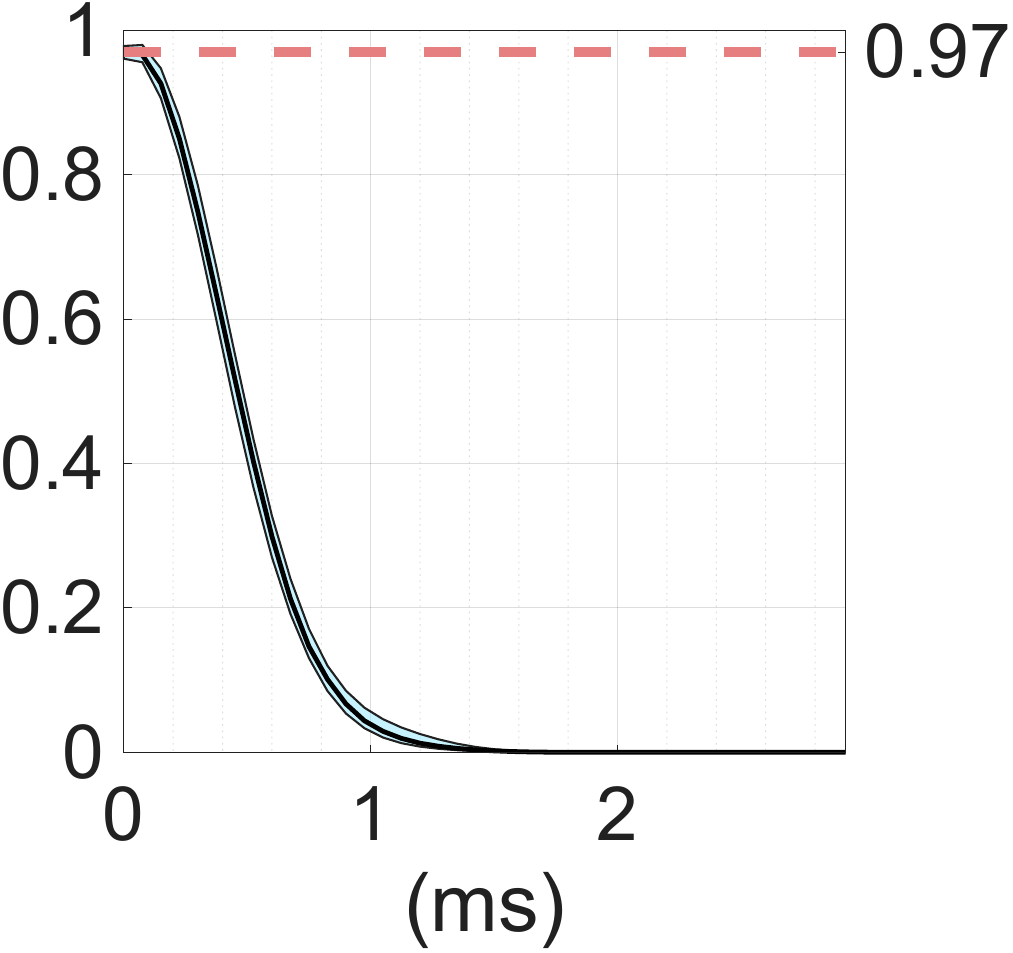}
\end{minipage}     
\begin{minipage}{2cm}
\centering
\includegraphics[trim={0 0cm 0 0},clip, width=1.7cm]{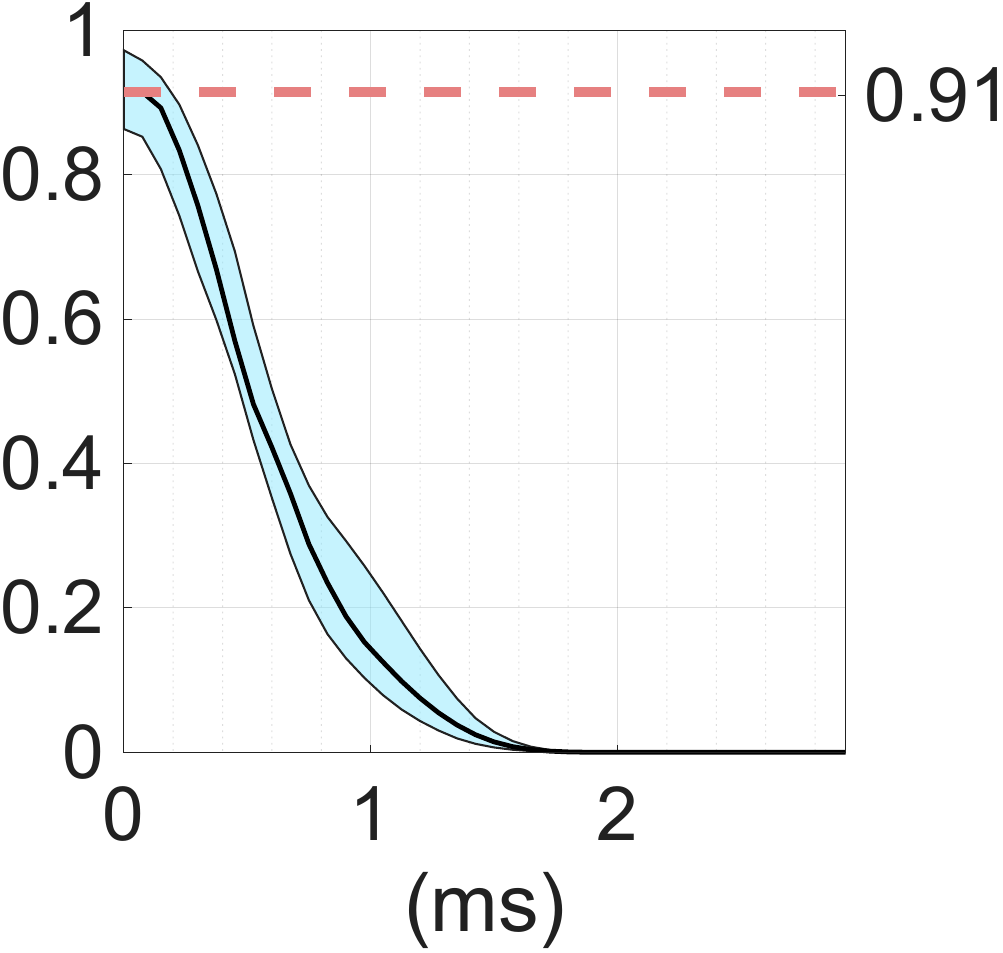}
\end{minipage}   \\ \vskip0.2cm
\begin{minipage}{0.5cm}
\rotatebox{90}{SSKF}
\end{minipage}
\begin{minipage}{2cm}
\centering
\includegraphics[trim={0 0cm 0 0},clip, width=1.7cm]{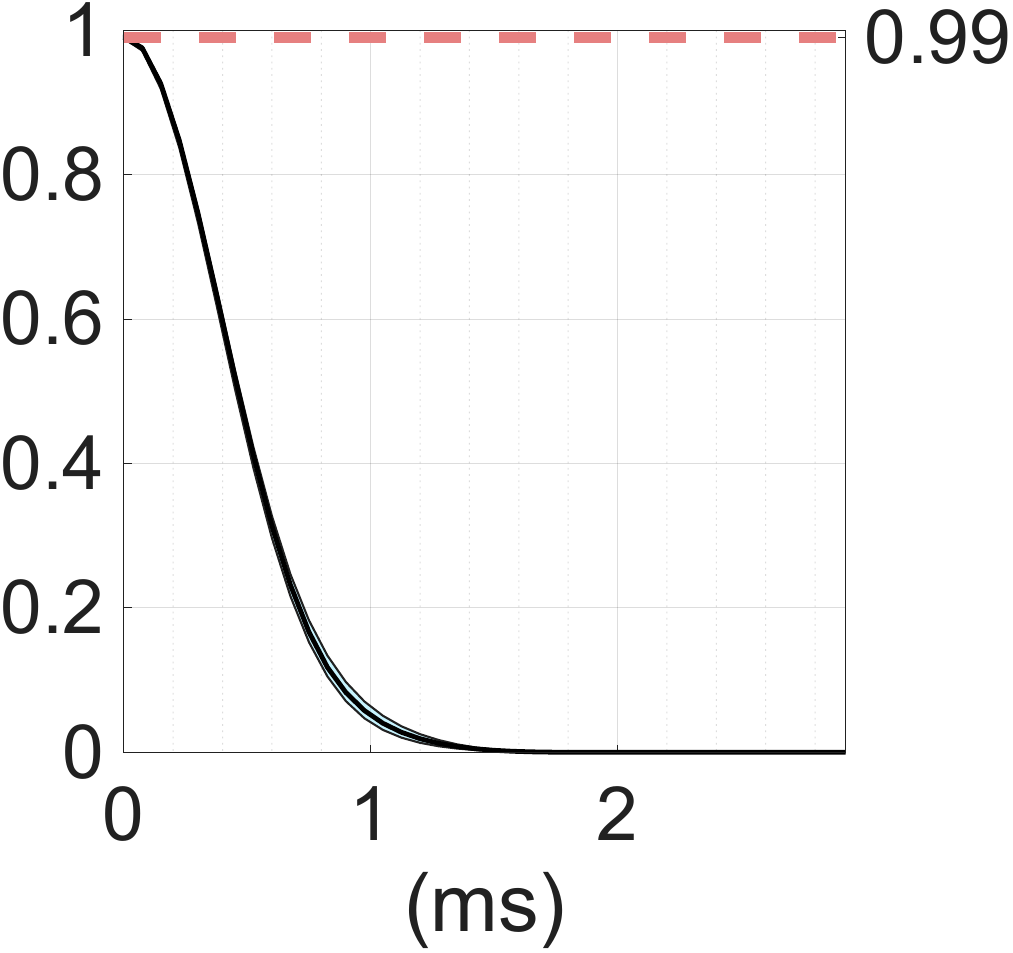}
\end{minipage}
\begin{minipage}{2cm}
\centering
\includegraphics[trim={0 0cm 0 0},clip, width=1.7cm]{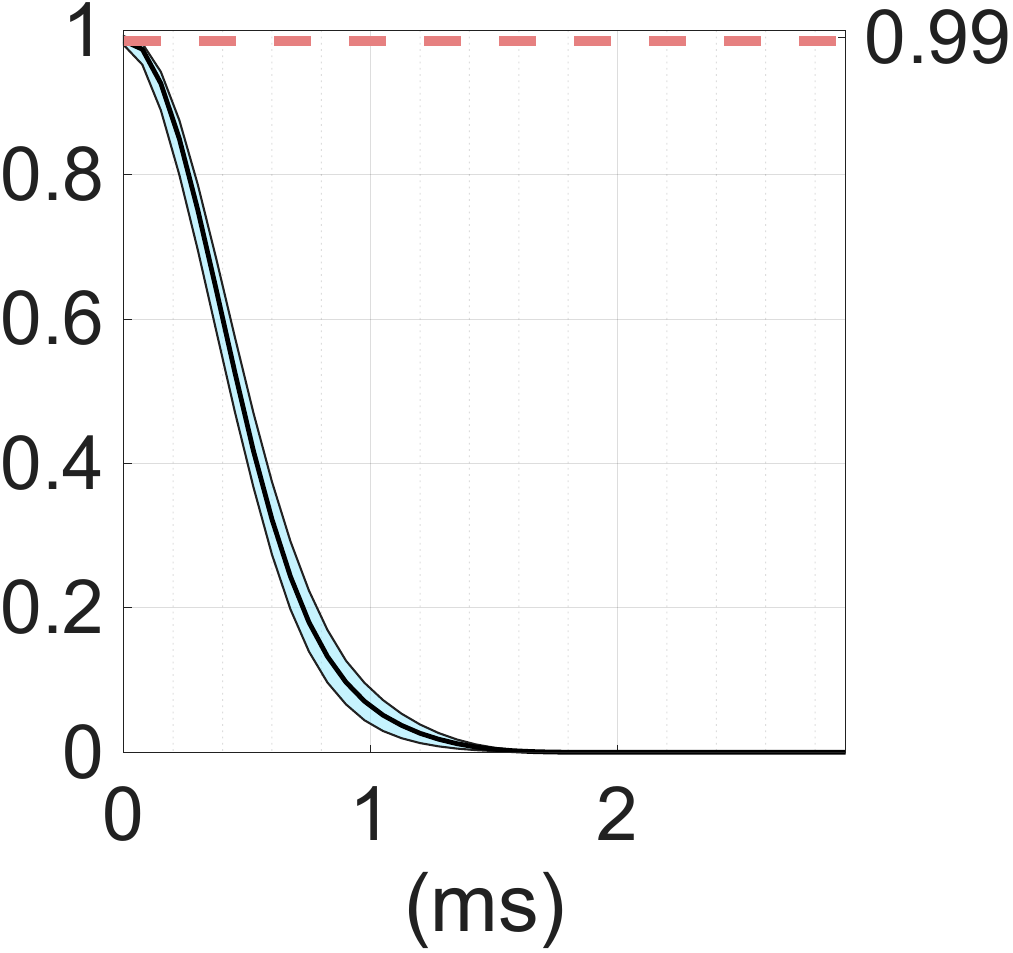}
\end{minipage}     
\begin{minipage}{2cm}
\centering
\includegraphics[trim={0 0cm 0 0},clip, width=1.7cm]{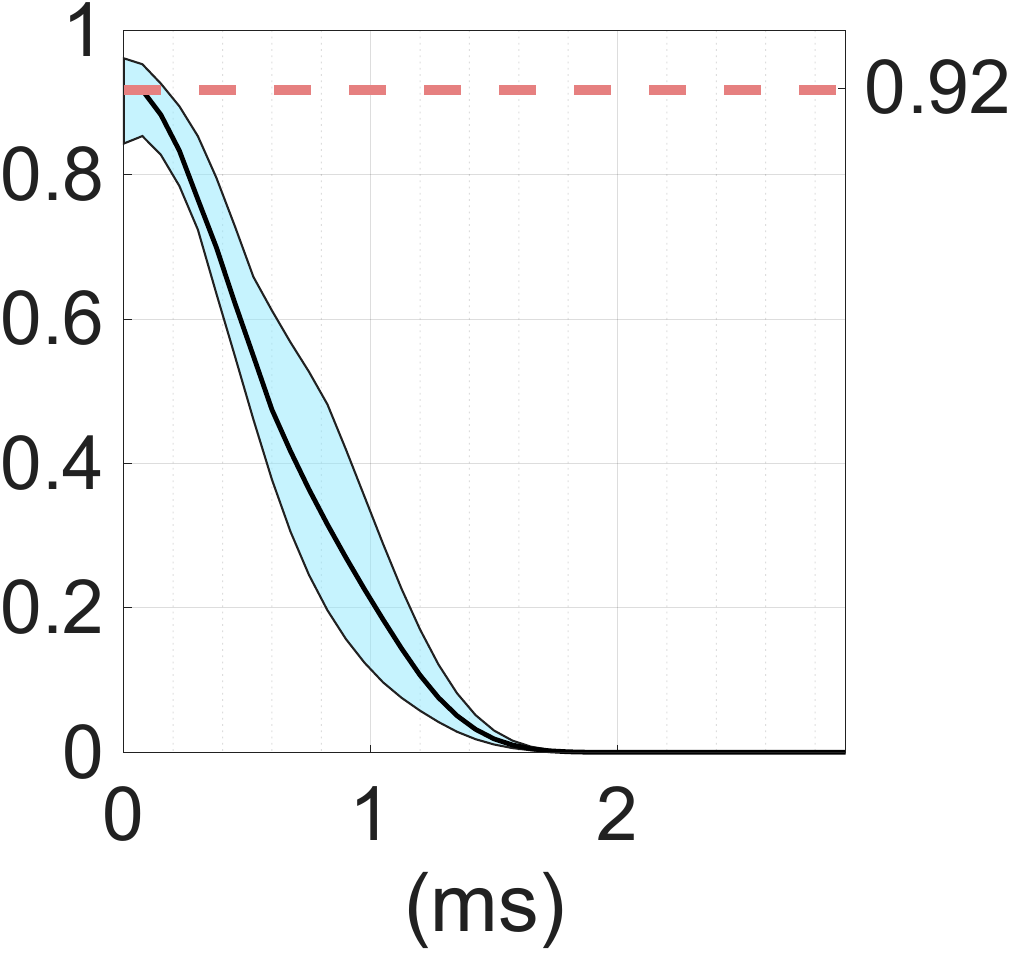}
\end{minipage}   \\ \vskip0.2cm
\begin{minipage}{0.5cm}
\rotatebox{90}{SKF}
\end{minipage}
\begin{minipage}{2cm}
\centering
\includegraphics[trim={0 0cm 0 0},clip, width=1.7cm]{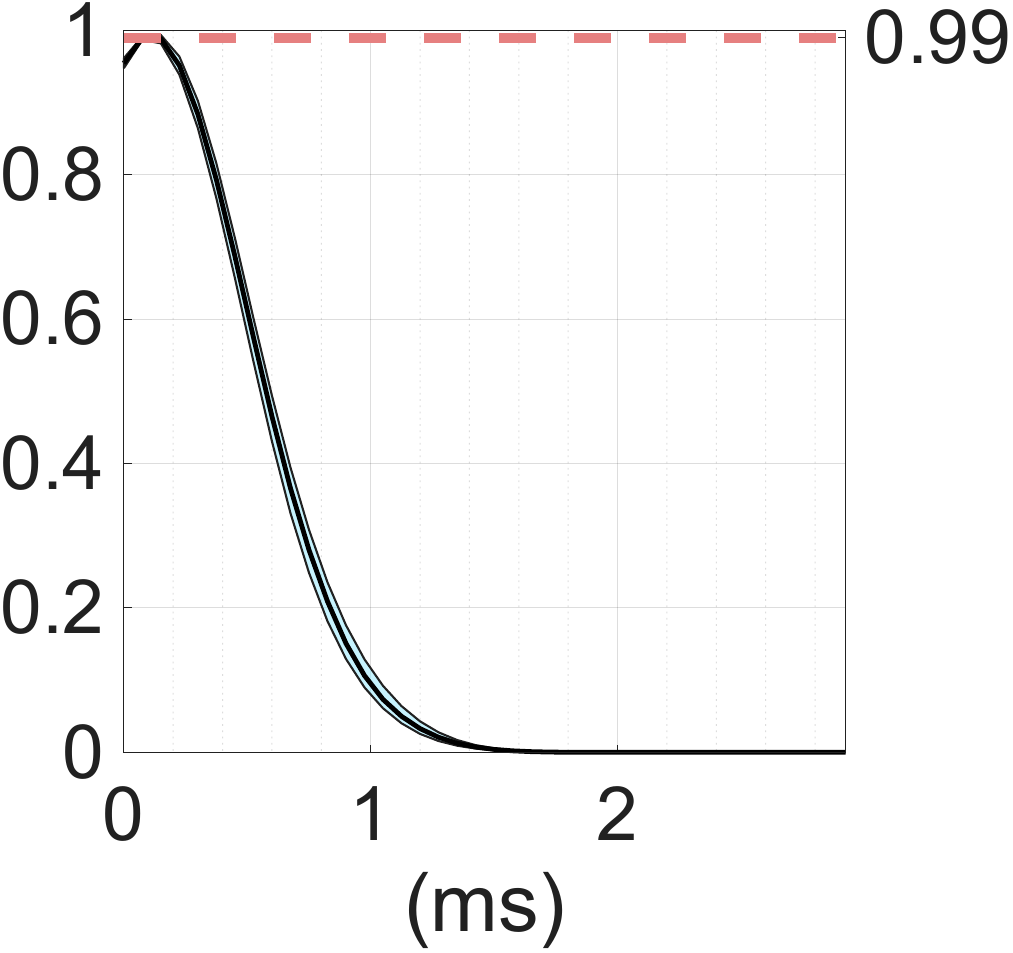}
\end{minipage}
\begin{minipage}{2cm}
\centering
\includegraphics[trim={0 0cm 0 0},clip, width=1.7cm]{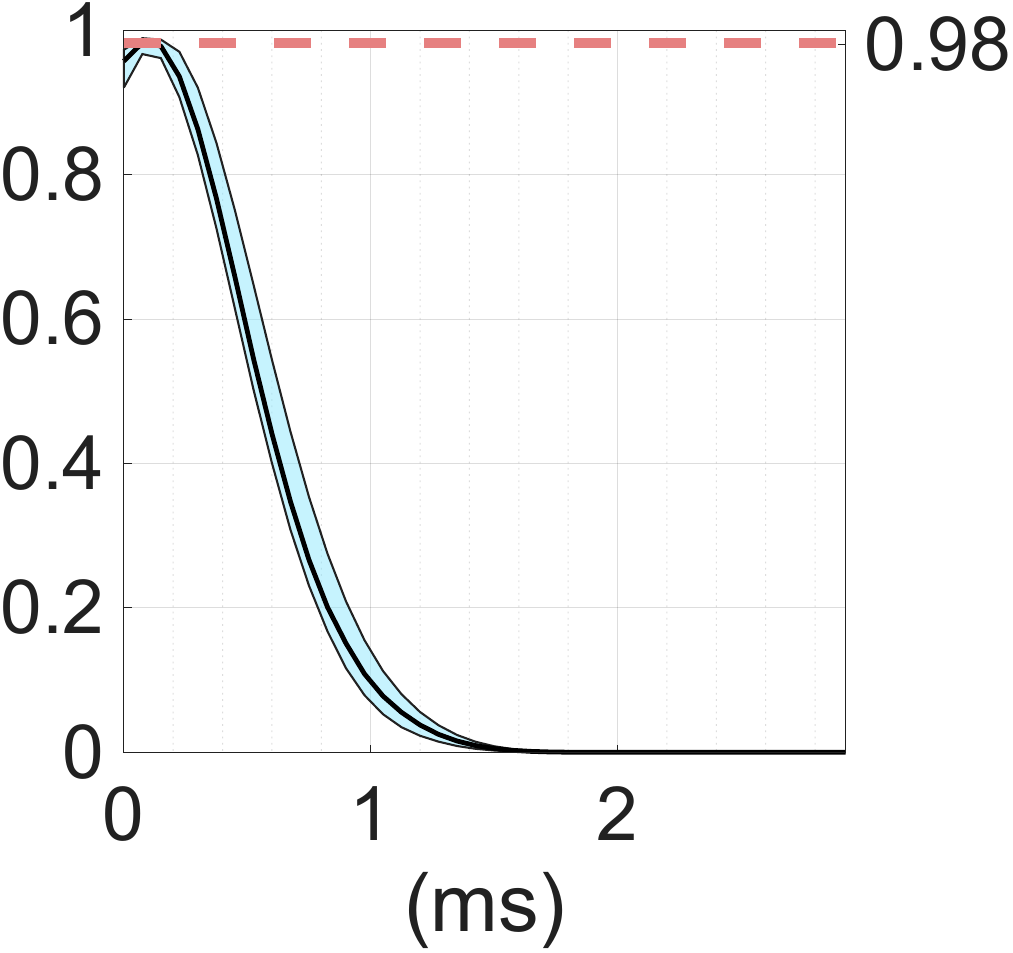}
\end{minipage}     
\begin{minipage}{2cm}
\centering
\includegraphics[trim={0 0cm 0 0},clip, width=1.7cm]{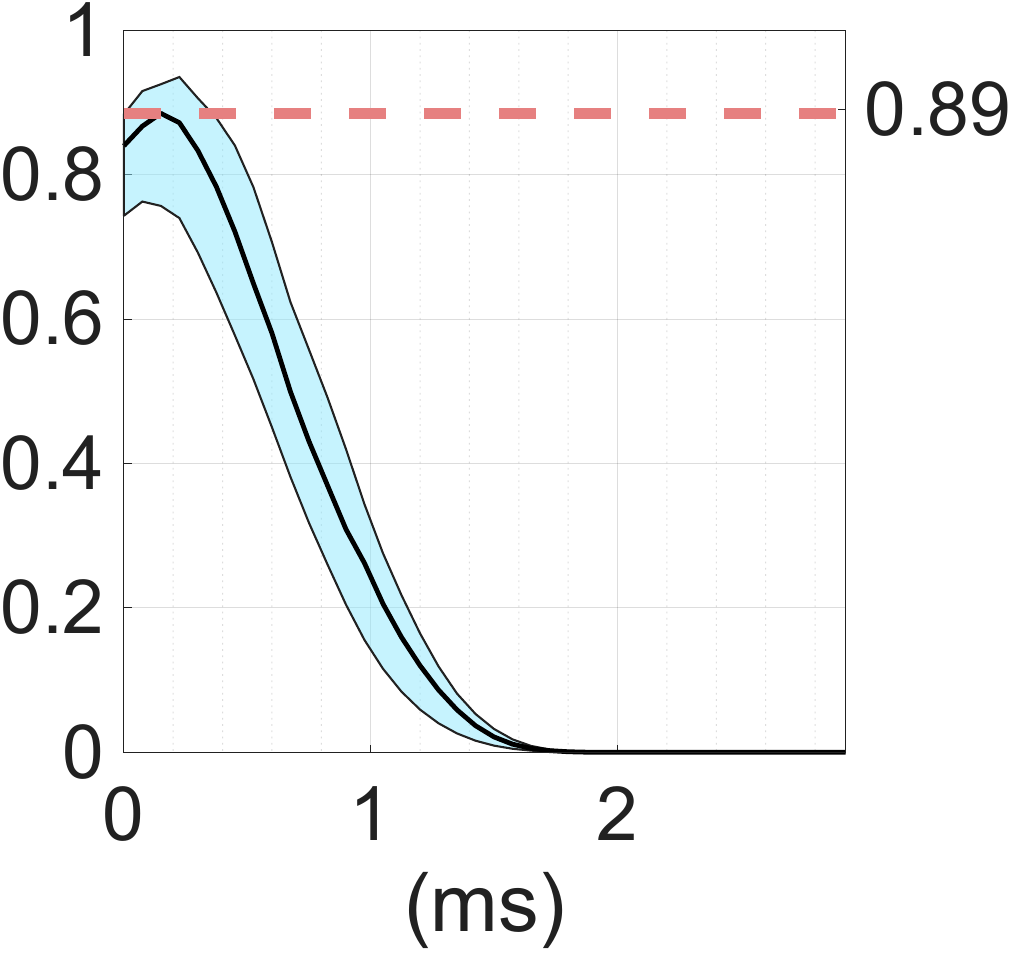}
\end{minipage}\\ 
\end{footnotesize}
    \caption{Normalized cross-correlations of the estimated tracks of the cortical activity (presented by blue curves in Figure \ref{fig:DynamicKFTracks}) computed against the true cortical activity evolution track shown in Figure \ref{fig:trueevolution}. The colored area covers the 10 and 90 \% interval around the median. A bold, solid curve presents the median track. The dashed light red horizontal line is the highest point of the median cross-correlation. This maximum value is presented numerically on the right side of the graphs.}
    \label{fig:DynamicKFCrossCorrTrue2}
\end{figure}

\begin{figure}[h!]
    \centering
    \includegraphics[width=0.25\linewidth]{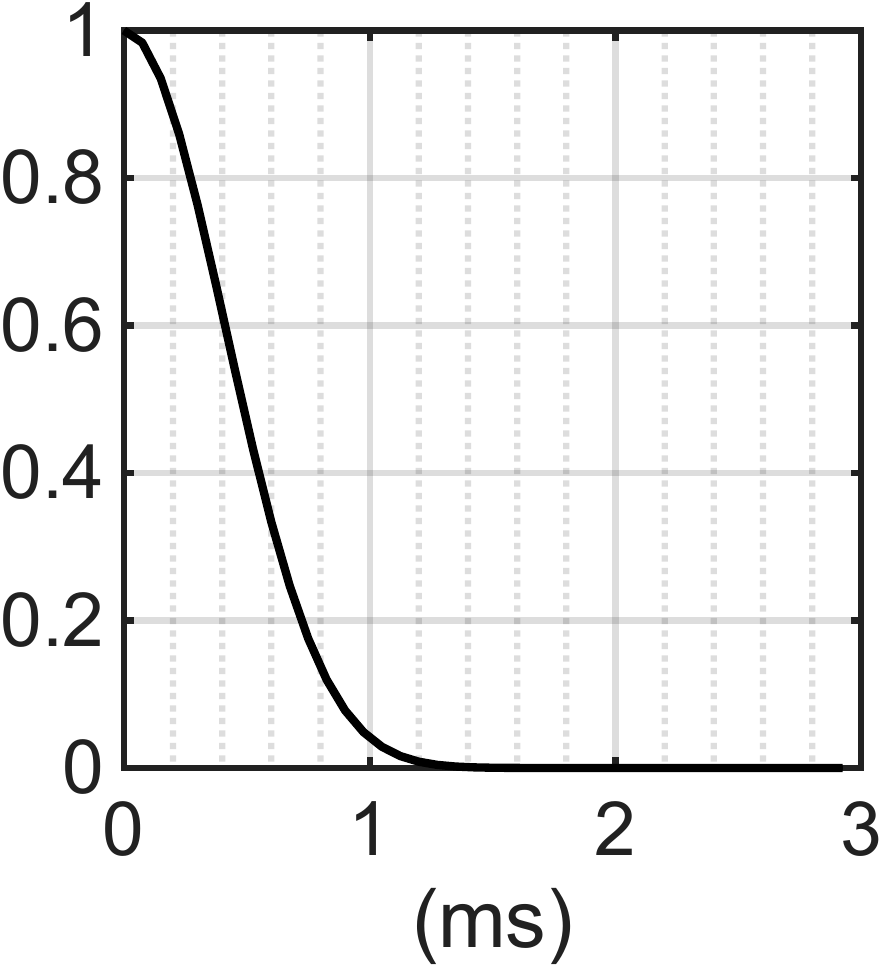}
    \caption{Ideal cross-correlation for Gaussian pulse time evolution against itself, i.e., autocorrelation. The results presented in Figures \ref{fig:DynamicKFCrossCorrTrue1} and \ref{fig:DynamicKFCrossCorrTrue2} are compared against this.}
    \label{fig:idealCrossCorr}
\end{figure}

\begin{figure*}[h!]
\centering
\begin{minipage}{0.7\textwidth}


   {\bf 3-DSKF} \\
\vskip0.2cm

\begin{minipage}{0.01\textwidth}
\rotatebox{90}{\small{30 dB}}
\end{minipage}\hspace{0.2cm}\begin{minipage}{0.13\textwidth}
    \centering
    \includegraphics[trim={0 0 1cm 0},clip,height=1.4cm]{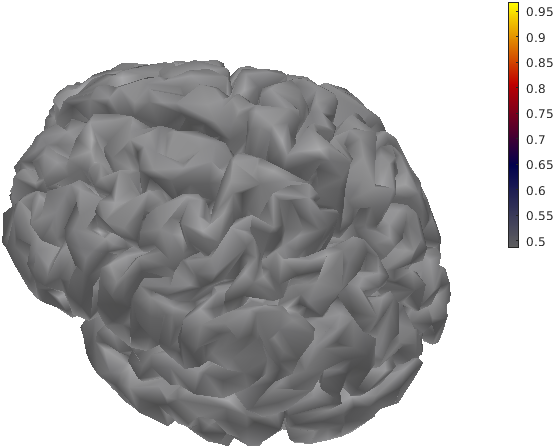}
\end{minipage}
\hspace{0.2cm}\begin{minipage}{0.13\textwidth}
    \centering
    \includegraphics[trim={0 0 1cm 0},clip,height=1.4cm]{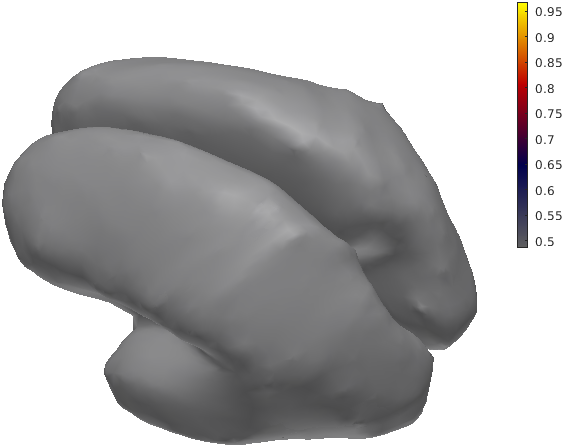}
\end{minipage}
\hspace{0.2cm}\begin{minipage}{0.1\textwidth}
    \centering
    \includegraphics[trim={0 0 3cm 0},clip,height=1.4cm]{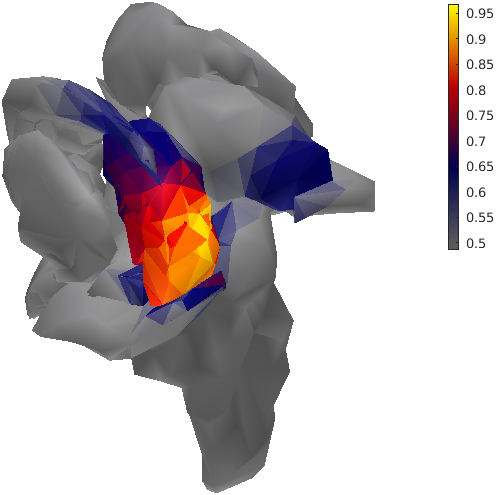}
\end{minipage}\hspace{1.3cm}\begin{minipage}{0.13\textwidth}
    \centering
    \includegraphics[trim={0 0 1cm 0},clip,height=1.4cm]{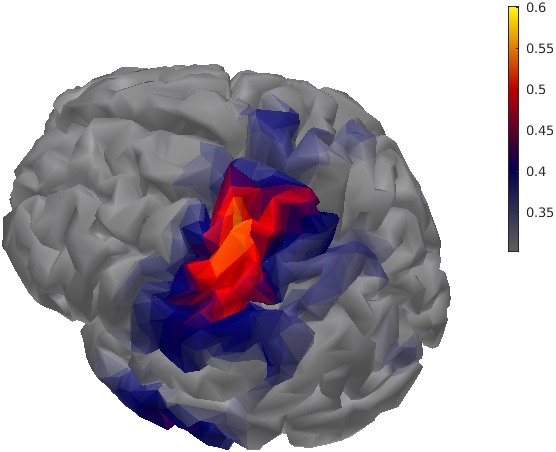}
\end{minipage}
\hspace{0.2cm}\begin{minipage}{0.13\textwidth}
    \centering
    \includegraphics[trim={0 0 1cm 0},clip,height=1.4cm]{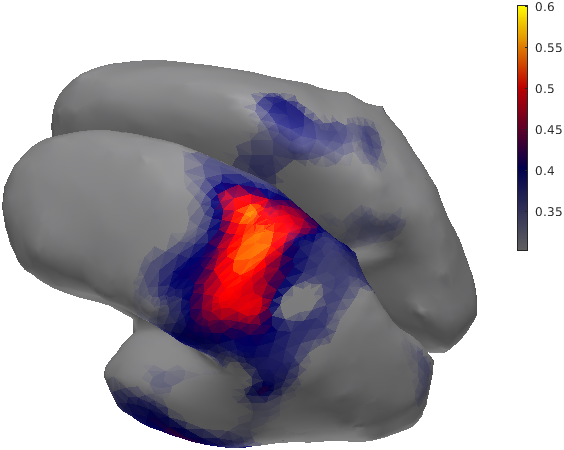}
\end{minipage}
\hspace{0.2cm}\begin{minipage}{0.1\textwidth}
    \centering
    \includegraphics[trim={0 0 3cm 0},clip,height=1.4cm]{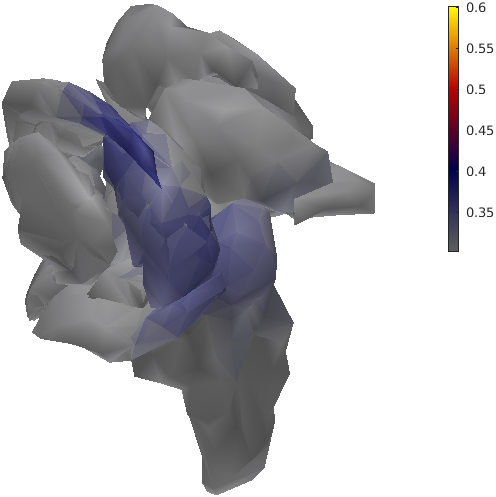}
\end{minipage}

\begin{minipage}{0.01\textwidth}
\rotatebox{90}{\small{20 dB}}
\end{minipage}\hspace{0.2cm}\begin{minipage}{0.13\textwidth}
    \centering
    \includegraphics[trim={0 0 1cm 0},clip,height=1.4cm]{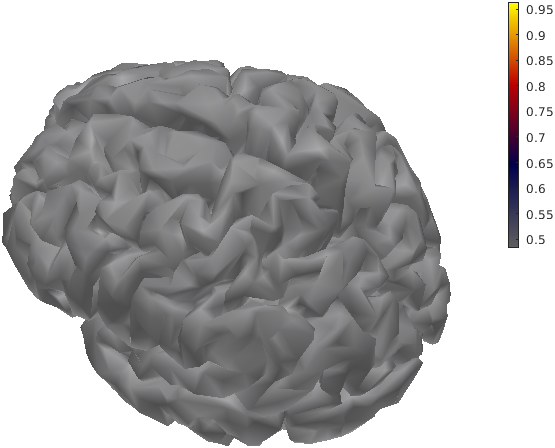}
\end{minipage}
\hspace{0.2cm}\begin{minipage}{0.13\textwidth}
    \centering
    \includegraphics[trim={0 0 1cm 0},clip,height=1.4cm]{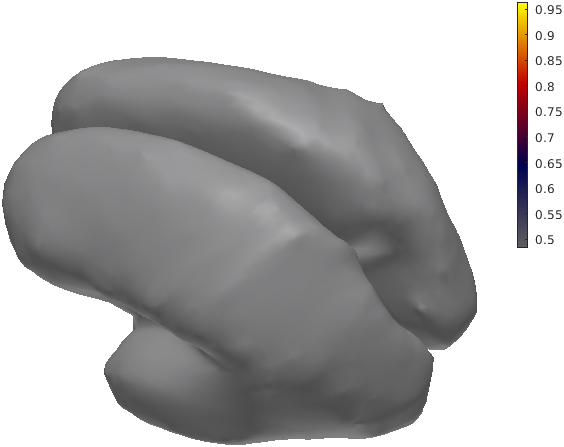}
\end{minipage}
\hspace{0.2cm}\begin{minipage}{0.1\textwidth}
    \centering
    \includegraphics[trim={0 0 3cm 0},clip,height=1.4cm]{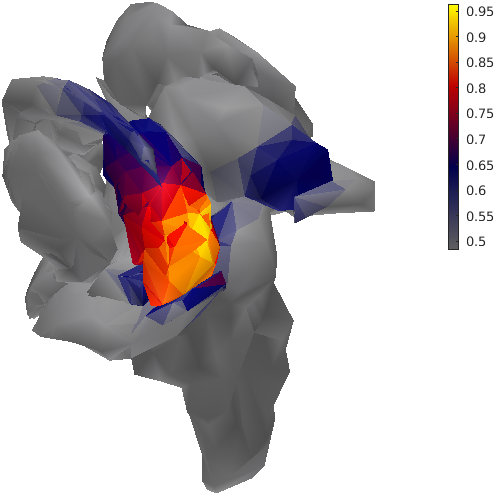}
\end{minipage}\hspace{1.3cm}\begin{minipage}{0.13\textwidth}
    \centering
    \includegraphics[trim={0 0 1cm 0},clip,height=1.4cm]{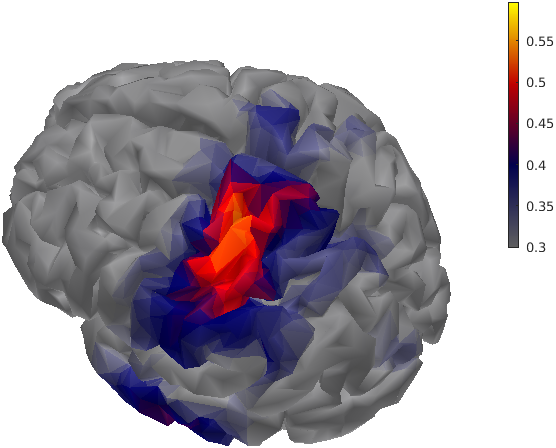}
\end{minipage}
\hspace{0.2cm}\begin{minipage}{0.13\textwidth}
    \centering
    \includegraphics[trim={0 0 1cm 0},clip,height=1.4cm]{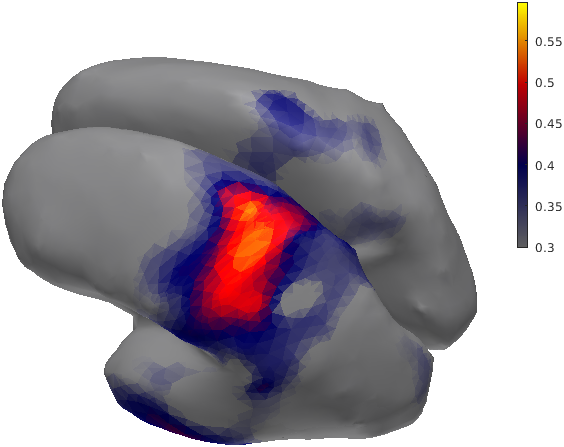}
\end{minipage}
\hspace{0.2cm}\begin{minipage}{0.1\textwidth}
    \centering
    \includegraphics[trim={0 0 3cm 0},clip,height=1.4cm]{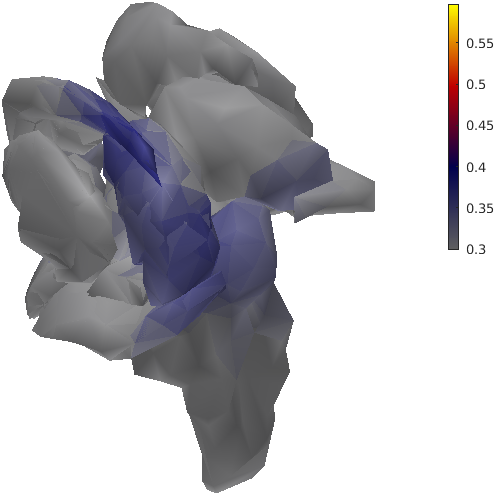}
\end{minipage}

\begin{minipage}{0.01\textwidth}
\rotatebox{90}{\small{10 dB}}
\end{minipage}\hspace{0.2cm}\begin{minipage}{0.13\textwidth}
    \centering
    \includegraphics[trim={0 0 1cm 0},clip,height=1.4cm]{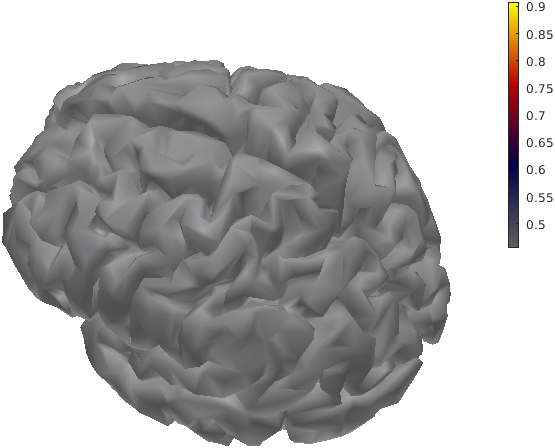}
\end{minipage}
\hspace{0.2cm}\begin{minipage}{0.13\textwidth}
    \centering
    \includegraphics[trim={0 0 1cm 0},clip,height=1.4cm]{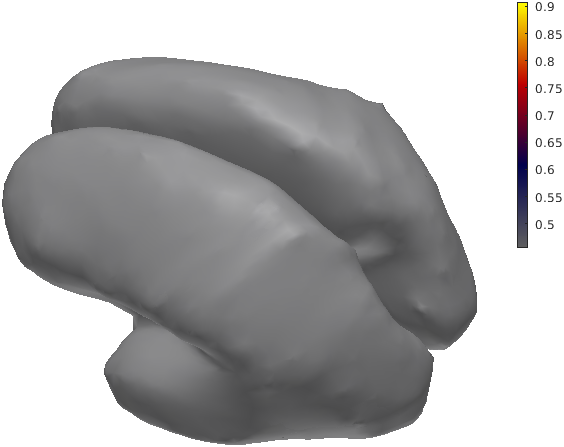}
\end{minipage}
\hspace{0.2cm}\begin{minipage}{0.1\textwidth}
    \centering
    \includegraphics[trim={0 0 3cm 0},clip,height=1.4cm]{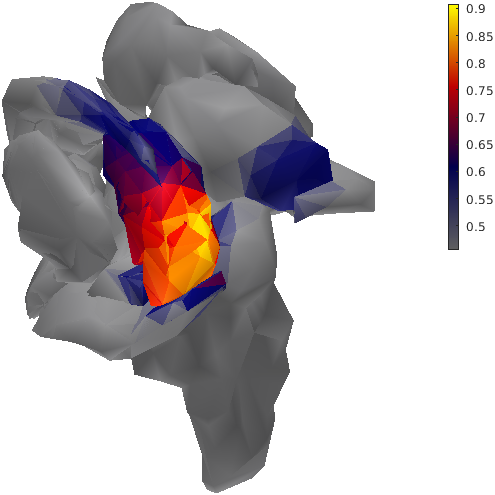}
\end{minipage}\hspace{1.3cm}\begin{minipage}{0.13\textwidth}
    \centering
    \includegraphics[trim={0 0 1cm 0},clip,height=1.4cm]{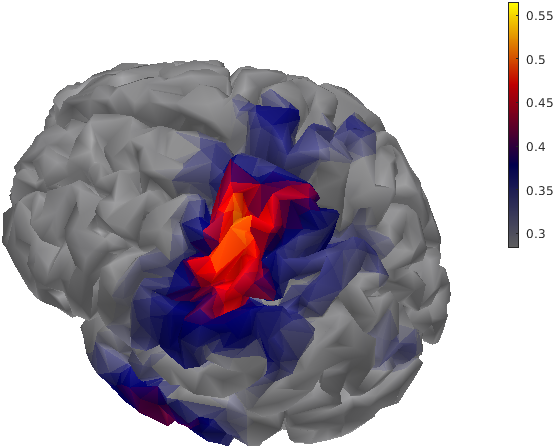}
\end{minipage}
\hspace{0.2cm}\begin{minipage}{0.13\textwidth}
    \centering
    \includegraphics[trim={0 0 1cm 0},clip,height=1.4cm]{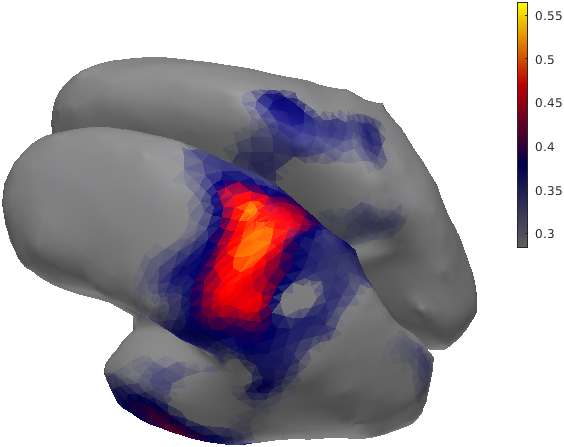}
\end{minipage}
\hspace{0.2cm}\begin{minipage}{0.1\textwidth}
    \centering
    \includegraphics[trim={0 0 3cm 0},clip,height=1.4cm]{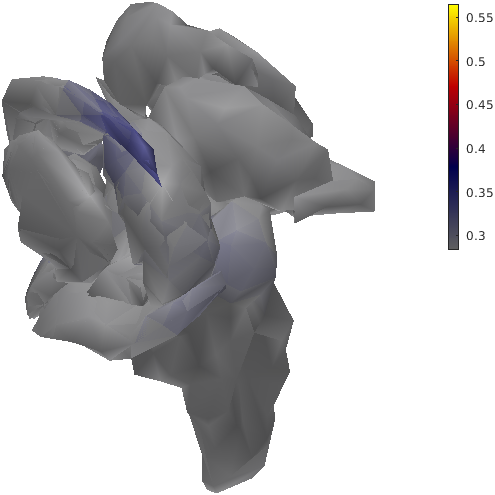}
\end{minipage}\\
\vskip0.2cm


   {\bf 2-DSKF} \\
\vskip0.2cm

\begin{minipage}{0.01\textwidth}
\rotatebox{90}{\small{30 dB}}
\end{minipage}\hspace{0.2cm}\begin{minipage}{0.13\textwidth}
    \centering
    \includegraphics[trim={0 0 1cm 0},clip,height=1.4cm]{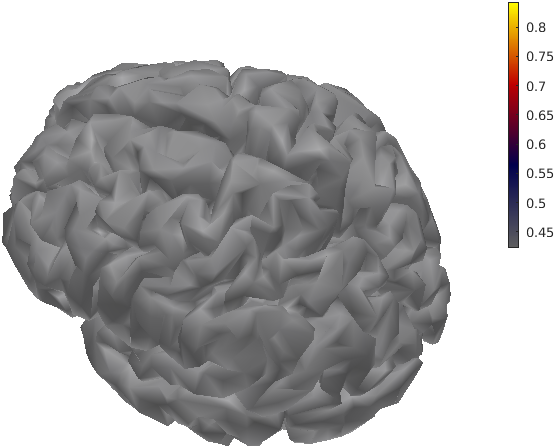}
\end{minipage}
\hspace{0.2cm}\begin{minipage}{0.13\textwidth}
    \centering
    \includegraphics[trim={0 0 1cm 0},clip,height=1.4cm]{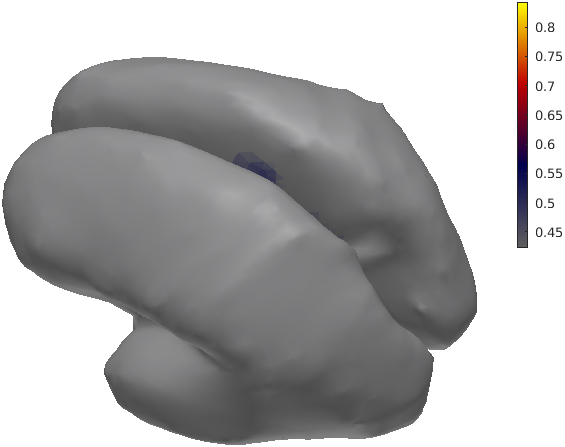}
\end{minipage}
\hspace{0.2cm}\begin{minipage}{0.1\textwidth}
    \centering
    \includegraphics[trim={0 0 3cm 0},clip,height=1.4cm]{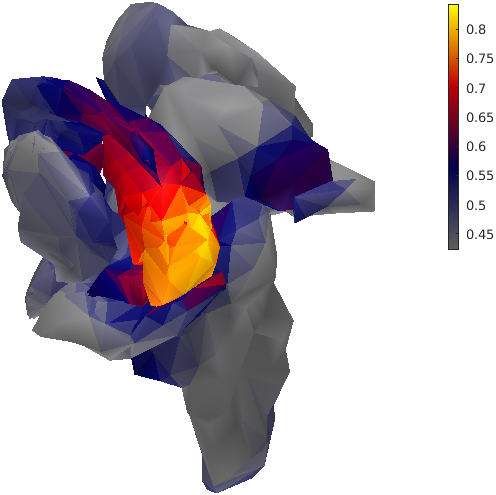}
\end{minipage}\hspace{1.3cm}\begin{minipage}{0.13\textwidth}
    \centering
    \includegraphics[trim={0 0 1cm 0},clip,height=1.4cm]{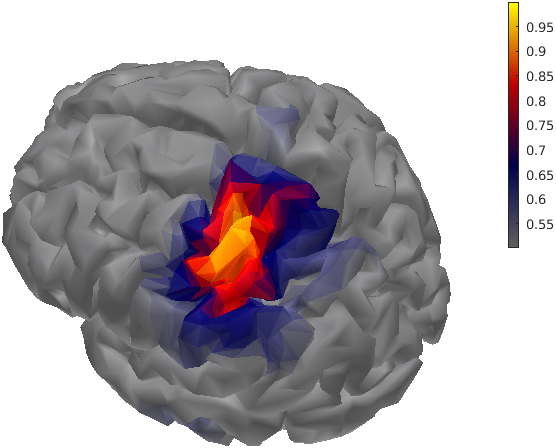}
\end{minipage}
\hspace{0.2cm}\begin{minipage}{0.13\textwidth}
    \centering
    \includegraphics[trim={0 0 1cm 0},clip,height=1.4cm]{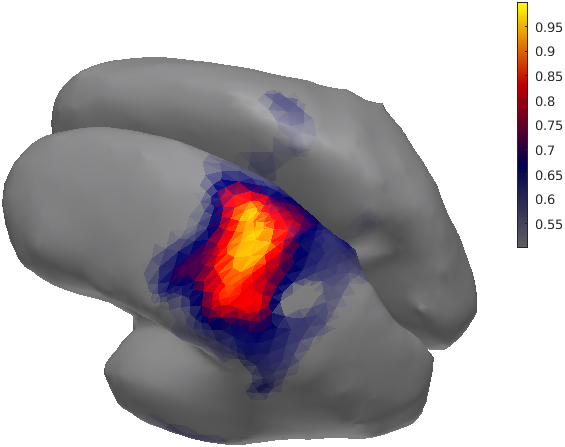}
\end{minipage}
\hspace{0.2cm}\begin{minipage}{0.1\textwidth}
    \centering
    \includegraphics[trim={0 0 3cm 0},clip,height=1.4cm]{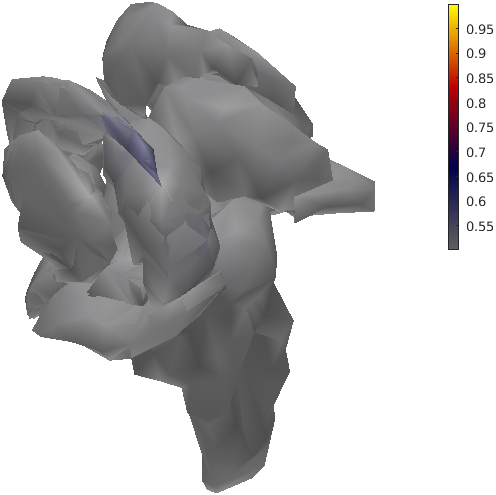}
\end{minipage}

\begin{minipage}{0.01\textwidth}
\rotatebox{90}{\small{20 dB}}
\end{minipage}\hspace{0.2cm}\begin{minipage}{0.13\textwidth}
    \centering
    \includegraphics[trim={0 0 1cm 0},clip,height=1.4cm]{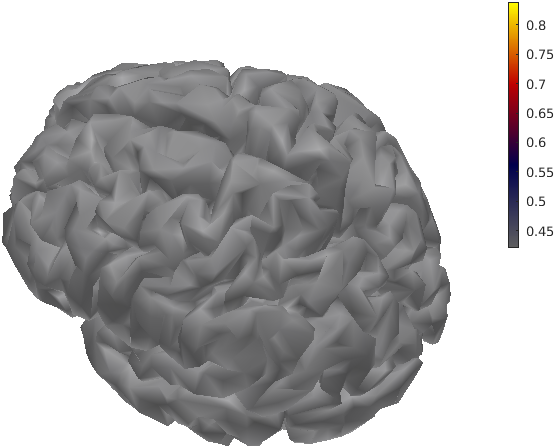}
\end{minipage}
\hspace{0.2cm}\begin{minipage}{0.13\textwidth}
    \centering
    \includegraphics[trim={0 0 1cm 0},clip,height=1.4cm]{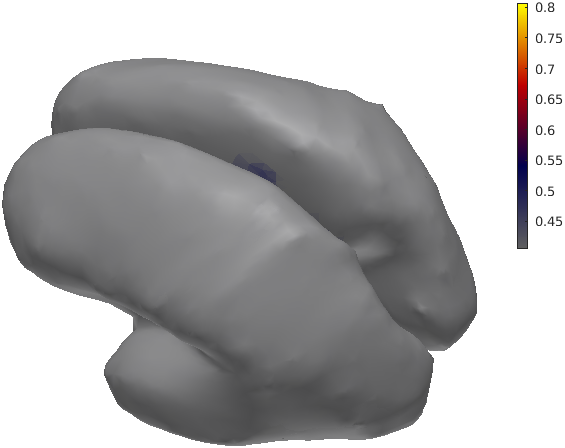}
\end{minipage}
\hspace{0.2cm}\begin{minipage}{0.1\textwidth}
    \centering
    \includegraphics[trim={0 0 3cm 0},clip,height=1.4cm]{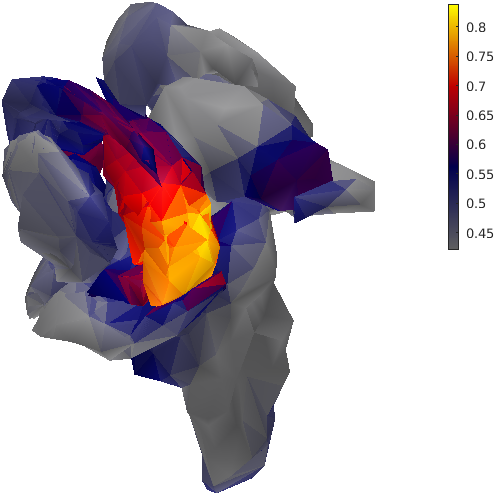}
\end{minipage}\hspace{1.3cm}\begin{minipage}{0.13\textwidth}
    \centering
    \includegraphics[trim={0 0 1cm 0},clip,height=1.4cm]{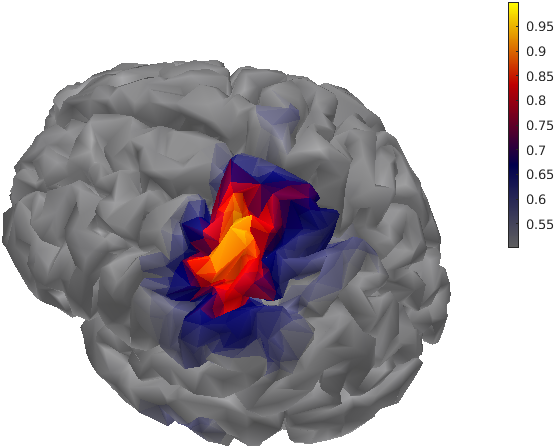}
\end{minipage}
\hspace{0.2cm}\begin{minipage}{0.13\textwidth}
    \centering
    \includegraphics[trim={0 0 1cm 0},clip,height=1.4cm]{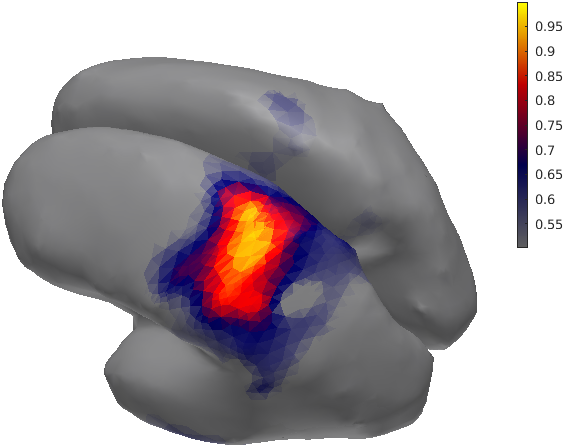}
\end{minipage}
\hspace{0.2cm}\begin{minipage}{0.1\textwidth}
    \centering
    \includegraphics[trim={0 0 3cm 0},clip,height=1.4cm]{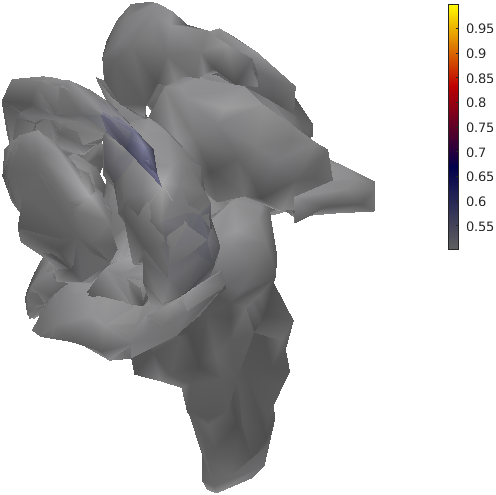}
\end{minipage}

\begin{minipage}{0.01\textwidth}
\rotatebox{90}{\small{10 dB}}
\end{minipage}\hspace{0.2cm}\begin{minipage}{0.13\textwidth}
    \centering
    \includegraphics[trim={0 0 1cm 0},clip,height=1.4cm]{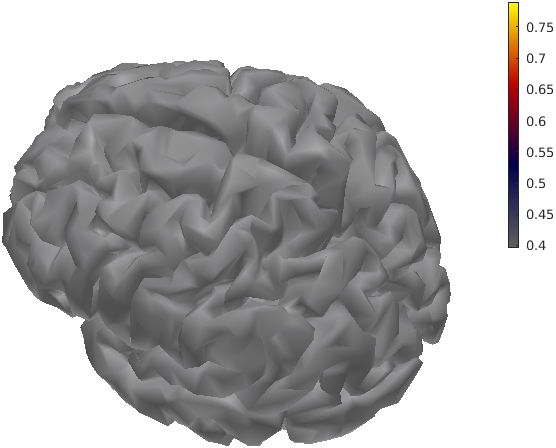}
\end{minipage}
\hspace{0.2cm}\begin{minipage}{0.13\textwidth}
    \centering
    \includegraphics[trim={0 0 1cm 0},clip,height=1.4cm]{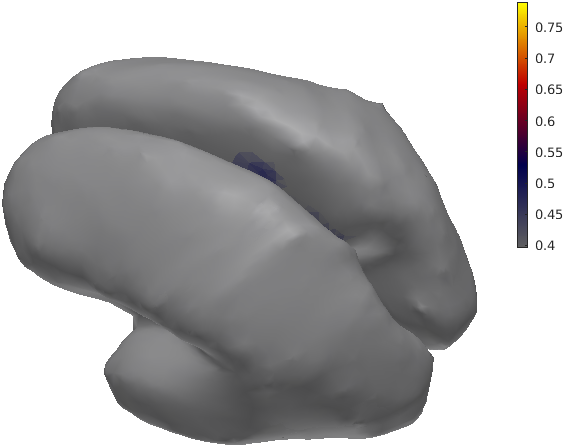}
\end{minipage}
\hspace{0.2cm}\begin{minipage}{0.1\textwidth}
    \centering
    \includegraphics[trim={0 0 3cm 0},clip,height=1.4cm]{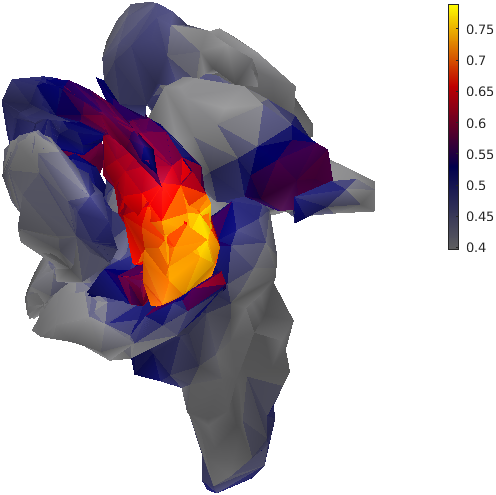}
\end{minipage}\hspace{1.3cm}\begin{minipage}{0.13\textwidth}
    \centering
    \includegraphics[trim={0 0 1cm 0},clip,height=1.4cm]{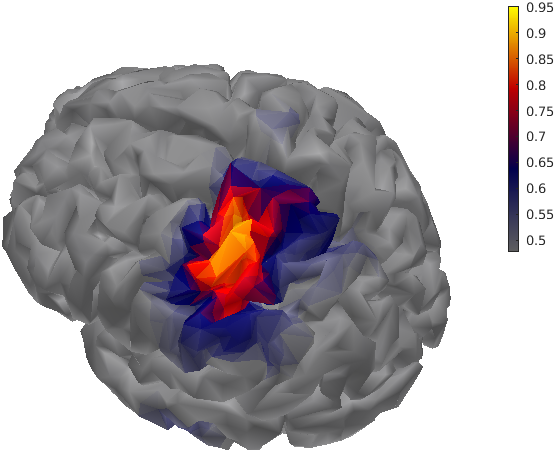}
\end{minipage}
\hspace{0.2cm}\begin{minipage}{0.13\textwidth}
    \centering
    \includegraphics[trim={0 0 1cm 0},clip,height=1.4cm]{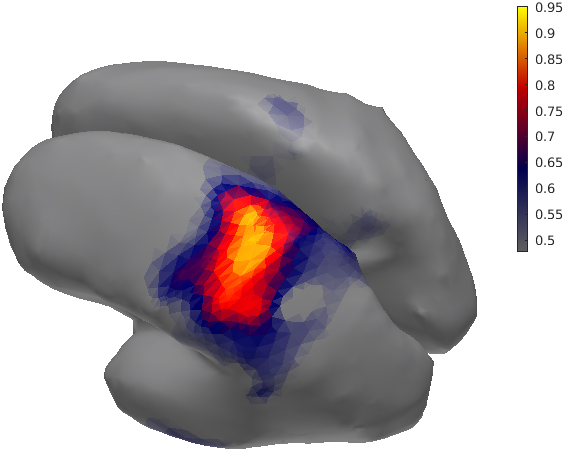}
\end{minipage}
\hspace{0.2cm}\begin{minipage}{0.1\textwidth}
    \centering
    \includegraphics[trim={0 0 3cm 0},clip,height=1.4cm]{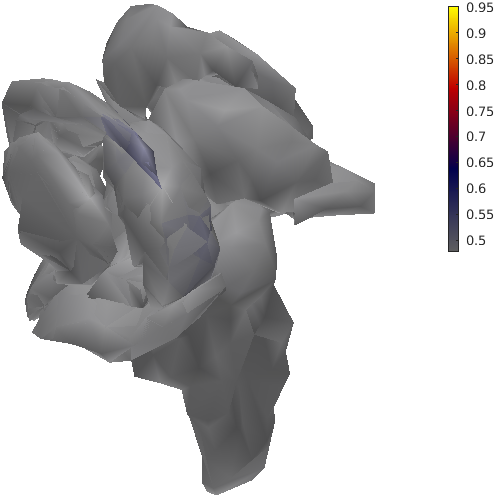}
\end{minipage}\\
\vskip0.2cm

    {\bf SSKF} \\
\vskip0.2cm

\begin{minipage}{0.01\textwidth}
\rotatebox{90}{\small{30 dB}}
\end{minipage}\hspace{0.2cm}\begin{minipage}{0.13\textwidth}
    \centering
    \includegraphics[trim={0 0 1cm 0},clip,height=1.4cm]{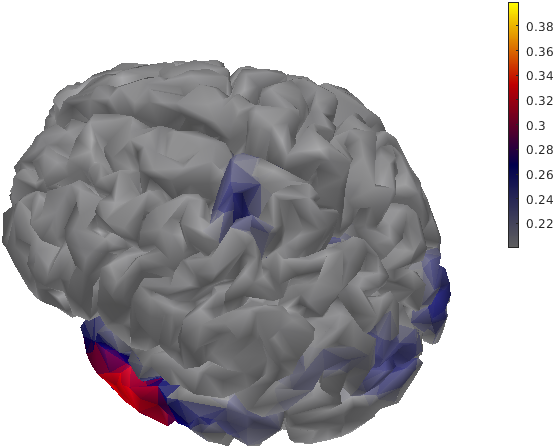}
\end{minipage}
\hspace{0.2cm}\begin{minipage}{0.13\textwidth}
    \centering
    \includegraphics[trim={0 0 1cm 0},clip,height=1.4cm]{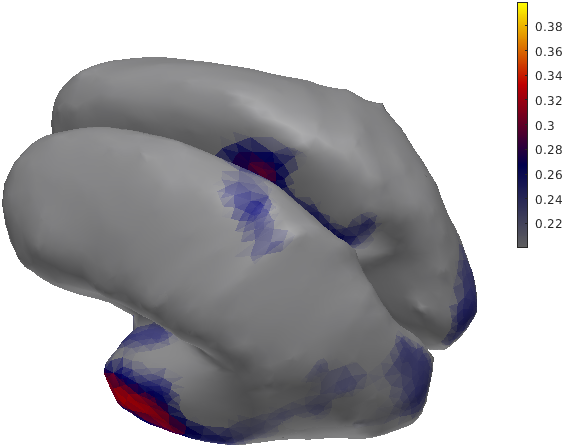}
\end{minipage}
\hspace{0.2cm}\begin{minipage}{0.1\textwidth}
    \centering
    \includegraphics[trim={0 0 3cm 0},clip,height=1.4cm]{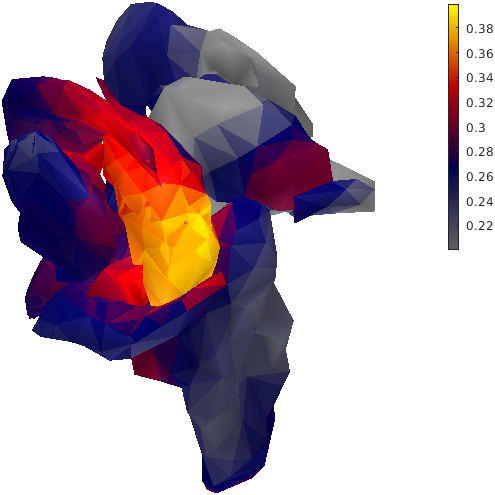}
\end{minipage}\hspace{1.3cm}\begin{minipage}{0.13\textwidth}
    \centering
    \includegraphics[trim={0 0 1cm 0},clip,height=1.4cm]{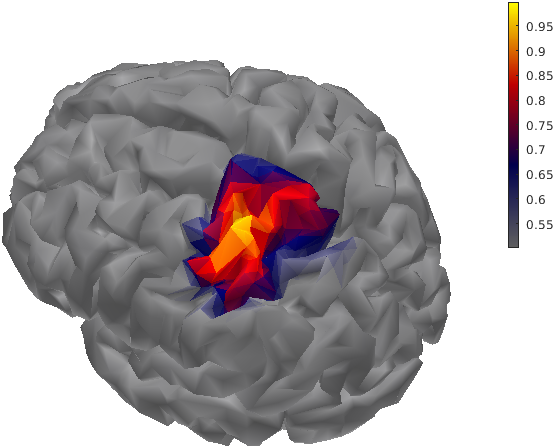}
\end{minipage}
\hspace{0.2cm}\begin{minipage}{0.13\textwidth}
    \centering
    \includegraphics[trim={0 0 1cm 0},clip,height=1.4cm]{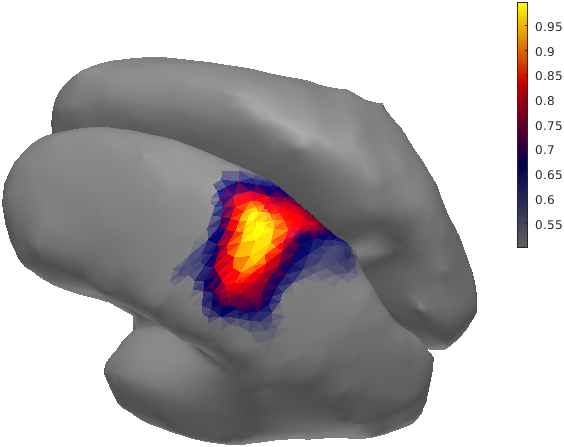}
\end{minipage}
\hspace{0.2cm}\begin{minipage}{0.1\textwidth}
    \centering
    \includegraphics[trim={0 0 3cm 0},clip,height=1.4cm]{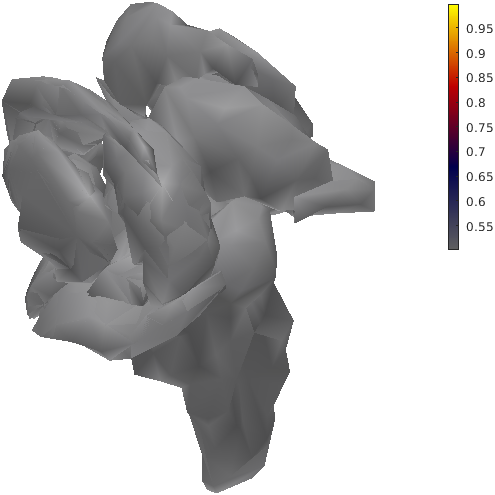}
\end{minipage}

\begin{minipage}{0.01\textwidth}
\rotatebox{90}{\small{20 dB}}
\end{minipage}\hspace{0.2cm}\begin{minipage}{0.13\textwidth}
    \centering
    \includegraphics[trim={0 0 1cm 0},clip,height=1.4cm]{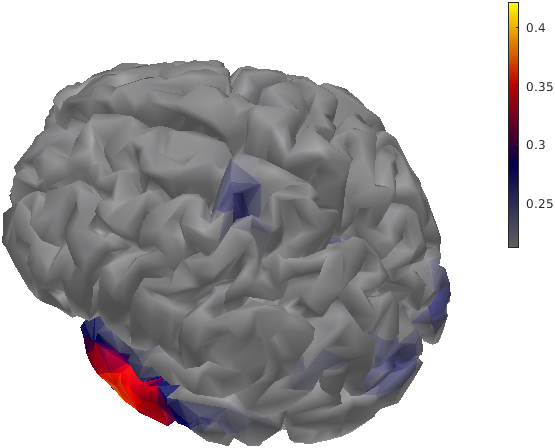}
\end{minipage}
\hspace{0.2cm}\begin{minipage}{0.13\textwidth}
    \centering
    \includegraphics[trim={0 0 1cm 0},clip,height=1.4cm]{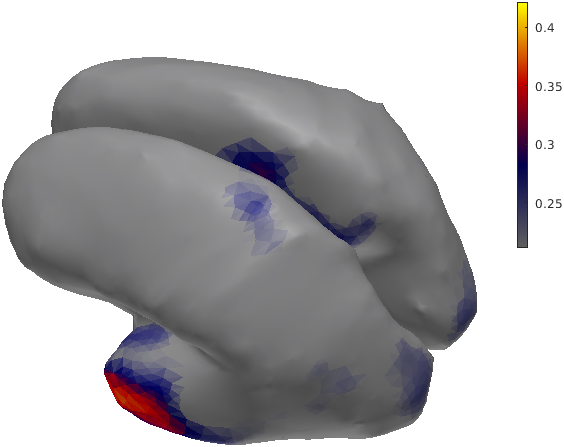}
\end{minipage}
\hspace{0.2cm}\begin{minipage}{0.1\textwidth}
    \centering
    \includegraphics[trim={0 0 3cm 0},clip,height=1.4cm]{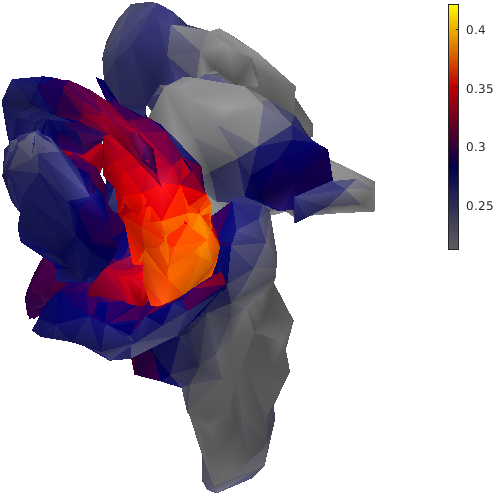}
\end{minipage}\hspace{1.3cm}\begin{minipage}{0.13\textwidth}
    \centering
    \includegraphics[trim={0 0 1cm 0},clip,height=1.4cm]{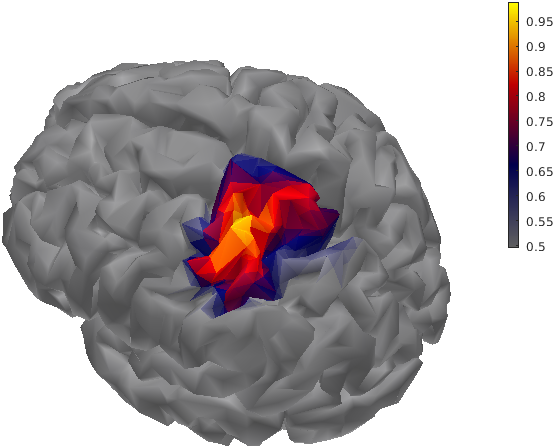}
\end{minipage}
\hspace{0.2cm}\begin{minipage}{0.13\textwidth}
    \centering
    \includegraphics[trim={0 0 1cm 0},clip,height=1.4cm]{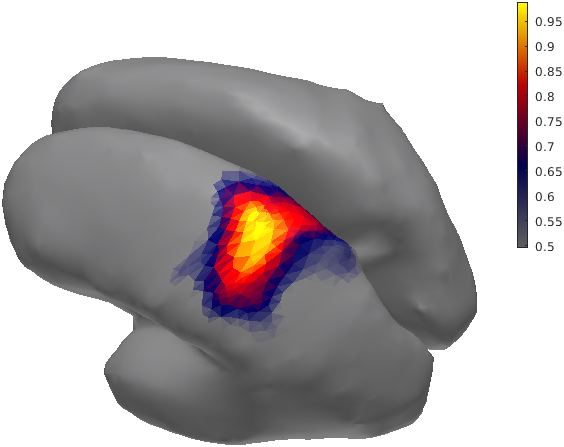}
\end{minipage}
\hspace{0.2cm}\begin{minipage}{0.1\textwidth}
    \centering
    \includegraphics[trim={0 0 3cm 0},clip,height=1.4cm]{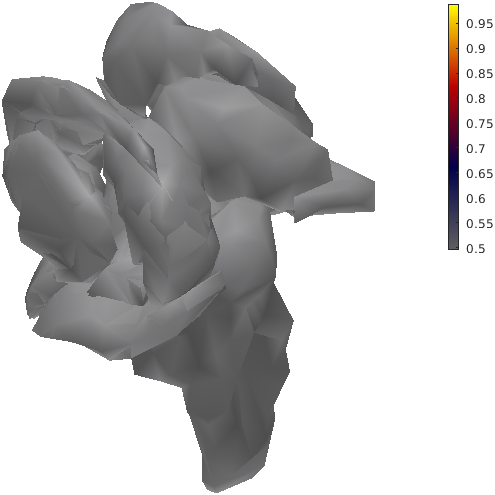}
\end{minipage}

\begin{minipage}{0.01\textwidth}
\rotatebox{90}{\small{10 dB}}
\end{minipage}\hspace{0.2cm}\begin{minipage}{0.13\textwidth}
    \centering
    \includegraphics[trim={0 0 1cm 0},clip,height=1.4cm]{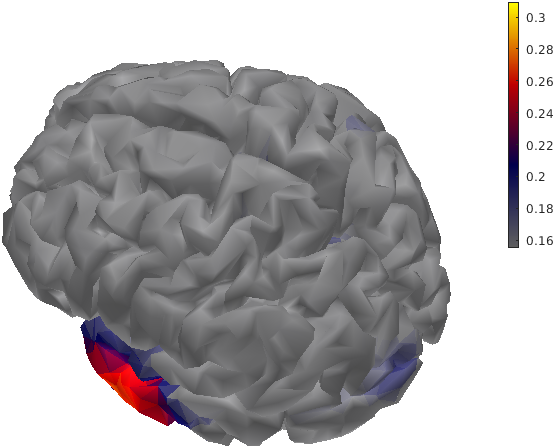}
\end{minipage}
\hspace{0.2cm}\begin{minipage}{0.13\textwidth}
    \centering
    \includegraphics[trim={0 0 1cm 0},clip,height=1.4cm]{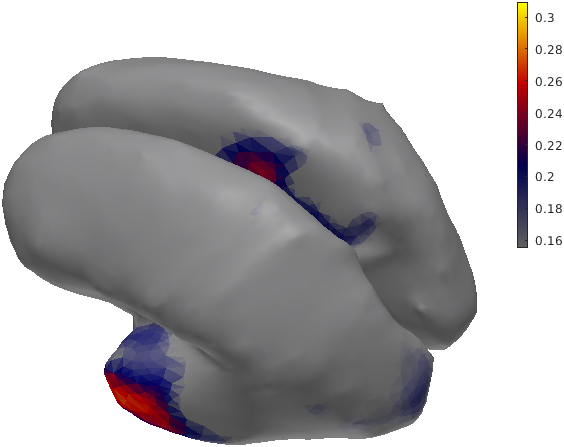}
\end{minipage}
\hspace{0.2cm}\begin{minipage}{0.1\textwidth}
    \centering
    \includegraphics[trim={0 0 3cm 0},clip,height=1.4cm]{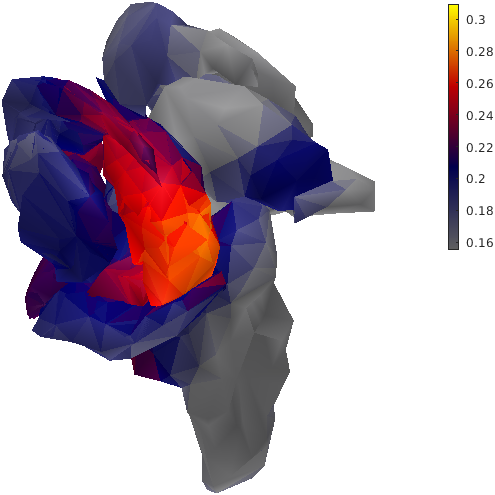}
\end{minipage} \hspace{1.3cm}\begin{minipage}{0.13\textwidth}
    \centering
    \includegraphics[trim={0 0 1cm 0},clip,height=1.4cm]{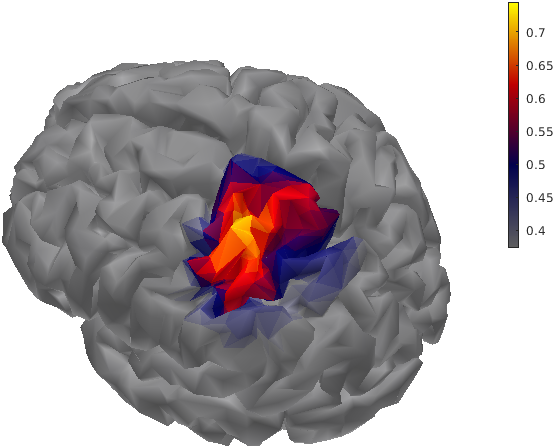}
\end{minipage}
\hspace{0.2cm}\begin{minipage}{0.13\textwidth}
    \centering
    \includegraphics[trim={0 0 1cm 0},clip,height=1.4cm]{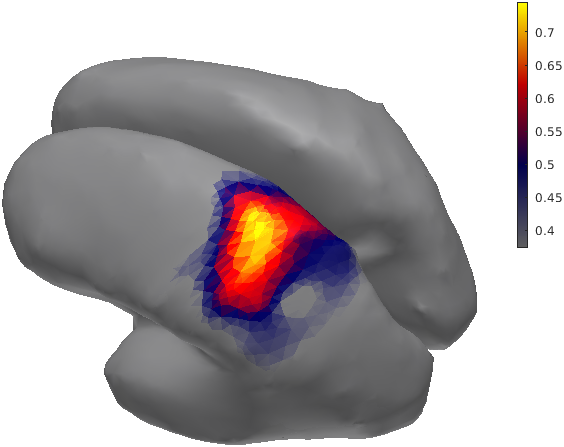}
\end{minipage}
\hspace{0.2cm}\begin{minipage}{0.1\textwidth}
    \centering
    \includegraphics[trim={0 0 3cm 0},clip,height=1.4cm]{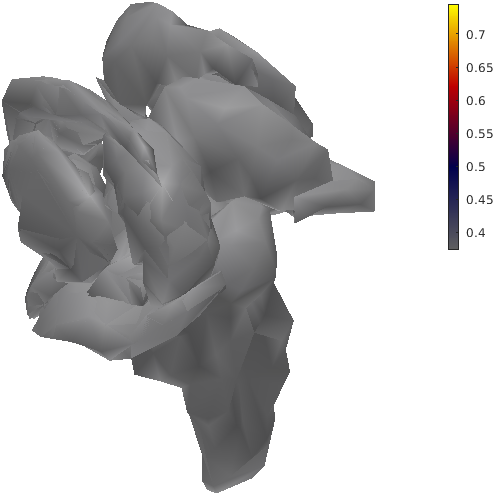}
\end{minipage}\\
\vskip0.2cm


   {\bf SKF} \\
\vskip0.2cm

\begin{minipage}{0.01\textwidth}
\rotatebox{90}{\small{30 dB}}
\end{minipage}\hspace{0.2cm}\begin{minipage}{0.13\textwidth}
    \centering
    \includegraphics[trim={0 0 1cm 0},clip,height=1.4cm]{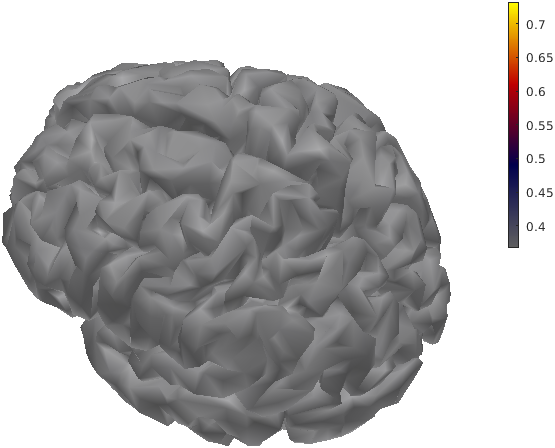}
\end{minipage}
\hspace{0.2cm}\begin{minipage}{0.13\textwidth}
    \centering
    \includegraphics[trim={0 0 1cm 0},clip,height=1.4cm]{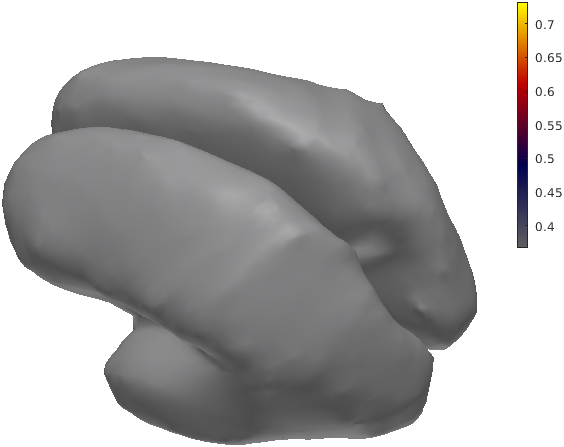}
\end{minipage}
\hspace{0.2cm}\begin{minipage}{0.1\textwidth}
    \centering
    \includegraphics[trim={0 0 3cm 0},clip,height=1.4cm]{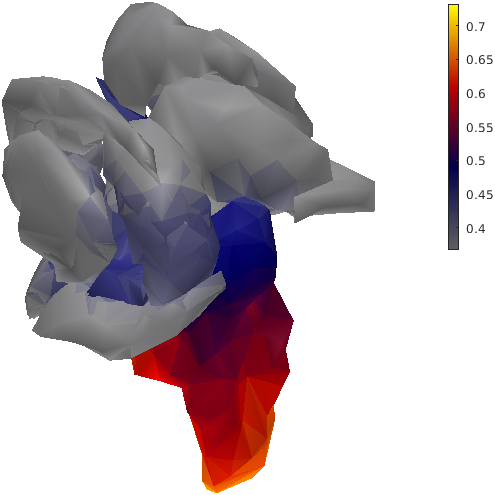}
\end{minipage}\hspace{1.3cm}\begin{minipage}{0.13\textwidth}
    \centering
    \includegraphics[trim={0 0 1cm 0},clip,height=1.4cm]{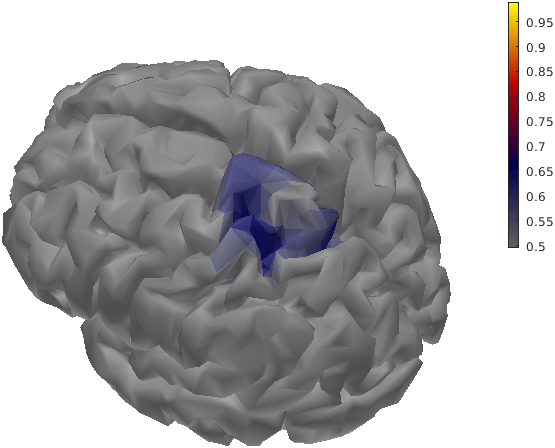}
\end{minipage}
\hspace{0.2cm}\begin{minipage}{0.13\textwidth}
    \centering
    \includegraphics[trim={0 0 1cm 0},clip,height=1.4cm]{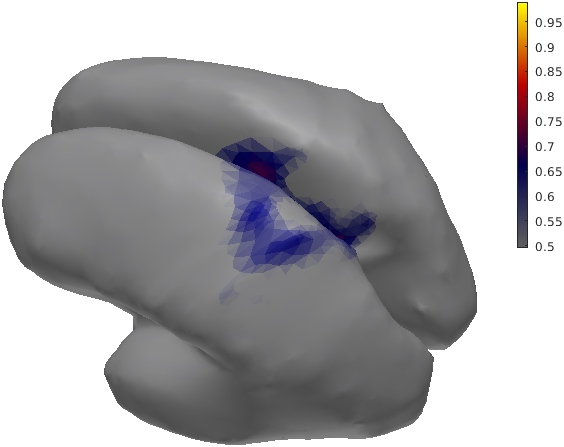}
\end{minipage}
\hspace{0.2cm}\begin{minipage}{0.1\textwidth}
    \centering
    \includegraphics[trim={0 0 3cm 0},clip,height=1.4cm]{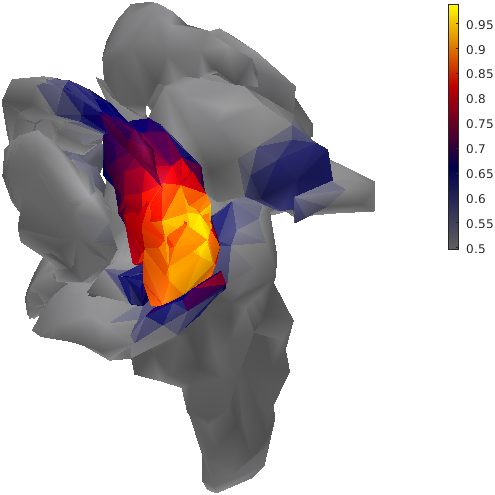}
\end{minipage}

\begin{minipage}{0.01\textwidth}
\rotatebox{90}{\small{20 dB}}
\end{minipage}\hspace{0.2cm}\begin{minipage}{0.13\textwidth}
    \centering
    \includegraphics[trim={0 0 1cm 0},clip,height=1.4cm]{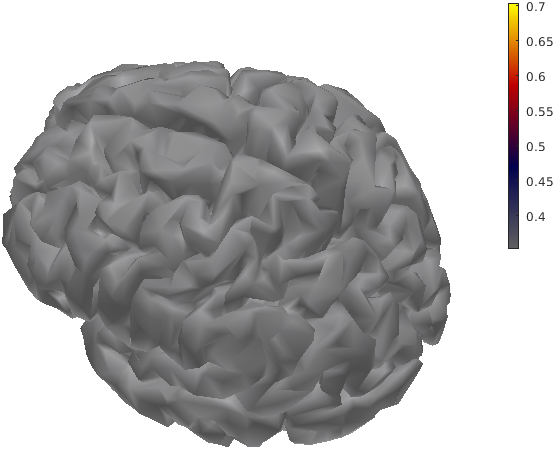}
\end{minipage}
\hspace{0.2cm}\begin{minipage}{0.13\textwidth}
    \centering
    \includegraphics[trim={0 0 1cm 0},clip,height=1.4cm]{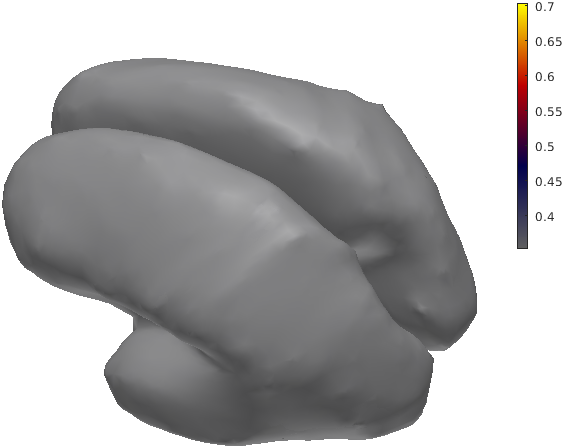}
\end{minipage}
\hspace{0.2cm}\begin{minipage}{0.1\textwidth}
    \centering
    \includegraphics[trim={0 0 3cm 0},clip,height=1.4cm]{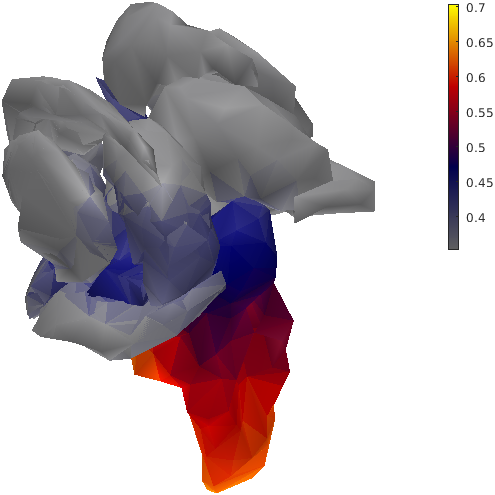}
\end{minipage}\hspace{1.3cm}\begin{minipage}{0.13\textwidth}
    \centering
    \includegraphics[trim={0 0 1cm 0},clip,height=1.4cm]{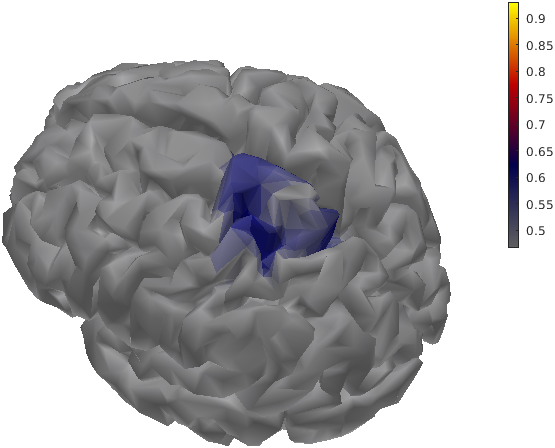}
\end{minipage}
\hspace{0.2cm}\begin{minipage}{0.13\textwidth}
    \centering
    \includegraphics[trim={0 0 1cm 0},clip,height=1.4cm]{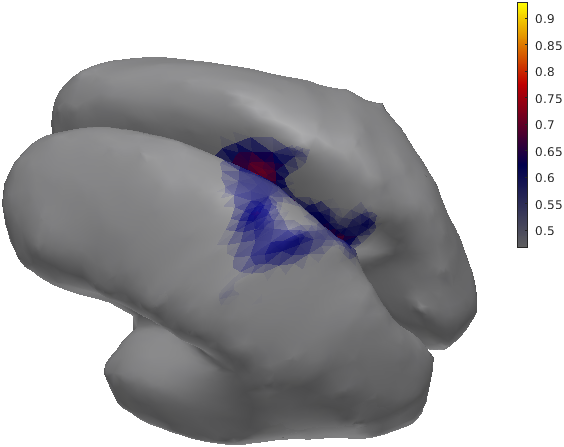}
\end{minipage}
\hspace{0.2cm}\begin{minipage}{0.1\textwidth}
    \centering
    \includegraphics[trim={0 0 3cm 0},clip,height=1.4cm]{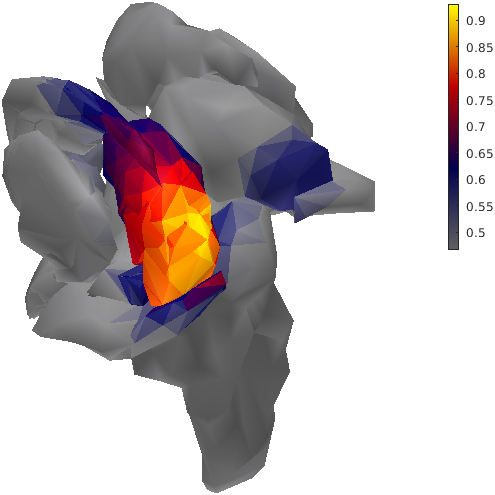}
\end{minipage}

\begin{minipage}{0.01\textwidth}
\rotatebox{90}{\small{10 dB}}
\end{minipage}\hspace{0.2cm}\begin{minipage}{0.13\textwidth}
    \centering
    \includegraphics[trim={0 0 1cm 0},clip,height=1.4cm]{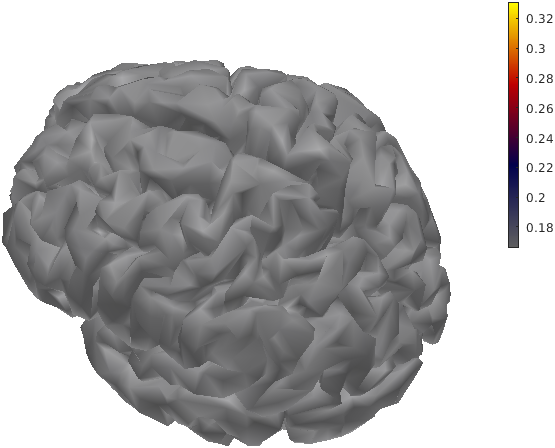}
\end{minipage}
\hspace{0.2cm}\begin{minipage}{0.13\textwidth}
    \centering
    \includegraphics[trim={0 0 1cm 0},clip,height=1.4cm]{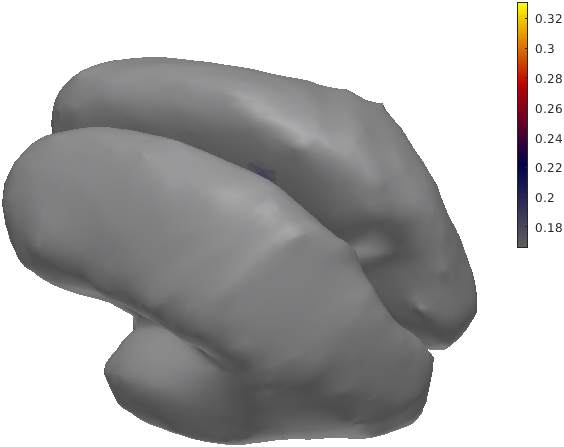}
\end{minipage}
\hspace{0.2cm}\begin{minipage}{0.1\textwidth}
    \centering
    \includegraphics[trim={0 0 3cm 0},clip,height=1.4cm]{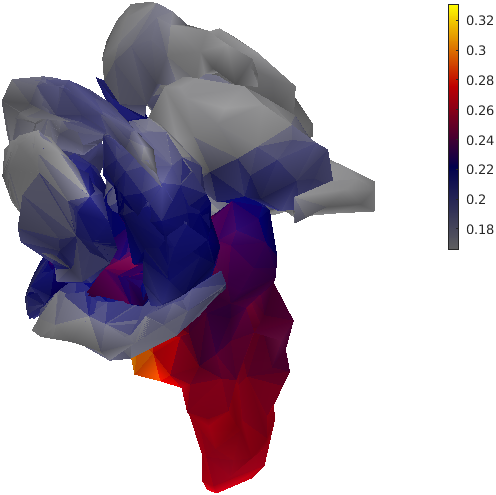}
\end{minipage}\hspace{1.3cm}\begin{minipage}{0.13\textwidth}
    \centering
    \includegraphics[trim={0 0 1cm 0},clip,height=1.4cm]{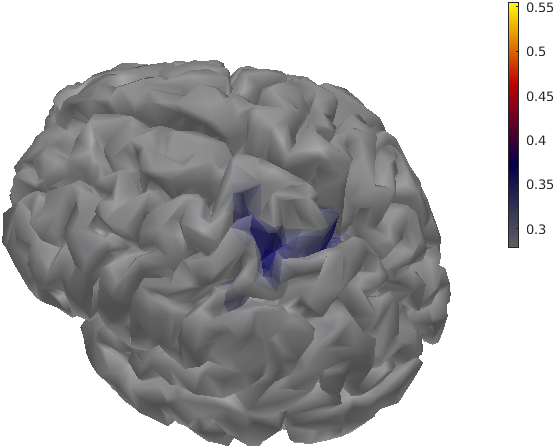}
\end{minipage}
\hspace{0.2cm}\begin{minipage}{0.13\textwidth}
    \centering
    \includegraphics[trim={0 0 1cm 0},clip,height=1.4cm]{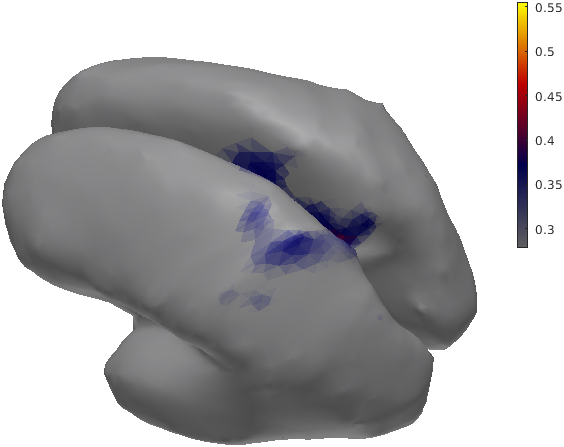}
\end{minipage}
\hspace{0.2cm}\begin{minipage}{0.1\textwidth}
    \centering
    \includegraphics[trim={0 0 3cm 0},clip,height=1.4cm]{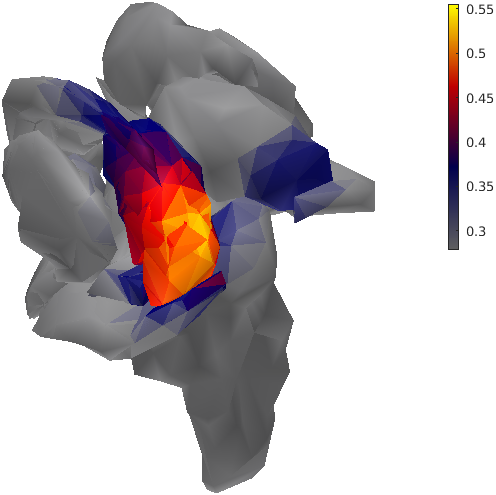}
\end{minipage} \\
\end{minipage}
\begin{minipage}{0.01\textwidth}
    \centering
    \includegraphics[trim={0 0 0 0},clip]{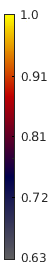}
\end{minipage}
\vskip0.2cm

    \caption{Spatial estimations of the activity at both peaks with Dynamical Standardized Kalman filter and the original Standardized Kalman filter (on top) and Standardized Kalman filter solution smoothed by RTS smoother (on bottom). The first three columns (starting from left to right) of Figure \ref{fig:dynamicKFBrains} show the average reconstructions at a time point 1.1 ms, and the last three columns show the reconstructions at 1.9 ms. The first and fourth columns show the cortical surface, the second and fifth columns show the inflated version of the cortical surface, and the third and sixth columns show the deep brain structures, including the thalamus and brainstem. The rows represented the estimation at the different noise levels: 30, 20, and 10 dB. The brain activity distributions have been thresholded to the relative amplitude interval from 0.63  to 1 (-4--0 dB) with respect to the maximum amplitude of the spatio-temporal reconstruction, i.e., the dynamical range of each distribution is 4 dB.}
    \label{fig:dynamicKFBrains}
\end{figure*}

\section{Simulations}

The numerical simulations of this study were conducted using the open Zeffiro Interface software package, in particular, its August2025 release\footnote{\url{https://doi.org/10.5281/zenodo.16839914}}, which includes open implementations of the forward and inverse methodology applied in this study. The scripts for Kalman filter can be found in the folder with relative path \texttt{/plugins/Kalman} w.r.t.\ the root directory. We used a realistic open head model obtained from the T1-weighted MRI data of a 49-year-old healthy male  \cite{piastra_2020_3888381}. The multi-compartment head mesh used for the simulations was constructed using Zeffiro's meshing tool as described \cite{galaz2023multi} with the following tissue conductivities \cite{dannhauer2011modeling}: 0.14 S/m for the white matter, 1.79 S/m for the CSF, 0.0064 S/m for the skull, 0.33 S/m for the skin, gray matter, and the remaining tissues. The lead field matrix $L$ was built using linear basis functions via Zeffiro's lead field generator. A set of 515 point sources was placed randomly on the gray matter layer in the neocortex and deep structures below the white matter layer.

Two dipolar sources were used to simulate short-latency somatosensory evoked potentials (SEP). This choice is based on the fact that this source setup has a known and predictable activity pattern \cite{allison1991cortical} and has been used before in \cite{RezaeiA2021,Lahtinen2022}, making it well-suited for validation studies. In particular, the cortical source was oriented along the normal of the sulcus wall, while the deep source was oriented towards the scalp surface. 
Figure \ref{fig:sourcesetup} (top row) depicts the locations and orientations of the left thalamic source and cortical source at the somatosensory area. Both sources (cortical and thalamic) had a peak amplitude of 10 nAm, which is a typical activity strength obtained by EEG according to \cite{hamalainen1993}. The time evolution of both simulated sources followed a Gaussian pulse with a length of 2 milliseconds (ms). More precisely, the first source (thalamic) has peak activity at 1.1 ms, and the cortical one has peak activity at 1.9 ms. The total duration of the simulation of the activity was 3 ms, recorded with a sampling frequency of 20000 Hz, and 10, 20, or 30 dB white Gaussian noise was added to the EEG measurements. The time interval \(\Delta t\) used in the filtering was \(3\text{ms}/40\), dividing the time sequence into 40 time steps.
Furthermore, we generated 20 different additive measurement noise realizations to produce the sets of EEG data. An example of these measurement data without noise and with these noise levels is presented in Figure \ref{fig:butterfly}. The true source strength curves, that we aim to track, are presented in Figure \ref{fig:trueevolution}.

Then, we calculated the estimations over the 3 ms time course using the Standardized Kalman filter (SKF), the Standardized Kalman filter with RTS smoother (SSKF), and the SKF with two different kinematic evolution models that are the "activity-velocity" Dynamical SKF (2-DSKF) and "activity-velocity-acceleration" Dynamical SKF (3-DSKF); both proposed in the current work.
The estimated source strengths were averaged within the regions of interest shown in the second row of Figure \ref{fig:sourcesetup} at each time frame. This way, we could get a temporal profile of the reconstructed source activity (strength) over the whole 3 ms time interval, both for the cortical and thalamic regions. To assess the statistical characteristics of these profiles, we used the reconstructions from the EEG sets with the 20 different noise realizations. Specifically, we computed the mean of the estimated averages for each region at each time frame, along with the first quartile interval around the mean. (This is illustrated in Figure~\ref{fig:DynamicKFTracks}). To quantitatively assess the correctness of the activity evolution curves, we computed normalized cross-correlations, i.e.,
\begin{equation}
    (f\star g)[s]=\sum_{n=0}^\infty \frac{f[n]g[n+s]}{\left\|f\right\|\, \left\|g\right\|},
\end{equation}
for thalamic activity against cortical activity plotted in Figure \ref{fig:DynamicKFCrossCorr}, and the estimated thalamic evolution against the true thalamic evolution in Figure \ref{fig:DynamicKFCrossCorrTrue1} and the estimated cortical evolution against the true cortical evolution in Figure \ref{fig:DynamicKFCrossCorrTrue2}. As the evolution curves for cortical and deep activity are identical Gaussian pulses, both of them have the same ideal cross-correlation curve presented in Figure \ref{fig:idealCrossCorr}. So, the more truthful the estimated track is, the more its cross-correlation against the true track represents the ideal curve. Therefore, by requiring independent tracks for both activities, we can compute the following $L2$-error measure for the tracks:
\begin{equation}\label{eq:correrr}
    \sqrt{\sum_{n=0}^N \left[\left(\Delta f_{\mathrm{Thalamic}}[n]\right)^2+\left(\Delta f_{\mathrm{Somato}}[n]\right)^2\right]},
\end{equation}
where $\Delta f_{\mathrm{Thalamic}}$ is the discrete difference curve of the thalamic cross-correlation curve and the ideal cross-correlation curve, and $\Delta f_{\mathrm{Somato}}$ is a similar difference for cross-correlation curves of the cortical track. The cross-correlation error values are presented in Table \ref{tb:L2error}.

To assess the spatial distribution, the accuracy of localization, and the focality of the estimation, i.e., the degree of concentration around the true source, we visualized the average of the reconstructions calculated from the 20 sample EEG data.

In addition, three alternate experiments were conducted to verify the validity of the experiments. As the first alteration, we inverted the activation order of the thalamic and somatosensory sources. As the second alteration, we removed the thalamic source, and only the somatosensory source was placed. As the third alteration, we changed the position of the somatosensory source and ROI to the visual cortex located at the back of the brain (third and fourth rows of figure \ref{fig:sourcesetup}).

\begin{figure*}[h!]
\centering
\begin{footnotesize}
\begin{minipage}{0.5cm}
\mbox{}
\end{minipage}
\begin{minipage}{4.5cm}
\centering
30 dB 
\end{minipage} 
\begin{minipage}{4.5cm}
\centering
20 dB 
\end{minipage} 
\begin{minipage}{4.5cm}
\centering
10 dB 
\end{minipage} \\ \vskip0.2cm
\begin{minipage}{0.5cm}
\rotatebox{90}{Inverted}
\end{minipage}
\begin{minipage}{4.5cm}
\centering
\includegraphics[trim={0 3cm 0 0},clip, width=4.2cm]{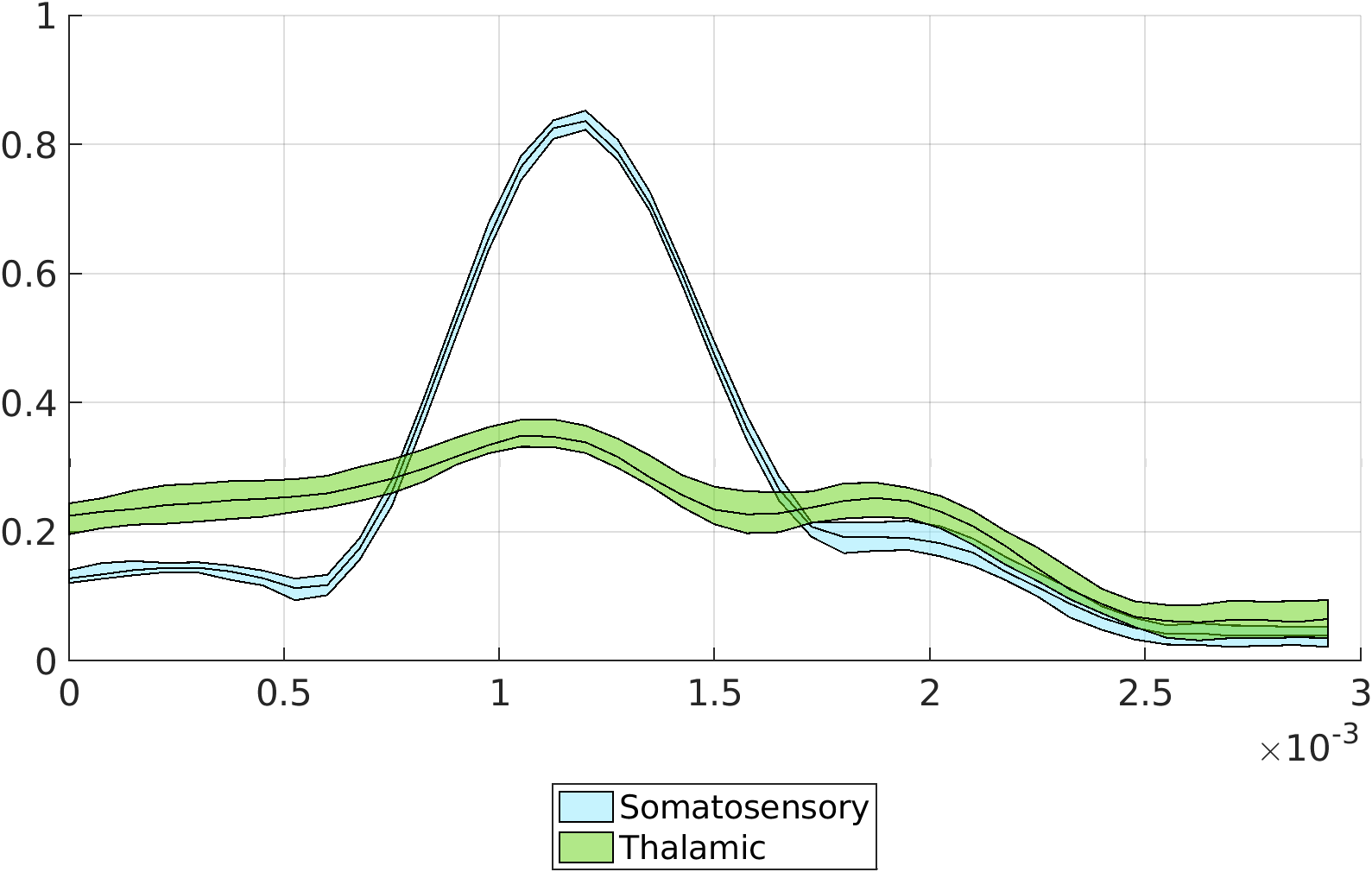}
\end{minipage}
\begin{minipage}{4.5cm}
\centering
\includegraphics[trim={0 3cm 0 0},clip, width=4.2cm]{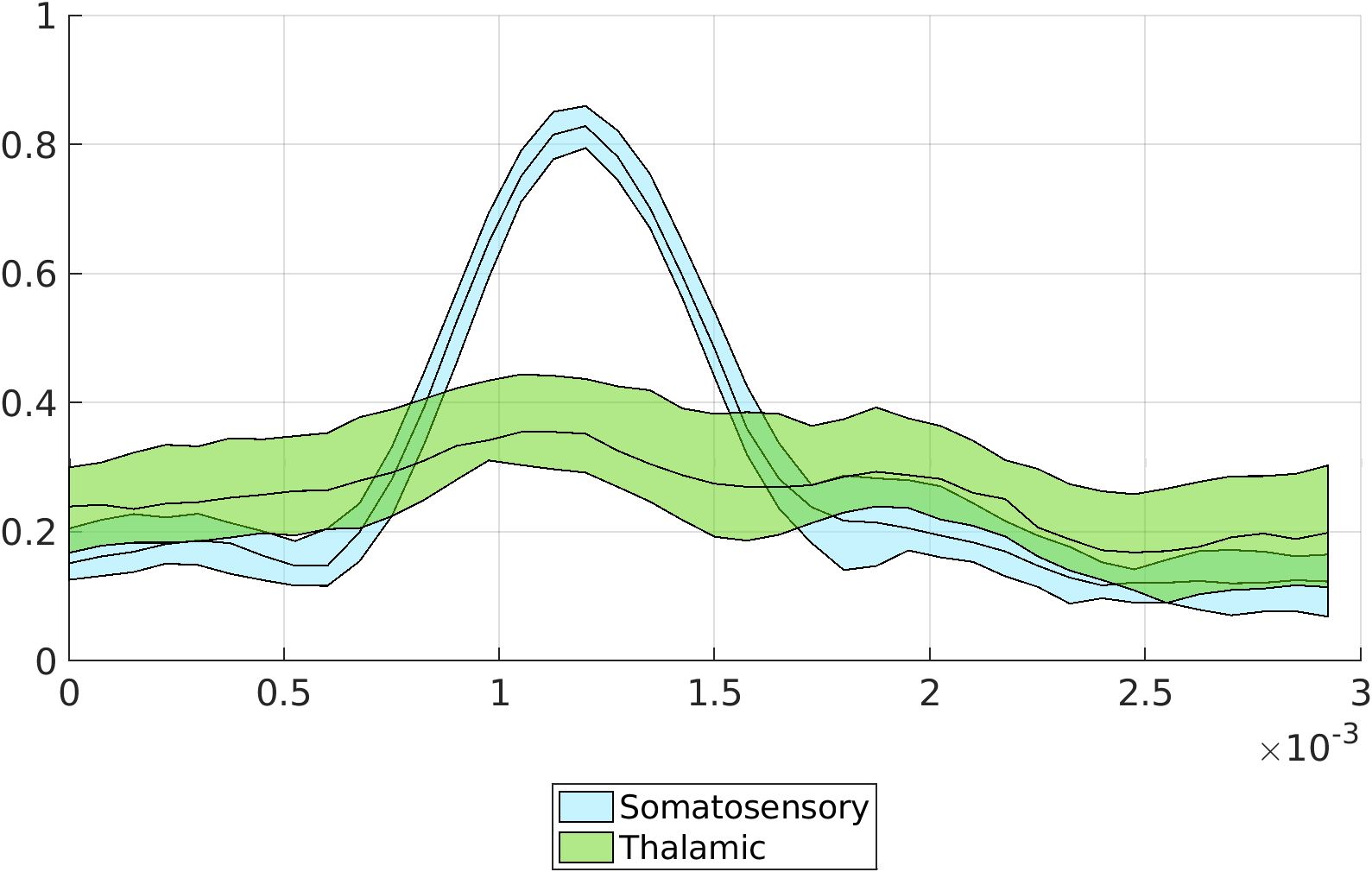}
\end{minipage}     
\begin{minipage}{4.5cm}
\centering
\includegraphics[trim={0 3cm 0 0},clip, width=4.2cm]{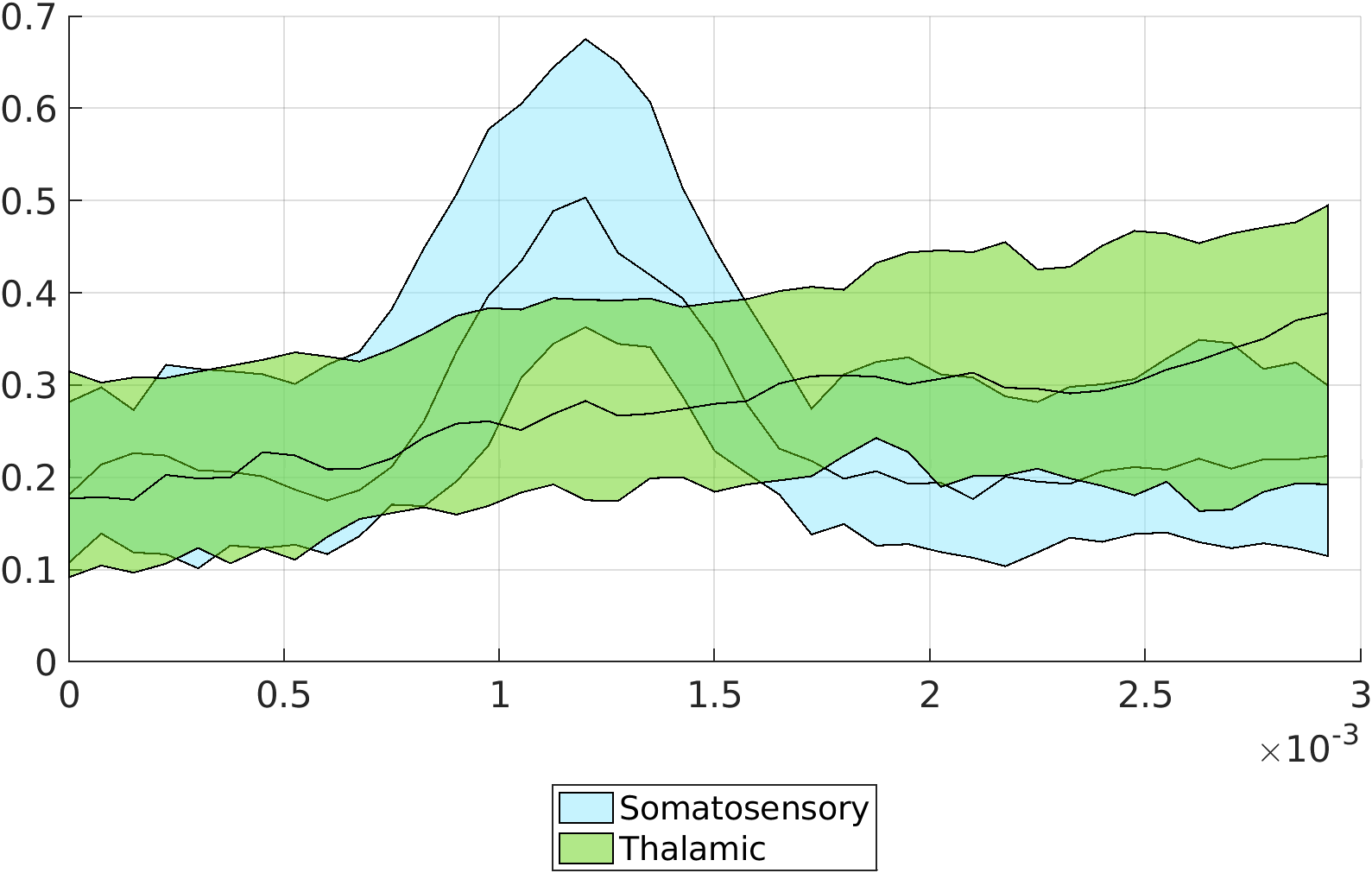}
\end{minipage}   \\ \vskip0.2cm
\begin{minipage}{0.5cm}
\rotatebox{90}{Single source}
\end{minipage}
\begin{minipage}{4.5cm}
\centering
\includegraphics[trim={0 3cm 0 0},clip, width=4.2cm]{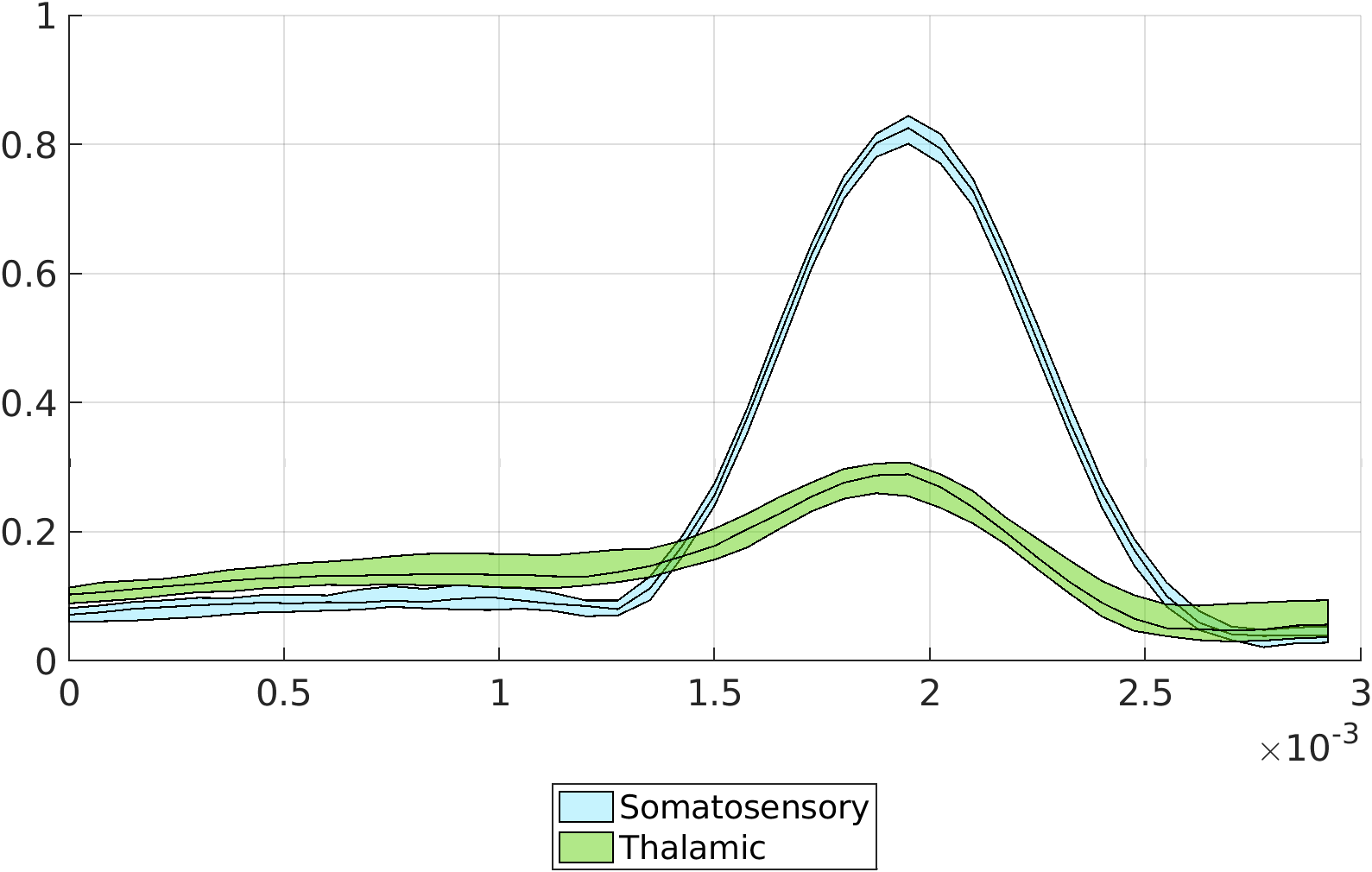}
\end{minipage}
\begin{minipage}{4.5cm}
\centering
\includegraphics[trim={0 3cm 0 0},clip, width=4.2cm]{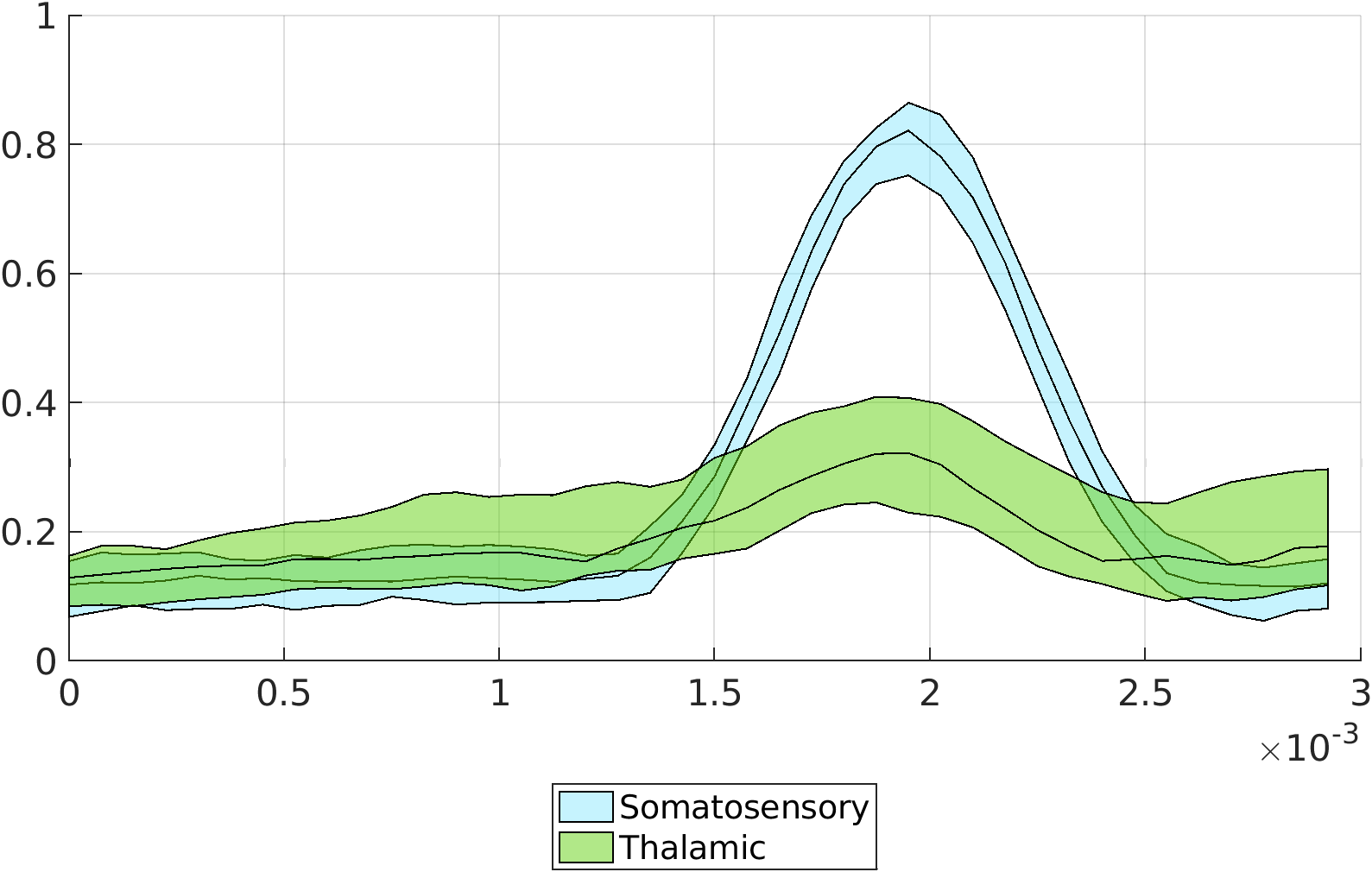}
\end{minipage}     
\begin{minipage}{4.5cm}
\centering
\includegraphics[trim={0 3cm 0 0},clip, width=4.2cm]{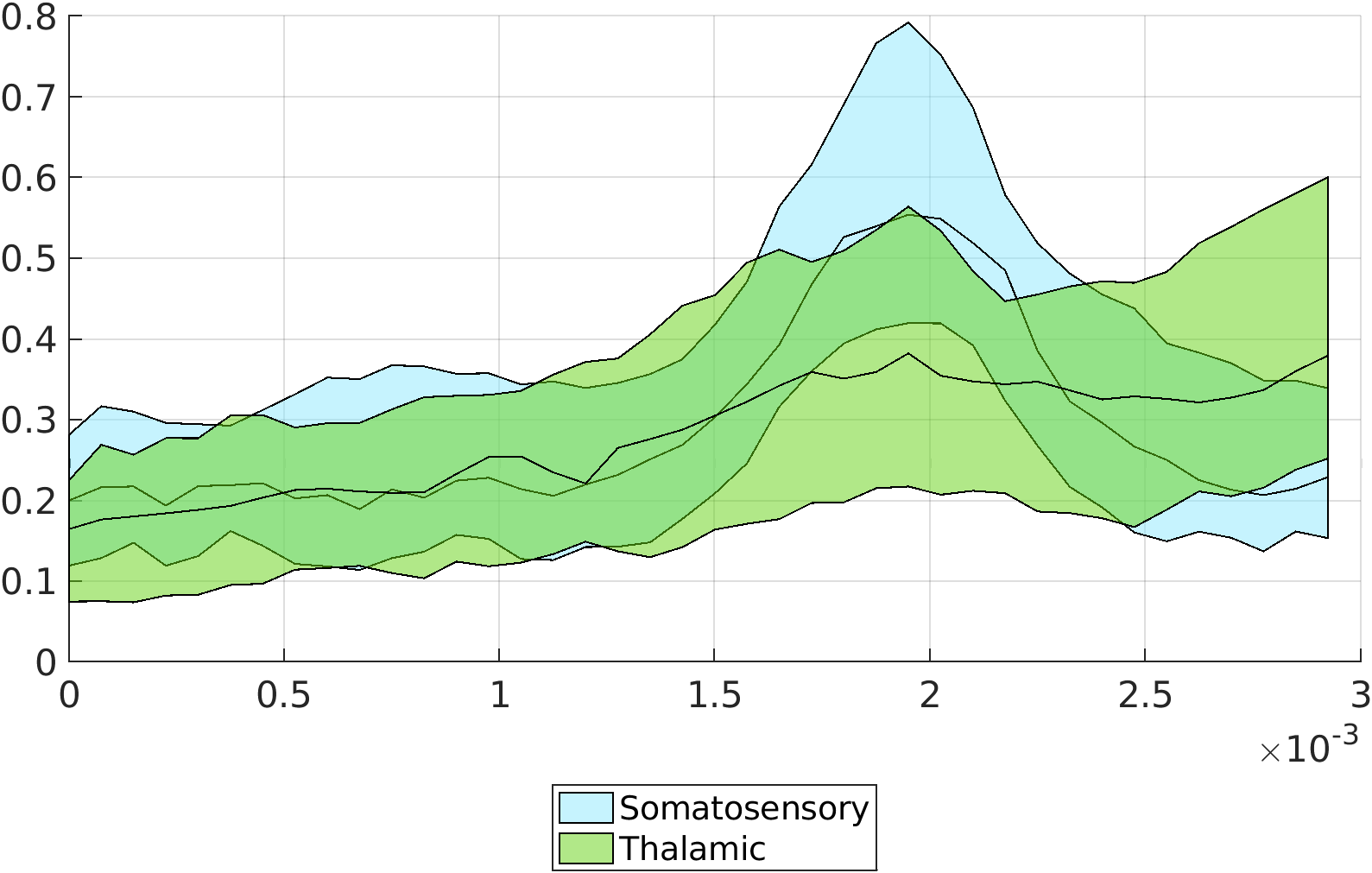}\\
\end{minipage}
\begin{minipage}{12cm}
\centering
\includegraphics[trim={0 0 0 15cm},clip, width=12cm]{VeikkaTimeseries/time_series_nl_10_evp_10_kf_1_sl_1_rts_0.png}
\end{minipage}\\ \vskip0.2cm
\begin{minipage}{0.5cm}
\rotatebox{90}{Visual}
\end{minipage}
\begin{minipage}{4.5cm}
\centering
\includegraphics[trim={0 3cm 0 0},clip, width=4.2cm]{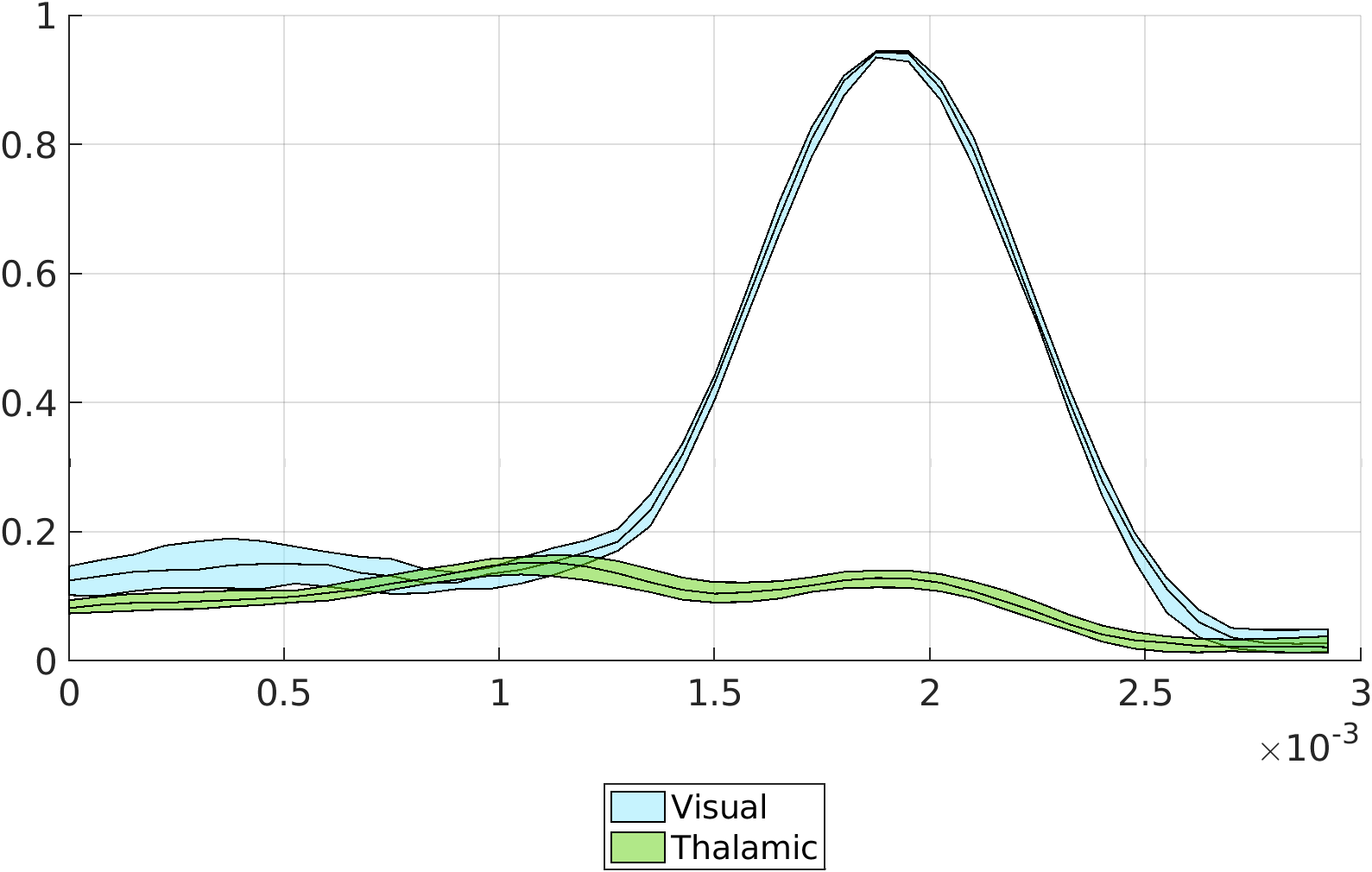}
\end{minipage}
\begin{minipage}{4.5cm}
\centering
\includegraphics[trim={0 3cm 0 0},clip, width=4.2cm]{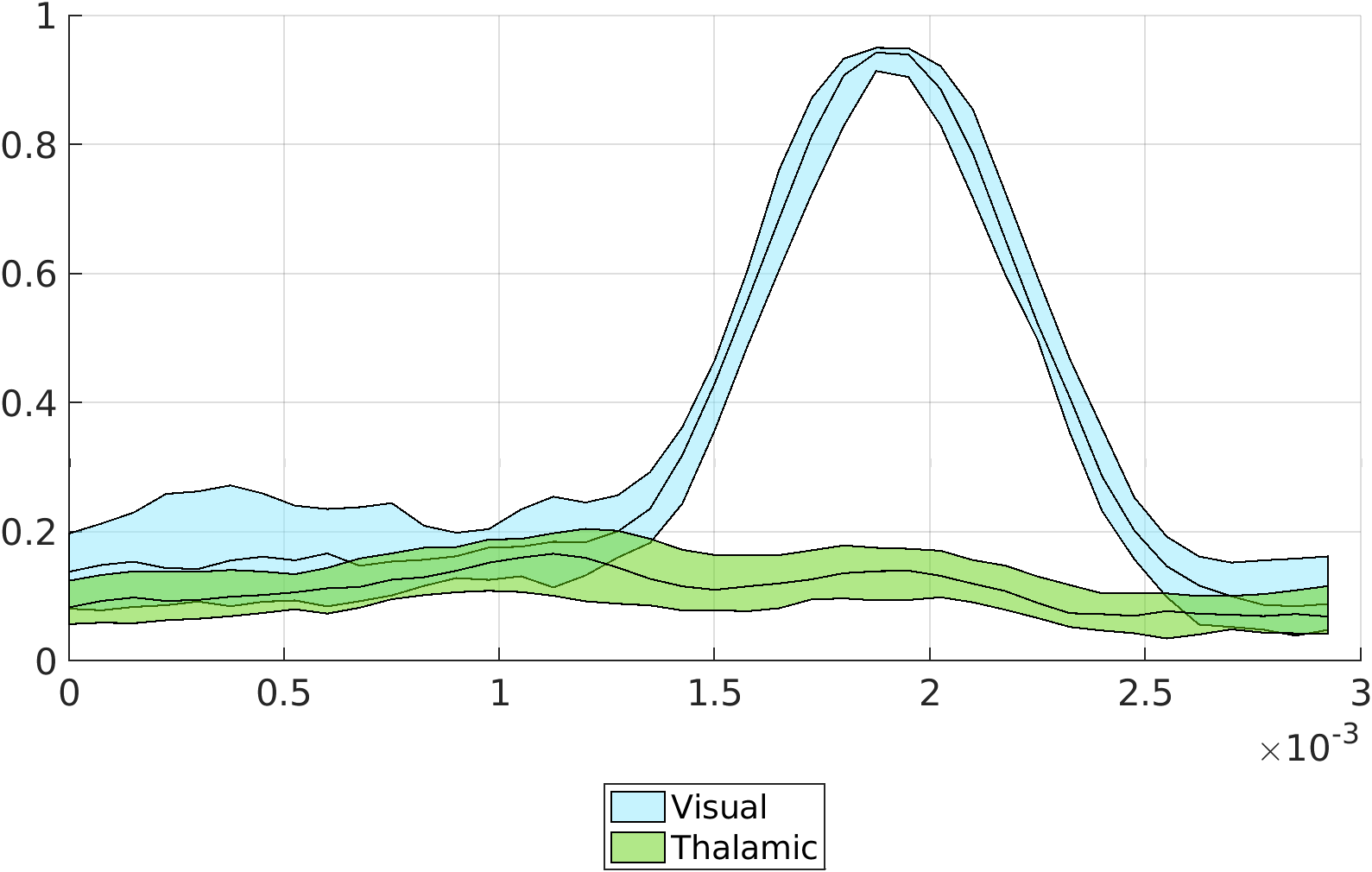}
\end{minipage}     
\begin{minipage}{4.5cm}
\centering
\includegraphics[trim={0 3cm 0 0},clip, width=4.2cm]{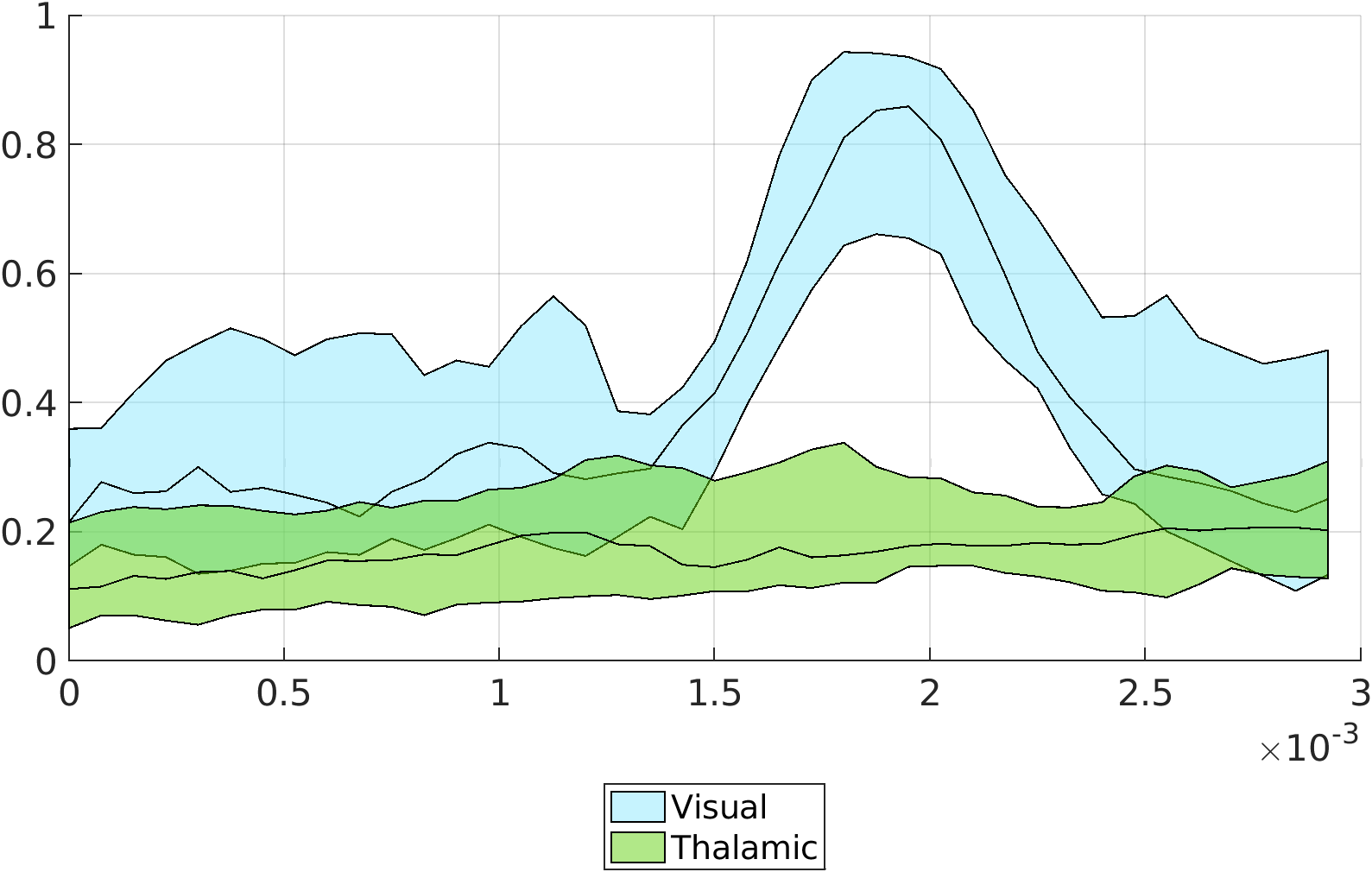}
\end{minipage}
\begin{minipage}{12cm}
\centering
\includegraphics[trim={0 0 0 15cm},clip, width=12cm]{VeikkaTimeseries/inversion_data_kalman_prior_nl_10_evp_30_pmsnr_0_ft_6_st_3_bi_20_visual.png}
\end{minipage}
\end{footnotesize}
    \caption{ Results of the three alternate experiments. Brain activity strength evolution in thalamus (green),  somatosensory cortex and visual cortex (blue). The solid, darker line represents the mean over 20 estimations with different measurement noise realizations. Time runs on the x-axis, presented in milliseconds, and the y-axis shows the strength. }
    \label{fig:TempLabel}
\end{figure*}

\section{Results}

Figure \ref{fig:DynamicKFTracks} presents the reconstructed source time courses for the three algorithms presented in different rows named by the abbreviations of the method names. Each column displays the results for 30, 20, and 10 dB noise from left to right. The x-axis shows the time in milliseconds, and the y-axis shows the estimation strength time courses. The curve surrounded by green represents the estimation strength averaged over the region of interest (ROI) around the thalamic source, and the curve with turquoise area represents the averaged strength over the cortical ROI. The ROIs are presented in the bottom row of Figure \ref{fig:sourcesetup}. The curve progressing in the middle of the colored area is the median strength estimation over the reconstructions of the 20 different data realizations. The colored area represents the second and third quartile range.

\begin{table}[h!]
\caption{Cross-correlation curve error against the ideal cross-correlation of the individual Gaussian pulses presented in Figure \ref{fig:idealCrossCorr} presented for all the compared methods and the used noise levels of 30, 20, and 10 \unit{\decibel}.}
\label{tb:L2error}
\setlength{\tabcolsep}{3pt}
    \centering
    \begin{tabular}{|p{30pt}|p{40pt}|p{40pt}|p{40pt}|}
    \hline
    Method & 30 dB & 20 dB & 10 dB \\ \hline
        3-DSKF &  $0.98\pm 0.11$ & $1.04\pm 0.37$ & $1.73\pm 0.84$\\
        2-DSKF & $0.96\pm 0.18$ & $1.23\pm 0.34$ & $2.24\pm 0.86$\\
        SSKF & $1.27\pm0.18$ & $1.46\pm 0.52$ & $2.19\pm 0.82$\\
        SKF & $1.92\pm 0.12$ & $1.91\pm 0.34$ & $2.06\pm 0.52$ \\
        \hline
    \end{tabular}
\end{table}

The 2-DSKF and 3-DSKF algorithms accurately identify the deep-source peak with correct timing. 3-DSKF provides almost equally high amplitudes and preserves the distinct superficial peak. However, the thalamic activity peaks a few frames before the actual activity peak. On the other hand, correct peak timings are obtained with 2-DSKF, except in the case of 10 dB noise, where the thalamic track becomes indistinguishable from the background noise.

In comparison, SSKF overly smooths the deep-source activity, attenuating its peak, whereas SKF captures the deep activity with reasonable intensity but exhibits higher enmeshment of the curves, i.e., the deep activity emerges while only cortical activity should be present. Moreover, the peaking of the cortical activity is delayed by a couple of frames.

At the lowest SNR (10dB), 3-DSKF is the only method that successfully separates the deep and superficial peaks, maintaining distinct temporal profiles. Additionally, the time courses estimated by 3-DSKF are visibly smoother and less noisy than those produced by 2-DSKF, SKF, and SSKF across all noise levels.

Figure \ref{fig:dynamicKFBrains} depicts the spatial distribution of the averaged reconstructed activity in the cortical and thalamic regions,  at the time points of peak activity (shown in Fig. \ref{fig:trueevolution}),  
 for three noise levels (30, 20, and 10 dB) and for each method: Standardized Kalman Filter (SKF), Standardized Kalman Filter with RTS smoother (SSKF), and Dynamical SKFs (2-DSKF and 3-DSKF). In particular, the first three columns (starting from left to right) of Figure \ref{fig:dynamicKFBrains} show the average reconstructions at time point 1.1 ms and the last three columns the reconstructions at 1.9 ms. Based on this, the proposed DSKF algorithms accurately localize both the subcortical (thalamic) and cortical (somatosensory) activity across all noise levels, maintaining focal and appropriate reconstructions. At higher noise levels, a mild echo of the deep activity appears at the time of the cortical peak, with additional spread into the thalamus and caudate nucleus. The SSKF also correctly localizes the deep activity to the posterior of the left thalamus, but produces broader estimates that extend into the temporal lobe across all SNRs. Notably, SSKF avoids false deep echoes at the time of the cortical peak. In contrast, the SKF method mislocalizes the deep activity at the bottom of the brainstem and shows reduced intensity for the cortical peak, which is the most focal of the methods compared. Additionally, SKF exhibits consistent and strong false deep activations (echoes) during the cortical activity peak across all noise levels.

Considering the cross-correlations of the estimated tracks (cortical and subcortical against each other) in Figure \ref{fig:DynamicKFCrossCorr}, ideally, we should see one Gaussian pulse similar to the blue curve in \ref{fig:trueevolution}. We can see that the 3 and 2-DSKF provide almost such curves, except for the tails at both ends for 30 and 20 dB. With SSKF and SKF, we see a significant distortion at the end part, signaling unwanted correlations between both cortical and subcortical pulses. Further examining the cross-correlations of individual thalamic activity (Figure \ref{fig:DynamicKFCrossCorrTrue1}) and cortical activity (Figure \ref{fig:DynamicKFCrossCorrTrue2}) against the true tracks, we should see in both cases half of the Gaussian pulse, as presented in Figure \ref{fig:idealCrossCorr}. For 3-DSKF, the thalamic cross-correlation curve represents the ideal result the most at all noise levels considered. And by eye, the cross-correlation of the somatosensory tracks of 2-DSKF and 3-DSKF are ideal. In Table \ref{tb:L2error}, we can see that if we consider the mean error, computed using (\ref{eq:correrr}), it is lowest for 2-DSKF at 30 dB with a value of 0.96, but otherwise 3-DSKF has the lowest mean errors (0.98, 1.04, and 1.73). However, in light of the 10 and 90 \% error bounds, the bound is distinctly lower, 0.11, as it is 0.18 for 2-DSKF; therefore, we could consider 3-DSKF providing more truthful results in the majority of the sampled scenarios. 

In addition, to further assess the adaptability of the current method, we repeated the simulations with 2-DSKF using (i) reversed activation order, (ii) only a cortical source, and (iii) a cortical source relocated to the visual cortex. We chose 2-DSKF for these trials as it offers a good balance between accuracy and computational efficiency. Across all scenarios considered, 2-DSKF maintained correct temporal sequencing, suppressed false activations in the single-source case, and localized activity accurately to the intended cortical region, confirming that the method generalizes beyond the original setup (Figure \ref{fig:TempLabel}).

\section{Discussion and Conclusions}

This work introduces a Dynamical Standardized Kalman Filter (DSKF) that combines temporal dynamics of elevated order with power-weighted depth standardisation for EEG source localisation. Across a range of signal-to-noise ratios, DSKF consistently achieved more accurate and less enmeshed spatial localisation than the Standardised Kalman Filter (SKF) and its Rauch-Tung-Striebel (RTS) smoothed variant (SSKF), while also producing smoother and more temporally distinct source time courses. In simulations involving concurrent deep and superficial generators, DSKF preserved temporal separation and anatomical specificity even under high noise, while reducing spurious deep activations commonly seen in existing methods. 

By extending the state vector to include not only dipole activity but also its first and second temporal derivatives (velocity and acceleration), DSKF better reflects the inherently dynamic nature of EEG signals, which often involve rapid changes such as event-related potentials or oscillatory bursts. This higher-order state formulation provides two key advantages: it yields smoother and more physiologically plausible source estimates by preventing abrupt, unrealistic changes, and it improves predictive capability by anticipating future dynamics. These properties are particularly valuable under low-SNR conditions or with limited sensor coverage, where traditional models struggle to disambiguate signal from noise.

Previous Kalman-based EEG source estimation methods have relied primarily on zero-order dynamics, which limit temporal precision or require sensitive tuning of process noise \cite{GalkaAndreas2004KalmanEEG, Hamid2021}. While the Standardised Kalman Filter (SKF) addressed depth bias through post-hoc standardization weighting \cite{Lahtinen2024SKF}, its temporal evolution model remained restricted to a simple random-walk model. 

The results show improvements in recovering the correct temporal evolution of simultaneous subcortical (at the thalamus) and cortical activity (at the somatosensory cortex) through kinematic evolution modeling. However, the challenge of obtaining completely independent cortical activity without a false echo at subcortical locations remains. The reason could be that the observation model produces highly similar electrical potentials on the surface of the model for both sources, causing a weak dissolution for the sources. 

Assuming that we cannot eliminate the false activity, we can compute the height difference between the curves at the peak point, as shown in the following table for the results of this paper.

\begin{table}[h!]
\caption{Track height difference for the cortical peaking moment 1.1 \unit{\milli\second} and the thalamic peak moment at 1.9 \unit{\milli\second} presented for all the compared methods and the used noise levels of 30, 20, and 10 \unit{\decibel}. The results for 10 and 90 \% quantiles are used to infer the error bounds.}
\label{table}
\setlength{\tabcolsep}{3pt}
    \centering
    \begin{tabular}{|p{30pt}|p{25pt}|p{25pt}|p{25pt}|}
    \hline
    \multicolumn{4}{|c|}{\bf Deep peak}\\ \hline
    Method & 30 dB & 20 dB & 10 dB \\ \hline
        3-DSKF &  0.18 & 0.18 & 0.16\\
        2-DSKF & 0.16 & 0.14 & 0.02\\
        SSKF & 0.02 & 0.01 & -0.05\\
        SKF & 0.15 & 0.13 & 0.03 \\
        \hline
    \multicolumn{4}{|c|}{\bf Cortical peak}\\ \hline
    3-DSKF &  0.13 & 0.14 & 0.11\\
        2-DSKF & 0.54 & 0.53 & 0.2\\
        SSKF & 0.64 & 0.63 & 0.41\\
        SKF & 0.16 & 0.15 & 0.06 \\
        \hline
    \end{tabular}
\end{table}

We can obtain a balance between the peak height differences for the Standardized Kalman filter utilizing the first-order kinematic model (2-DSKF). This could be interpreted as a secondary approach to evaluate whether the measurement data contains the firing of multiple separate sources, either subsequently or simultaneously. By considering the 3-DSKF method, which was otherwise best by the quantitative measures, the estimated evolution curve allows us to easily believe that the deep activity evolves slowly and gradually vanishes in two stages. Although the original Standardized Kalman filter algorithm and its smoothed counterpart clearly indicate two separate subcortical activities to emerge. Furthermore, the alternate experiments carried out using the 2-DSKF method demonstrate that the proposed approach generalizes well to varied source configurations.

Future work will focus on validating DSKF using real EEG data, including evoked potentials with known timing and location, as well as multimodal datasets such as EEG-fMRI or intracranial EEG. Extensions to adaptive or switching-dynamics models may allow the method to better accommodate non-stationary neural activity. In addition, data-driven or Bayesian strategies for selecting the optimal power-weighting and process noise parameters could further enhance robustness and reliability across subjects and recording environments.


\appendices

\section{Rationale for the Standardization Exponent $p$}
\label{app:exponent}

The Standardized Kalman Filter (SKF) \cite{Lahtinen2024SKF}, introduced a time dependant post hoc weighting approach on the Kalman filtering estimates with a dynamical standardization matrix, $W_t = \mathrm{Diag}\left(P_{t\mid t-1}^{-1/2} K_t S_t K_t^\mathrm{T} P_{t\mid t-1}^{-1/2} \right)^{-p} P_{t\mid t-1}^{-1/2}$ and the standardized Kalman estimate at time step t, ${\bf z}_{t\mid t}$, is obtained by applying this weighting on the posterior mean ${\bf x}_{t\mid t}$. In the original work, the standardization exponent $p$ is 0.5 as the standardization is applied to the amplitude of the estimate, likewise in \cite{PascualMarqui2002}. In this study, instead, we apply $p = 1$ (standardization of estimation power).  This is to avoid suppressed amplitudes which are likely to occur, especially, when Rauch--Tung--Striebel (RTS) smoother is applied after filtering, as the smoother itself is not standardized. A greater standardization exponent can be interpreted to compensate for the bidirectional nature of the smoother, combining past and future information, thereby restoring the relative amplitude of the smoothed estimates. The choice $p = 1$, i.e., power standardization, is thus a necessary extension of the standard $p = 0.5$. In this study, we use it systematically in all of our numerical experiments.





\bibliographystyle{IEEEtran}
\bibliography{references}

\end{document}